\documentclass[11pt,letterpaper]{amsart}
\usepackage{amsfonts}
\usepackage{amssymb}
\usepackage{amsthm}
\usepackage{amsmath}
\usepackage{amscd}
\usepackage[latin2]{inputenc}
\usepackage{t1enc}
\usepackage[mathscr]{eucal}
\usepackage{indentfirst}
\usepackage{graphicx}
\usepackage{bbm}
\usepackage{graphics}
\usepackage{pict2e}
\usepackage{epic}
\usepackage[capitalize]{cleveref}   %for \cref

\makeatletter
\newcommand{\addresseshere}{%
  \enddoc@text\let\enddoc@text\relax
}
\makeatother
\DeclareMathOperator*{\dom}{Dom}

\numberwithin{equation}{section}
\usepackage[margin=1in]{geometry}

\usepackage[
backend=biber,
natbib,
style=authortitle,
citestyle=authoryear
]{biblatex}
\usepackage{epstopdf}
\theoremstyle{plain}
\newtheorem{Th}{Theorem}[section]
\newtheorem{lemma}[Th]{Lemma}

\newtheorem{prop}[Th]{Proposition}

\newtheorem{Def}[Th]{Definition}

\newtheorem{?}[Th]{problem}

\newcommand{\eq}[1]{\begin{align}#1\end{align}}

\newcommand{\Z}{\mathbb{Z}}
\newcommand{\RRR}{\mathbb{R}}

\newcommand{\PP}{{\underline{P}}}
\newcommand{\PPP}{{F^{n}}}

\makeatletter
\renewcommand\subsection{\@startsection{subsection}{2}%
  \z@{.5\linespacing\@plus.7\linespacing}{-.5em}%
  {\normalfont\scshape}}

\def\specialsection{\@startsection{section}{1}%
  \z@{\linespacing\@plus\linespacing}{.5\linespacing}%
  {\normalfont}}% NEW
\def\section{\@startsection{section}{1}%
  \z@{.7\linespacing\@plus\linespacing}{.5\linespacing}%
  {\normalfont\scshape}}% NEW
\makeatother

\usepackage{xcolor}% for colors

\title{A Stochastic Model of Mathematics and Science}

\author{David H. Wolpert}
\address{Santa Fe Institute, Santa Fe, New Mexico \\
Complexity Science Hub, Vienna\\
Arizona State University, Tempe, Arizona\\
\texttt{http://davidwolpert.weebly.com}}
\author{David B. Kinney}
\address{Princeton University, Princeton, New Jersey\\
\texttt{http://davidbkinney.com}}

\begin{document}
\begin{abstract}
We introduce a framework that can be used to model both mathematics and  human reasoning about mathematics.
This framework involves \textit{stochastic mathematical systems} (SMSs), which
are stochastic processes that generate pairs of questions and associated answers (with no explicit referents). We
use the SMS framework  to define normative conditions for mathematical reasoning,
by defining a ``calibration'' relation between a pair of SMSs. The first SMS is the human reasoner, and
the second is
an ``oracle'' SMS that can be interpreted as deciding whether the question-answer pairs of the reasoner SMS are valid.
To ground thinking, we understand the answers to questions given by this oracle to be the answers that would be given by an SMS representing the entire mathematical community in the infinite long run of the process of asking and answering questions.
We then introduce a slight extension of SMSs to allow us to model both the physical universe
and human reasoning about the physical universe. We then
define a slightly different calibration relation appropriate for the case of scientific reasoning.
In this case the first SMS represents a human scientist
predicting the outcome of future experiments, while the second SMS
represents the physical universe in which the scientist is embedded, with the question-answer pairs of that SMS
being specifications of the experiments that will occur and the outcome of those experiments, respectively.
Next we derive conditions justifying two important patterns of inference in both mathematical and scientific reasoning: i) the practice of increasing one's degree of belief in a claim as one observes increasingly many lines of evidence for that claim, and ii) abduction, the practice of inferring a claim's probability of
being correct from its explanatory power with respect to some other claim that is already taken to hold for independent reasons.
\end{abstract}

\maketitle

\section{Introduction}
%This paper concerns two questions that are both
%central to the philosophy of science, the philosophy of mathematics, and philosophy more broadly. The first question is,
%%a question about reasoning:\
%how should we represent the formal reasoning in mathematics and science? The second question is,
%%concerns representation and ontology:\
%how should we model the physical universe?
%
%Here, we introduce a formalism that can resolve both questions. Specifically, we propose to

%\dhwc{Go through introduction, making sure its emphasis parallels that in the abstract.}

In this paper, we introduce a framework that extends both conventional formulations of
mathematics and conventional formulations of the structure of the physical universe. Because of how
it extends those formulations, this new framework can
also be used to represent (and investigate) human reasoning about mathematics, i.e., the reasoning of individual
mathematicians, or the collective reasoning of communities of mathematicians. Similarly,
our new framework can also be  used to represent (and investigate) human reasoning about the structure of
the physical universe, i.e., the reasoning of individual
scientists, or the collective reasoning of communities of scientists.

The central motivation behind our new framework is simple: One of the defining features of mathematicians
%(or communities thereof)
and of scientists
%(or communities thereof) i\
is that at any given moment, they have
what they perceive (rightly or wrongly) to be a set of ``(almost) incontrovertibly true'' propositions. (Examples of such propositions for mathematicians
include accepted proofs, and examples for scientists include
observations or the results of past experiments.) They then map those
``true'' propositions  to beliefs concerning other
propositions that they are uncertain about.
That is, they use a corpus of settled facts to formulate beliefs about propositions that are not settled facts.
Our goal is to formalize this map
% that lies at the heart of both mathematical and scientific practice,
and start
to investigate its formal properties. We are particularly interested in presenting a framework
that captures how the inputs to that map change over time as the mathematician / scientist conducts
their research, i.e., how their set of established facts and their set of beliefs about non-settled
facts evolve over time. In creating this framework this we make the following contributions:
\begin{enumerate}
    \item We introduce a very broad formalization of  the set of accepted facts as well as
the beliefs
of mathematicians or scientists concerning propositions that they do not perceive as accepted fact.
There is a long tradition within Bayesian epistemology of using probability theory
%~\cite{kolmogorov1933grundbegriffe},
to represent the degrees of belief that epistemic agents (including mathematicians and scientists) assign to propositions that are not in their corpus of accepted facts.
%`Bayesian epistemology' in this context is the broad family of approaches to epistemological modeling in which an agent's degree of belief in some proposition is represented by a probability assigned to an element of some event algebra.
Here, we embrace this probabilistic approach to the epistemological modeling of mathematicians and scientists. However,
we do not want to limit ourselves by assuming that the settled facts and propositions of math
and science are all represented in some single particular formal system --- an assumption which would be patently
false for any choice of such a formal system, since real mathematicians and scientists formulate the propositions
they are investigating using a broad range of formal systems.
To avoid that assumption, we formulate both settled facts and unsettled
propositions as question-answer pairs, without restricting the formal system used to represent either
questions or answers.

\item In the real world, at any particular moment,
any given mathematician or scientist has a set of \textit{multiple} accepted facts,
i.e., multiple question-answer pairs, and a set of \textit{multiple} non-settled facts that they
have beliefs about. Accordingly, we represent agents' epistemic attitudes as conditional probability distributions
over possible answers, conditioned on specific associated questions and a set of possibly multiple (established facts represented as)
independent question-answer pairs, where
both the questions and answers are elements of some arbitrary (and
possibly infinite) associated space. For every question-answer pair that the mathematician / scientist
considers to be accepted fact, the conditional probability of that answer given that question equals
one. For all other pairs --- all other question-answer pairs that they have beliefs about, but do not
view as accepted facts --- that conditional probability is less than one.

\item To model how such distributions over sets of question-answer pairs evolve over time, we need to
capture the fact that the reasoning of both mathematicians and scientists is stochastic, both
in the questions that they choose to investigate and the answers they provide to those questions.
%Reflecting this, the form of Bayesian epistemology that we introduce
%is based on what we call a \textit{stochastic mathematical system} (SMS). This
Stated formally, such evolution is
a stochastic process, which generates successive sets of question-answer pairs.
%in which questions are iteratively posed and answered.
%Thus, we represent both mathematical and scientific reasoning, and the physical universe itself, as a stochastic process of asking and answering questions.
%In what follows, we define this approach in full mathematical detail, and explore some of its elementary implications.\par
%
%Our approach allows us to explore new territory in the space of conceptual approaches to reasoning, the physical universe, and the relationship between the two.
%Regarding mathematical reasoning, our approach allows us to express the thesis that mathematics itself may be a non-trivially stochastic process; that is, u
%The reasoning processes of both real human mathematicians and real human scientists are stochastic,
%not only in the questions that real mathematicians and scientists pose, but also in the answers that
%they generate to those questions.
We refer to such a stochastic process as a \textit{stochastic mathematical system} (SMS).
SMSs are the central element of our approach.

%Indeed, one of the important ways we break from the tradition within Bayesian epistemology
%is by allowing the answers in such pairs to themselves be distributions, over answers to yet other questions.
%In this way we capture the fact
%%\textit{objects} degree-of-belief attitude (i.e., we add structure to the
%%formal representation of the objects of degrees of belief.
%%We hold
%that for many purposes, beliefs about unsettled facts in mathematics and science can
%be formalized \textit{distributions over answers to questions}.
%, if we allow
%the space of both questions and answers to be broad enough.\par

% \dhwc{Perhaps the possibility of answers being predictions should only be introduced below? After all,
% we don't require it (which the preceding paragraph implies). And
% it's secondary to a lot of the other points in the introduction.}

\item Our interrogative approach to the probabilistic modeling of reasoning in mathematics and science allows us to capture several
incontrovertible aspects of said reasoning in the real world. (Note that we make no claims about how that
reasoning itself \textit{should} take place.) First, it is a simple fact that for real-world mathematicians (or communities thereof) and real-world scientists (or communities thereof), if the probability distribution over answers to some questions is to be representationally accurate, then it cannot be a delta function. In both mathematics and science, there are questions that amount to \textit{open problems}, such that even the most advanced practitioners are unsure of the correct answers.
If you ask a computer scientist (for example), ``Does $\mathrm{P} = \mathrm{NP}$?'' (one of the most famous open questions in all
of mathematics and science), many will refuse to give a definitive answer. If instead you ask them,
``What is the \textit{probability} that $\mathrm{P} = \mathrm{NP}$?'', now almost all will be happy to provide a response.
This is captured in the SMS framework by allowing answers to questions to specify probability distributions
over answers to other questions. Thus, the SMS framework allows us to model both the essentially \textit{interrogative} nature of mathematical and scientific reasoning and the inherent stochasticity of real human reasoners.

\item We use
the SMS framework to address the question of how to formalize the foundational nature of mathematics \textit{itself},
as opposed to formalizing only the nature of mathematical reasoning. Our conceptual starting point is to take seriously the slogan that ``mathematics is what mathematicians do.'' This leads us to use the SMS framework to formulate
mathematics itself. Using the SMS framework this way allows
for the possibility that there is no unique answer to a given mathematical question,
but instead a non-degenerate objective probability distribution over possible answers.
%(\textit{Mutatis mutandi} for
%questions concerning events in the physical universe.)
It allows for the possibility that \textit{mathematics itself} may
inherently and irrevocably contain contradictions. This is in contrast to conventional formulations of mathematics
in terms of formal systems, category theory, ZFC set theory, etc., in which one set of
hypotheses either does or does not imply another set, with no role for stochasticity,
and according to which ``mathematical existence is freedom from contradiction'' (\cite{hilbert1928grundlagen}). In this respect, we build on work by \citet{bueno2011inferential} and \citet{mccullough2020representing}; they argue that structuralist approaches in philosophy of mathematics (e.g., Shapiro \citeyear{shapiro1997philosophy}) are limited by their inability to represent the possibility that mathematics is inherently contradictory. However, our framework departs radically from previous work in the foundations of mathematics by holding that not only can mathematics be fundamentally contradictory, it can also be fundamentally \textit{stochastic}, such that the answers to mathematical questions are fundamentally generated by sampling from a non-delta-function probability distribution. Thus, we use our framework to put forward a novel position in the foundations of mathematics.

\item An important normative motivation for mathematicians is almost always that their probability distributions over answers to possible questions be very close to those that would arise if they had access to some oracle. To ground the reader's thinking, one can follow
\citet{Peirce1878-PEIHTM}, and take the oracle to be the far-future community of all mathematicians,
located anywhere in the universe. Many mathematicians would feel they have made a ``good'' prediction for the answer to a question if it is the same one that would be given by a far-future community. In what follows, we show how the SMS framework is able to formalize this crucial normative feature of mathematical reasoning
by considering appropriate expected values of the divergence between the probability distribution of a mathematician SMS
and that of mathematics itself. To do this, we first suppose
that some  ``oracle'' SMS stipulates what the answers are to arbitrary questions.
We then further suppose that the normative goal of any mathematician is for the process that they use to answer
to arbitrary mathematical questions to correspond to the greatest degree possible with the process used by that oracle.
Specifically, we first postulate that both an individual mathematical reasoner and such an oracle
can be represented as a separate SMS. We then define that mathematician-SMS to be
``calibrated'' with the SMS representing the oracle if the expected accuracy of the answer that the
mathematician gives to a particular question agrees with that of
the oracle, up to a given threshold.
%One would expect that community not to generate question-answer pairs
%that are broadly considered to be open issues, or that are considered
%only likely to be true, rather than considered to be unequivocally established.
% The goal of a present-day mathematician
% would be to provide answers in a way that coheres to the greatest degree possible with the answers of that far-future community
% of mathematicians. We will in fact adopt this interpretation of the ``oracle'' in our
% discussion below, to try to hone intuition. However, it is important to note
% though that the SMS approach does not require that interpretation in any way; though our formalization of the normative constraints on mathematicians is precise, its interpretation is highly general.
%In this way, we are able to give a precise statement of the sense in which a given mathematician ought to answer a particular question in a particular way.\par

\item Moving to the case of empirical science and the nature of the physical universe, our approach is similar to the one that we use to model mathematical reasoning and the ``mathematical universe'', but with some important differences. As is the case with mathematicians and mathematics,
we use the SMS framework to formulate
humans reasoning stochastically about the physical universe.  Similarly, we can also use
the SMS framework to formulate stochastic models of the physical universe itself. At one level,
this allows us to capture simple truth, that physical facts about the future, governed the
present, are random. As an example, \textit{conditioned on the data a given climate scientist has}
(their set of ``established facts''), the precise global temperature five years from now is governed by
a probability distribution.\footnote{Note that in this regard, formulating the physical universe
as having inherent stochasticity might make more sense than formulating mathematics that way.} But modeling the physical universe using the SMS framework allows us to go further than this, to
%%we can use the SMS framework
%%not only to represent humans reasoning about the physical universe, to formulate stochastic models of the physical universe.
%These stochastic models of the physical universe include Copenhagen-style stochasticity in the events generated by the
%``laws of physics'' as a special case. However, they also include
%stochasticity in the ``laws of physics'' themselves, as another special case.
entertain the ontological thesis that the physical universe ``is'' an SMS, as opposed to the epistemological theses
just predictions about some of physical variables conditioned on knowledge of other physical variables are
necessarily uncertain.
In this respect, our approach has affinities with ontic structural realist approaches in philosophy of science (e.g., Ladyman and French \citeyear{Ladyman1998-LADWIS-2}; Ladyman \citeyear{ladyman2007scientific}), as well as arguments from theoretical physics that the universe is inherently mathematical (e.g., Barrow \citeyear{barrow1991theories}, \citeyear{barrow2011godel}; Schmidhuber \citeyear{schmidhuber1997computer}; Tegmark \citeyear{tegmark1998theory}, \citeyear{tegmark2008mathematical}, \citeyear{tegmark2014our}). However, the SMS framework allows us to represent different ways in which nature might be fundamentally stochastic with a high degree of generality.
For example (and as we will show in what follows), our framework allows us represent the possibility that the {events} of the universe are stochastic (as in some interpretations of
quantum mechanics, for example). However, it also allow us to investigate the implications if the \textit{laws} of
the physical universe themselves are stochastic. This includes as a special case the possibility, widely
entertained in modern cosmology, that the physical constants in the laws of physics vary randomly
across the universe.
%%, not just the \textit{events} of the universe (as in some interpretations of
%%quantum mechanics, for example).
%% for the possibility that physical universe is inherently stochastic.\par
%%As for the relationship between scientific reasoning and the physical universe,
%In this way, the SMS approach
%allows for the possibility that there is a non-degenerate objective probability distribution over
%the outcomes of physical experiments.

\item %by modelling the physical universe as an SMS, we invite the realist inference to the ontological claim that the physical universe \textit{is} an SMS. While
As in the mathematical case, we use the SMS framework to
formalize the normative aims of physical reasoning, specifically the idea that the answers given by some scientists to some questions are ``physically correct.''
To do this we first suppose that the human scientist and the physical universe are SMSs.
An extra complication in this case, not present in the case of mathematics, is that such scientists are themselves
physical systems. So formally, they are answers to certain questions posed to the very physical universe-SMS about
which those self-same human scientists are reasoning.
To capture this we define a special ``embedding'' relationship
that might hold between the answers of the scientist-SMS and those of the embedding universe-SMS.
(This can also be viewed as a way to provide a ``physical meaning'' to the answers generated by the scientist-SMS.)
%physical universe of which they themselves are a part, and that this relationship allows for the accurate prediction of phenomena in at least some domains of inquiry. Specifically, we define an embedding relation between the SMS representing a particular scientist's reasoning and the SMS representing the physical universe, which codifies the relationship between a scientist and the universe in which they live.
We then define a scientist-SMS to be ``embed-calibrated'' with the universe-SMS in which they are
embedded if the expected accuracy of the scientist's answers to some question(s) is the same as that of the embedding physical universe, up to some threshold. This allows us to represent accurate physical scientists in much the same way that we represent accurate mathematicians; the embedding physical universe plays an analogous
role for scientists as the far-future community of mathematicians plays for mathematicians. (As elaborated
below, ultimately the difference between the definitions of mathematical accuracy and scientific accuracy
is due to the fact that scientists
are embedded in the SMS about which they are asking questions, whereas mathematicians ask questions about a
wholly independent SMS.)

\item Formalizing these normative constraints on mathematical and scientific reasoning leads to
one of our primary goals in this paper, which is to
provide conditions that vindicate the use
%As briefly mentioned above, the SMS framework
%also allows us to present conditions under which
two types of heuristic reasoning commonly used by real world, human researchers. Both of these heuristics have resisted straightforward vindication within the traditional frameworks of Bayesian epistemology:
\begin{enumerate}
    \item The first heuristic involves assigning greater degree of belief to a hypothesis when that hypothesis is supported by multiple lines of evidence. That this heuristic is normatively justified is sometimes called the ``variety of evidence thesis.'' Historically (and particularly in light of results due to \cite{Bovens2002-BOVBNA}), this heuristic has been viewed as incompatible with Bayesian epistemology, at least if the heuristic is to be vindicated as a general maxim of epistemic rationality. Indeed, \citet{Landes2020-LANVOE} declares that ``the Bayesian quest for a vindication of the [value of evidence thesis] in full generality has failed'' (p.\ 185). Such claims are made in the context of a traditional Bayesian epistemology. By contrast, we show that on our interrogative approach to modeling both descriptive and normative aspects of reasoning, a highly general version of the variety of evidence thesis holds for both scientific and mathematical reasoning.

    \item The second heuristic is abduction, i.e., the inference from the explanatory power of a claim to its probability of being correct. Here too, there is a widespread perception that this heuristic is not fully consistent with Bayesian epistemology. Indeed, in his recent book-length treatment of abduction, \citet{douven2022art} states that his goal is to convince the reader ``that abduction is interestingly different from Bayes' rule, and that there is nothing intrinsically untoward about it'' (p.\ 25). By `Bayes' rule' here, Douven implicitly invokes a version of Bayesian epistemology that lacks the interrogative structure of the SMS framework. In what follows, we show that within the SMS framework, one can provide a highly general vindication of abduction that is wholly consistent with the use of probability theory, and only probability theory, to represent the degree to which agents believe a given proposition. Moreover, this defense of abduction does not require that an agent adopt a specific, putatively correct prior as some Bayesian approaches (e.g., \cite{weisberg2009locating}) do. Instead, we impose on the weaker constraint that the agent be calibrated with some arbiter of normatively justified reasoning. As discussed above, this crucial notion of calibration is defined in formal detail below.
\end{enumerate}

\end{enumerate}

%  In fact, we would argue that any particular approach to modeling
% real-world mathematical and scientific reasoning must be consistent with the SMS framework, if it
% involves Bayesian epistemology. It may
% be a special case of that framework, e.g., under some particular choice for the stochastic process. But to match the real-world
% behavior of human reasoners, it cannot contradict that framework. \dhwc{That's a very strong ending statement I
% point in --- what do you think?}
\par
\noindent
We take these contributions to demonstrate the fecundity of our approach within philosophy of science and formal epistemology.

\subsection{Roadmap}

In \cref{sec:sms}, we first provide a minimal description of SMSs, which will serve as the basis
for the analysis in this paper.
(See Appendix~\ref{sec:appendixfullframework}
for a more complete, fully formal set of definitions, which may be useful in other contexts.)
Next, in \cref{sec:distributions}, we introduce a set of distributions concerning a given SMS that
play a central role in our analysis. In \cref{sec:calibration} we use these distributions to
formally define the notion of calibration between two SMSs. The material in this section is presented
in terms of human mathematicians who cohere (or don't) with a putative far-future community of mathematicians.
In the next section, we present similar definitions for embed-calibration, which describes how a scientist
can cohere (or not) with the physical universe that they are both embedded in and that they are predicting.
Then in \cref{sec:math} we describe how calibration can be related to various proposed notions of (mathematical) truth.

Following these formal preliminaries, in \cref{sec:multiplelines} we present our first application of the SMS framework,
to show that if an SMS (set of human mathematical reasoners) is calibrated with a far-future community of mathematicians,
then having multiple lines of evidence for a mathematical hypothesis makes that hypothesis more likely.
(All proofs can be found in Appendix~\ref{sec:proofs}). We also
show in that section that the analogous result holds for scientists who have multiple lines of evidence for a hypothesis
about the outcome of a future experiment.

Next, in \cref{sec:abduction} we present our second application of the SMS framework,
to show that if an SMS (set of human mathematical reasoners) is calibrated with a far-future community of mathematicians,
then abductive reasoning can be applied to mathematical propositions that support one another.  We also
show in that section that the analogous result holds for scientists who use abductive reasoning
concerning the outcome of an experiment that has already happened but whose result they did not observe.
After that, in \cref{sec:earlier_work}, we describe how SMSs are related to other frameworks considered in earlier literature.
We end with a discussion of other applications of the SMS framework that we intend to pursue in future work.

As a final comment, it is important to emphasize that
we do \textit{not} claim that our approach can capture every possible vagary of mathematical and scientific reasoning.
We do not even claim to capture every aspect of human reasoning that can be formulated as
generating sequences of asked and answered questions. Rather, we aim to provide a compelling model that is useful in novel ways.

\section{Stochastic Mathematical Systems}\label{sec:sms}

\subsection{Fundamentals}\label{sec:fundamentals}

%U()At its core, an SMS is a stochastic process for generating
%mathematical claims, i.e., generating mathematical questions and associated answers.
To have the SMS framework be as broadly
applicable as possible, we impose no \textit{a priori} restrictions on the space of
possible ``mathematical questions'' or the space of ``associated answers.''  For example, a question could specify a particular
formal system as well as a well-formed formula (WFF) in that formal system, and ask whether that formula is a theorem of that system. In this case the answer would be a bit. As another example, a question could specify only the form of a desired solution to a mathematical problem, so that there are many possible answers.
%, rather than specify one particular possible solution and ask whether that solution is indeed a theorem.
For example, the question could be of the form, `what is an equation giving the roots of any single-real-variable quartic polynomial?'. Alternatively, a question could be as prosaic as `what is $485 + 923874$?'. Our framework can be applied with any of these types of questions.\par

We start with the following definitions:
\begin{Def}
The elements $q$ of an arbitrary set $\mathcal{Q}$ are known as \textbf{questions}. The elements
$v$ of an arbitrary set $\mathcal{V}$ are known as \textbf{answers}.
\end{Def}
\noindent
As mentioned above,
we impose no restrictions on $\mathcal{Q}$ or $\mathcal{V}$. In particular, we impose no syntactical restrictions;
% on questions and answers;
there is no need to specify a set of rules for what constitutes a well-formed question or answer in our formalism. We also do not require either $\mathcal{Q}$ or $\mathcal{V}$ to be countable. However, to ground thinking and to relate the framework here to the similar framework based on probabilistic Turing machines investigated in \citet{wolpertkinneyNDRmachines2020}, it may benefit the reader to think of the case where both $\mathcal{Q}$ and $\mathcal{V}$ are the set of all finite bit strings.\par

Given a set of questions $\mathcal{Q}$ and a set of answers $\mathcal{V}$,
%We write the set of all pairs of one element of each set
%as $\mathcal{C} = \{(q, v) : q \in \mathcal{Q}, v \in \mathcal{V}\}$. We also
we make the following associated definitions:
\begin{Def}
A \textbf{claim} $c$ is an arbitrary element of $\mathcal{C} = \{(q, v) : q \in \mathcal{Q}, v \in \mathcal{V}\}$.
%$\mathcal{C}$.
\end{Def}
\begin{Def}
A \textbf{claim vector} ${C}$ is a finite sequence of claims.
\end{Def}
\noindent
The set of all claim vectors is written as $\mathcal{C}^*$. We make analogous definitions for
question vectors, answer vectors, etc.
%In general, we use $|C|$ to indicate the number of claims in any claim vector $C$.
\begin{Def}
A \textbf{claim set} $\hat{C}$ is a finite (unordered) set of claims. The collection of all claim sets is written as $\hat{\mathcal{C}}$.
\label{def:claim_set}
\end{Def}
\noindent
%, where that distribution is defined via a measure on an implicit topology.\par

\begin{figure}
    \centering
    \includegraphics[scale=.4]{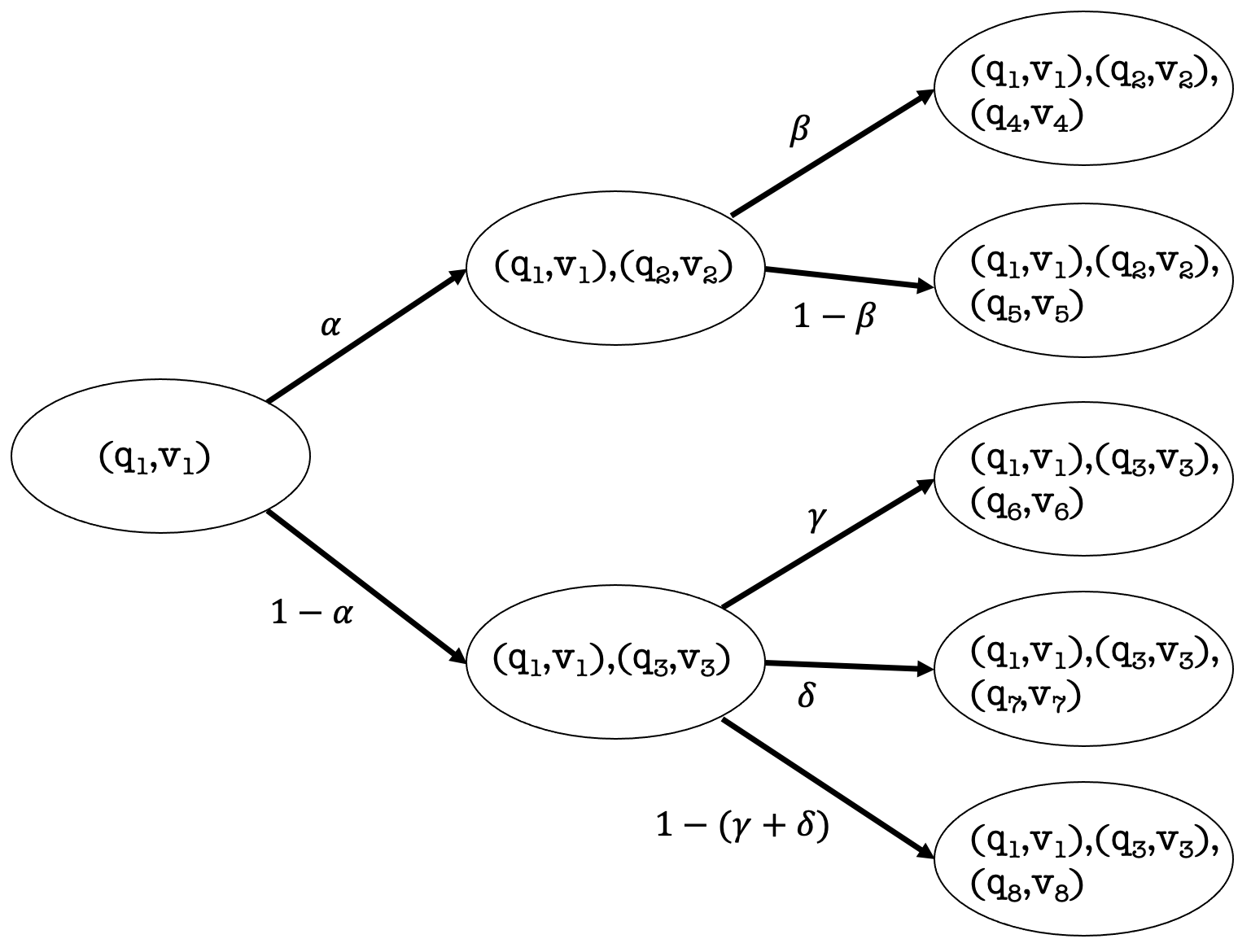}
    \caption{
%\dhwc{I think this figure is a holdover from our NDR paper, where we were considering Turing machines? There, bit strings made sense, but here I think it just confuses the reader to no benefit. Strings of question/answer pairs instead?} \dk{See new image and caption.}
A directed graph showing several of the possible evolutions of the claims that are output by an SMS. Labels on arrows show transition probabilities from each claim to the next, which are determined by the probability space $\mathfrak{C}$ of the SMS.}\label{fig:evolution}
\end{figure}

We suppose that claim vectors are iteratively generated via a discrete stochastic process, which we formalize with the following pair of definitions:
\begin{Def}
A \textbf{claim vector probability space} is a probability space  $\mathfrak{C}=(\Omega_{\mathfrak{C}},\Sigma_{\mathfrak{C}},P_{\mathfrak{C}})$ where $\Omega_{\mathfrak{C}},\Sigma_{\mathfrak{C}}$ and $P_{\mathfrak{C}}$ are the event space, sigma-field, and probability measure, respectively.
\end{Def}
%In general, we will write the space of all probability distributions over some set $Z$ as $\Delta_Z$.
\begin{Def}\label{def:sms}
A \textbf{stochastic mathematical system (SMS)} $\varphi$ is a pair $(\mathfrak{C},X)$ where $\mathfrak{C}$ is a claim vector probability space and the measurable function
$X:\mathbbm{Z}^{+} \times \Omega_{\mathfrak{C}} \rightarrow\mathcal{C}^*$ specifies a sequence of random variables (i.e., of measurable functions) taking values in $\mathcal{C}^*$.
\end{Def}
\noindent
We refer to the integer argument of $X(.)$ as the \textbf{step} of the SMS.
We use $P_{\mathfrak{C}}$ to indicate probabilities of events under  $(\mathfrak{C},X)$, and
use superscripts to indicate the step, e.g., writing $P_{\mathfrak{C}}^n(C)$ for the probability
that the SMS $(\mathfrak{C},X)$ produces the claim vector $C$ when the step is $n$.
Sometimes we will abuse notation and write for example $X(n)$ for some integer $n$ to mean
a sample of $P^n_{\mathfrak{C}}$.
Sometimes we will leave the precise specification of the SMS $\varphi$ implicit, and so for example
just write $P^n(C)$.
%\footnote{We could
%instead define the domain of an SMS in terms of all finite tuples of steps,
% $\cup_{k \in \mathbbm{Z}^+}[\mathbbm{Z}^+]^k$, rather than in terms of single timesteps $\mathbbm{Z}^{+}$,
%and / or introduce filtrations, etc. None of that extra complexity will be needed here.}
For simplicity, we assume that any SMS has zero probability of ever producing the empty set claim vector. We also assume there is zero probability of a claim vector that has multiple copies of the same claim. See Fig.~\ref{fig:evolution} for a visual representation of the process whereby an SMS produces claims.
%In general, we use subscripts to indicate lengths of claim vectors or claim sets, and superscripts to indicate steps or sequences of steps.
%In particular, $P^m(C^m)$ is the probability distribution over claim vectors $C^m$
%produced in the $m$'th step.

Whenever we (implicitly) invoke a topology on any of the sets ${\mathcal{C}}, {\mathcal{Q}}, {\mathcal{V}}, {\mathcal{C}}^*, {\mathcal{Q}}^*, {\mathcal{V}}^*$ (e.g., when discussing measures defined over such sets), that topology is
assumed from context. In particular, for countable $\mathcal{C}$ we assume the associated discrete topology.
%We make similar implicit assumptions when we (implicitly) invoke some measure in addition to $P_{\mathfrak{C}}$.

We can illustrate these definitions with some simple examples:
\begin{enumerate}
\item Choose $\mathcal{Q}$ to be all well-formed formulas (WFFs)
in some formal system, and choose $\mathcal{V} = \mathbb{B}$.
Then we could interpret the claim $(q, 1)$ to mean that $q$ is a theorem and $(q, 0)$ to mean that $q$ is not a theorem. The SMS could be a function that
iteratively produces the decidable WFFs in the formal system, along with
the bit of whether they are (not) a theorem, and adds those claims
to the set of all previously decided claims to produce a new, larger claim set.
\item
Fix some universal Turing machine (TM) $T$, and write $T(p)$ for the (perhaps undefined)
computation of $T$ starting with bit string $p \in \mathbb{B}^*$. Again choose $\mathcal{Q}$ to be all WFFs
in some formal system, but now choose $\mathcal{V} = \mathbb{B}^*$. Then we could interpret
the claim $(q, v)$ to mean that whether $q$ is a theorem or not is given by $T(v)$. In this case a
given answer $v$ paired with a question $q$ can be viewed as a ``translation'' of the WFF $q$ into
an input string to the TM. Alternatively, viewing $T$ as a (partial) function, $v$ can be viewed as an encoded
version of the bit of whether $q$ is or is not a theorem of the formal system. Whatever interpretation
we adopt, as in the previous example the SMS could be a function that
iteratively produces the decidable WFFs in the formal system, along with
the bit of whether they are theorems, and adds those
claims to the set of all previously decided claims to produce a new, larger claim set.
\item Using some appropriate (and somewhat arbitrary) definition of what it means for a specific question-answer
pair to be ``commonly accepted'' by the global community of human mathematicians,
we can view that community as having produced a sequence of sets of commonly accepted claims
concerning mathematical questions since 1900.
Each of those claims can be formulated as a question-answer pair, and the
successive sets of commonly accepted mathematical claims produced by
the community is a claim vector. Note that there has been randomness in which precise questions have interested the community of human mathematicians through time. This can be captured in an SMS model of the community of human mathematicians by introducing stochasticity with respect to which questions occur in each successive claim vector produced by that
SMS.\footnote{Of course, the community of mathematicians did much else besides produce sets of claims in those years, e.g., they produced working hypotheses, decided among formalisms to use to address a topic, etc. Those other aspects of the behavior
of that community do not concern us here; we are only concerned with the end products of their behavior, so to speak.}
\end{enumerate}

In what follows, we will not need to consider claim vectors directly, instead only directly considering sets of claims. However, distributions over
such sets of claims are defined in terms of distributions over claim vectors.
Consistently with \cref{def:claim_set} we use a hat above a multi-component vector to indicate we are considering the
un-ordered rather than ordered version of that vector. To make this precise
we define the \textbf{un-ordering} function
%$U : \mathcal{C}_n \rightarrow \hat{\                                                                                                                                                mathcal{C}}_n$
$U : \mathcal{C}^* \rightarrow \hat{\mathcal{C}}$ for any
%positive integer $n$ and any $C_n \in \mathcal{C}_n$,
$C \in \mathcal{C}^*$
by requiring that every element of $U(C)$ is a component of $C$, and vice-versa.\par

We define a claim vector $C$ to be \textbf{(non-)repeating} if it does (not) contain two claims that have the same question.
%(If it is not non-repeating, then it is \textbf{repeating}.)
We say that an SMS is non-repeating after step $k$ if with probability $1$ it only
produces non-repeating claim vectors  after step $k$. As an example, if the SMS of the community of mathematicians is non-repeating after step $k$, then in all such steps,
there might be hidden contradictions lurking in the set of all claims currently accepted by mathematicians, but there are not any \textit{explicit} contradictions, since no single question arises more than once in a claim vector
made at any step after $k$, and so cannot be given two different answers in the same step. Though we do not need to assume an SMS is non-repeating in most of what follows, there is no loss of generality if one makes this assumption.
%, which can be helpful guide to intuition.
In addition, in
almost all uses of an SMS to model reasoning, one would expect that SMS to be non-repeating after a certain step $k$.\par

\section{Important Distributions in the SMS Framework}\label{sec:distributions}

\subsection{Notation for distributions over claim sets}
Often we will not be concerned with the probability
that at some step $n$ an SMS outputs a particular claim vector $C$,
% or even any element of $U^{-1}(\hat{C})$ for a particular claim set $\hat{C}$,
but rather with the probability that at some step $n$ the SMS outputs any
vector $C$ such that its \textbf{unordering } $U(C)$ \textit{contains} some particular $\hat{C}$
as a subset. It may even be that the probability of the SMS producing a $C$ such that
$U(C) = \hat{C}$ is zero, even though the probability of a $C$ such
that $\hat{C} \subset U(C)$ is nonzero.

As an example, suppose that some mathematician has nonzero probability
of producing some set ${\hat{B}}'$ of multiple theorems (i.e., of multiple question-answer pairs) at some
specific timestep $n$,
%and also has nonzero probability of producing some set ${\hat{B}}''$ then, where both
where
${\hat{B}} \subset {\hat{B}}'$.
% and ${\hat{B}} \subset {\hat{B}}''$.
It may be that nonetheless,
the SMS has zero probability of
producing ${\hat{B}}$ at that time step. (Intuitively speaking, the theorems
comprising ${\hat{B}}'$ can only be derived
all at once, as a package.) Despite this specific character of that SMS,
to analyze that SMS we will sometimes be interested in distributions conditioned only
on the theorems in ${\hat{B}}$, not on the set of
all of the theorems in ${\hat{B}}'$.
% or ${\hat{B}}''$ .
%For example, if may be that there is another
%``package of theorems'', ${\hat{B}}'' \ne {\hat{B}}'$, where ${\hat{B}} \subset {\hat{B}}''$, and we do not care whether
%the package ${\hat{B}}'$ is generated at timestep $n$ or the package ${\hat{B}}''$, only that \textit{some} package of theorems
%containing ${\hat{B}}$ is generated.\par

Phrased formally (with some abuse of notation), we must distinguish the probability distribution
%\begin{linenomath}
\eq{
\label{eq:first_dist}
\PP^{n}(\hat{C}) &:=\sum_{C' \in\mathcal{C}^{*}:\hat{C}= U(C')} P^n(C')
}
from
\eq{
{P}^{n}(\hat{C}) &:= \sum_{C' \in\mathcal{C}^{*}:\hat{C} \subseteq U(C')} P^n(C')
\label{eq:second_dist}
}
%\end{linenomath}
The first distribution takes as its argument the claim set $\hat{C}$, and returns the probability that the (implicit) SMS produces, at step $n$, any claim vector whose unordering {is} the set $\hat{C}$. By contrast, the second distribution takes as its argument the
same claim set $\hat{C}$, and returns the probability that  at step $n$, the SMS produces any claim vector whose unordering \textit{is a superset of} the set $\hat{C}$. So,
\eq{
\label{eq:3.3}
P^n(\hat{C}) &= \sum_{\hat{C}' : \hat{C} \subseteq \hat{C}'} \PP^n(\hat{C}') \ge \PP^n(\hat{C})
}
%
%As discussed above, it may be that
%$\underline{P}^{n}(\hat{{C}}) \ne P^{n}(\hat{{C}})$ for some particular $\hat{C}$.
The proofs below use the distribution $P^{n}(\hat{{C}})$, not $\PP^{n}(\hat{C})$. This
reflects the fact that the associated results do not concern whether an SMS produces a particular set of question-answer outputs,
but rather whether it produces any set that contains a particular set as a subset. See \cref{app:alt_distribution} for some further
remarks concerning the distribution $\underline{P}^n$.

%We refer to the random variable
%$(\hat{\underline{C}}$ as a \textbf{thinned} claim set. Note that it is a member of $\hat{\mathcal{C}}$, and
%so we can consider its union with other (possibly non-thinned) claim sets, etc.

%\dhwc{Must check that everything below still holds even though now we are using a ``thinned''
%distribution. In particular, ensure that the response distribution
%for $\underline{P}^{n}(\hat{{C}})$ becomes a well-defined limit distribution
%for a non-repeating, backward consistent (i.e., ``convergent'') SMS.
%If can't do that, have a problem even in v7, that calibration is not defined for a convergent $\varphi_1$.}
%%the distribution in \cref{eq:second_dist} equals zero even though
%%the distribution in \cref{eq:first_dist} does not.
%%In this paper, we adopt the convention that $P^{n}(\hat{C})$ always refers to the
%%distribution in \cref{eq:first_dist}.

\subsection{Response Distributions}
In much of what follows, it will be useful to consider the probability that, in its $n$'th step, an SMS produces a
claim set that includes both a particular claim $(q,v)$ and a particular set of claims $\hat{C}$.
Abusing notation, we write this probability as follows:
%\begin{linenomath}
\begin{equation}\label{eq:3.44}
    P^{n}((q,v),\hat{{C}}) := P^{n}(\{(q,v)\} \cup \hat{{C}}) = \sum_{\{C' \in\mathcal{C}^{*}:\hat{C}\cup\{(q,v)\}\subseteq U(C')\}}P^n(C')
\end{equation}
%\end{linenomath}
Similarly, we again abuse notation to write
%\begin{linenomath}
\eq{
    P^{n}(q,\hat{{C}}) &:=
%P^{n}(\{\{(q,v)\} \cup \hat{{C}}   : v \in \mathcal{V}\} ) \\
%			&= \!\!\!\!\!\!
\sum_{\{C' \in\mathcal{C}^{*}, v \in \mathcal{V}:\hat{C}\cup\{(q,v)\}\subseteq U(C')\}}P^n(C')
\label{eq:3.55}
}
%\end{linenomath}
$ P^{n}(q,\hat{{C}})$ and $ P^{n}((q, v),\hat{{C}})$ are not properly normalized distributions over their arguments,
in general.\footnote{This is due to the fact that each $\hat{C}$ can occur in more than
one claim vector $C'$, and so summing over all  $\hat{C}$ means we are double-counting vectors
$C'$, in general.}
We refer to $ P^{n}(q,\hat{{C}})$ as a \textbf{question semi-distribution} for a question $q$ and claim set $\hat{C}$ and step $n$.
%This in turn allows us to define a \textbf{question semi-distribution} for a question $q$ and claim set $\hat{C}$ and step $n$:
%%\begin{linenomath}
%\begin{equation}
%    P^{n}(q,\hat{C}):=\sum_{v\in\mathcal{V}}P^{n}((q,v),\hat{C})
%\end{equation}
%%\end{linenomath}
%$P^{m}(q,\hat{C})$
is the probability that at step $n$, an (implicit) SMS outputs a set that contains both the claim set $\hat{C}$
and a claim that contains the question $q$.\footnote{Note that this differs from $P^{n}(\{\{(q,v)\} \cup \hat{{C}}   : v \in \mathcal{V}\} )$,
which is the probability that the SMS outputs a claim set that contains $\hat{C}$ and also contains every claim of the
form $(q, v)$ for some $v \in \mathcal{V}$.}

Combining these two definitions, we arrive at the following:
\begin{Def}
\label{def:responsedist}
For an SMS $\varphi$, step $n$, question $q$, and  claim set $\hat{C}$
where $P^n(q, \hat{C}) \ne 0$, the associated
\textbf{response distribution} is the map from all $v\in{\mathcal{V}}$ to the value
%\begin{linenomath}
\begin{equation*}
    P^n(v | q,\hat{C})=\frac{P^n((q, v), \hat{C})}{P^n(q, \hat{C})}.
\end{equation*}%\end{linenomath}
\end{Def}
\noindent
By comparing \cref{eq:3.44,eq:3.55}, we see that $P^n(v | q,\hat{C})$ is a properly normalized distribution
over $v \in \mathcal{V}$ (assuming $P^n(q, \hat{C}) \ne 0$).
Intuitively, we can interpret it as the distribution over possible answers to a question $q$ that might be given by a
particular SMS at a step $n$
in response to a question $q$, given that the SMS also produces a claim set that contains $\hat{C}$
at step $n$.
We generalize the notation in \cref{def:responsedist} in the obvious way: any probability we write of an indexed
set of $m$ answers conditioned on an indexed set of $m$ questions (perhaps also conditioned on a claim set) is defined in terms of $m$
claims that are ``broken in two'', with the answer components of those claims on one side of the conditioning bar
and the associated question components on the other.

We extend these definitions to the case of multiple claims in addition to $\hat{C}$ in the obvious way,
e.g.,
\eq{
    P^{n}((q_1,v_1), \ldots, (q_m,v_m),\hat{{C}})
			&= \!\!\!\!\!\! \sum_{\{C' \in\mathcal{C}^{*}:\hat{C}\cup\{(q_1,v_1), \ldots, (q_m,v_m)\}\subseteq U(C')\}}P^n(C')
\label{eq:3.44a}
}
and
\eq{
P^n(v_1, \ldots, v_m | q_1, \ldots, q_m, \hat{C})=\frac{P^{n}((q_1,v_1), \ldots, (q_m,v_m),\hat{{C}})}
					{P^n(q_1, \ldots, q_m, \hat{C})}
}
(Note that in this expression the $m$ questions and answers are indexed, in order to ``match up'' with one another,
while the questions and answers in $\hat{C}$ are not.)

Finally, we identify a special class of response distributions of a given question and claim set:
\begin{Def}\label{def:sure}
For an SMS $(\mathfrak{C},X)$, step $n$, question $q$, and  claim set $\hat{C}$
where $P^n(q, \hat{C}) \ne 0$, $P^{n}(v | q,\hat{C})$ is \textbf{sure} if and only if it is a delta function about some answer $v$.
\end{Def}
\noindent
To illustrate these definitions, consider the abstract model of a real-word human mathematician as an SMS.
%, discussed above.
In a given step
%in which they update their claim vector,
the mathematician may
carry forward all the question-answer pairs of the claim vector they produced at the previous step,
as the first claims in their new claim vector. They may then construct one or more new questions,
and then consider those questions, ending the step by (perhaps) producing answers.
Prior to such consideration of a new question, the human mathematician may have beliefs about
the possible answers to any single such question. We can interpret the response distribution of the mathematician to that question at the current step (conditioned on the previous claim vectors) as precisely such a set of beliefs. For some questions (e.g., `what is the value of $2+2$ according to the axioms of Peano arithmetic?'), this response distribution may be sure. For other questions, the mathematician's response distribution will likely \textit{not} be sure. As an example, consider a human being who aims to evaluate a difficult
definite integral using solely pen and paper. That is, they are posed the question, `what is your best guess for the value of $\int_a^b f(x)\text{d}x$?', where $f(x)$ is a highly complicated function. Even with significant mathematical training, if the stochastic process in which the question is posed is to represent an actual mathematical reasoner, then the response distribution of possible answers to that question is likely unsure, at least at some step in the SMS process of that reasoner.\par

\subsection{Backward-Consistent SMSs and the Limit Distribution}

The infinite-step limit of certain types of SMS will play a key role in our analysis.
This role is motivated by using an SMS to model the generation of the commonly accepted
body of mathematics by the community of human mathematicians.
It seems very unlikely that the actual community of human mathematicians would produce a claim
that is universally agreed to have been proven, and so included in the commonly accepted body of mathematics,
but later on that community determines that that precise claim is wrong, and so removes it from the
commonly accepted body of mathematics. To be clear, we do \textit{not} mean to say that mathematicians do not rely on \textit{assumptions} that are later found not to hold. Rather, we hold that at an certain level of being well-posed, the answers to mathematical questions are rarely overturned. To illustrate, the parallel line postulate was long thought to be an axiom of geometry in a very general way, until the discovery of non-Euclidean geometry. However, that discovery did \textit{not} lead mathematicians to revise the answers to any well-posed questions about whether a given proposition follows or does not follow from the parallel line postulate. Going forward, we will assume that mathematical questions are well-posed in just such a way that retrospective revision is unlikely, even when axioms presupposed in the posing of said questions are ultimately found not to hold in full generality.\par

We formalize this property by saying that an SMS is \textbf{backward-consistent} after step $\kappa$ if for all steps $j > \kappa$,
there is zero probability of the SMS producing a sequence of claim vectors $(X(1)=C^1, \ldots, X(n)=C^j)$ such
that for some $i$ where $j \ge i > \kappa$,
%which there are two integers $i, j$ with the following properties:
%\begin{enumerate}
%\item $n \ge j > i \ge \kappa$;
%\item
%%the vector
$U({C}^i) \not \subseteq U(C^j)$.
%differs from the vector of the first $|C^i|$ claims at the beginning of ${C}^j$, i.e.,
%for some $1 \le m \le |C^i|$, $C^i(m) \ne C^j(m)$.
%\end{enumerate}
If an SMS is backward-consistent after step $\kappa$,
then at any step after $\kappa$, with probability $1$ the claim set produced by the SMS
includes the claim set of the previous step. When the precise value of $\kappa$ does not matter,
we will sometimes say simply that an SMS is ``backward-consistent'' without specifying $\kappa$.

The following important property of backward-consistent SMSs is proven in the appendix:
\begin{lemma}\label{lem:limit}
For any SMS that is backward-consistent and any claim set $\hat{C}$,
$\lim_{n \rightarrow \infty} P^n(\hat{{C})}$ is well-defined.
\end{lemma}
\noindent
We refer to $\lim_{n \rightarrow \infty} P^n(\hat{{C})}$ as
% ${\overline{P}}$ as
% this distribution as
the \textbf{limit distribution} of the SMS, and as shorthand write it as ${\overline{P}}(\hat{C})$. The limit distribution
of an SMS defines associated
limit question semi-distributions and limit response distributions of an answer $v$ to a question $q$.
For any $m > 1$, we define limit response distributions
over $m$-tuples of answers conditioned on $m$-tuples of questions similarly.\par

\section{Calibration Between SMSs}\label{sec:calibration}

%\dk{I tried to make these first few paragraphs a bit more gentle and accessible.}
Much of our analysis below does not concern a single SMS, but rather the relationship between different SMSs.
In particular,
% as mentioned in Sec.~\ref{sec:earlier_work},
%Sec.~\ref{sec:fundamentals},
%often we will be interested in the case where
it is possible that the answer that one
SMS gives to a particular question can be accurately interpreted
as a specification of the \textit{probability distribution} over possible answers to a different question.
%In particular, it may be that the SMS $\varphi_2$ specifies the probability distribution over possible answers to a question $q_{1}$ that some \textit{other} SMS $\varphi_{1}$ might output when it considers $q_{1}$.
%This can be illustrated straightforwardly in ordinary language cases involving individual people; f
As a concrete example, consider two human reasoners:\ Alice and Bob. One can think of both Alice and Bob as SMSs, which we label $\varphi_{1}$ and $\varphi_{2}$, respectively. Bob (the SMS $\varphi_{2}$) might consider the question `what is the probability distribution over possible answers that Alice might give to the question $q_{1}$?', where $q_{1}$ is a question in Alice's question set.
Alternatively, instead of representing the human reasoner Alice, the SMS $\varphi_{1}$ might represent the physical universe itself. If Bob is a scientist, he might ask the question
`what is the probability distribution over possible outcomes of an experiment specified as $q_{1}$?' --- in this
special case, since the SMS $\varphi_1$ represents the physical universe, the answer that $\varphi_{1}$ gives to $q_1$ is
the result of the experiment.\footnote{Note that it is important that our framework allows
$\varphi_2$ to make predictions for distributions over the answers of $\varphi_1$, not
just for the answer directly. After all, scientists almost \textit{never}
directly give a single prediction for the outcome of an experiment, but rather give distributions over such
outcomes. }

%The most important such relation for our purposes is the \textit{calibration} relation that can hold between two SMSs.
In both of these kinds of cases, it is natural to begin thinking not just descriptively, but also normatively. For
example, we might want to know what it means for Bob to give \textit{accurate} answers about the probability distribution over possible answers that Alice might output in response to the question $q_{1}$. Moreover, we may be
interested in this kind of accuracy in the relationship between two SMSs when
$\varphi_2$ and $\varphi_1$ have both already generated some superset of a claim set $\hat{C}$.
We will refer to such relationships by saying that
the SMS $\varphi_{2}$ is ``calibrated'' at step $n$ with the SMS $\varphi_{1}$ for a given question $q$,
answer $v$, and  claim set $\hat{C}$.
%This means that, conditional on $\varphi_{2}$ having previously output the claim set $\hat{C}^{n-1}$, $\varphi_{2}$ at step $n$ assigns approximately the same probability to question $q$ having the answer $v$ as $\varphi_{1}$ does under its limit distribution.
To illustrate, for some  claim set $\hat{C}$, $\varphi_{2}$ and $\varphi_{1}$ are calibrated for the $\varphi_1$
question `what is the value of $1+1$?'
if, conditioned on both $\varphi_1$ and $\varphi_2$ having generated the claims in $\hat{C}^{}$
(perhaps along with other claims), the
answer that $\varphi_{2}$ provides at step $n$ to the question `what is the probability that $\varphi_1$ answers $z$
 in response to the question, `what does $1+1$ equal?'' is
close to the actual probability distribution over answers that are outputted by $\varphi_{1}$
%does indeed generate the answer `$2$'
%under its limit distribution
in response to the question, `what does $1+1$ equal?'.
%, conditioned on $\varphi_{1}$ having also produced the set $\hat{C}^{n-1}$.

We will define calibration broadly enough to capture the relationship between any two
SMSs. However, in the current section we will be concerned with SMSs engaged in explicitly mathematical (rather than physical) reasoning. Moreover, we will typically be concerned with calibration
between two specific kinds of SMSs. The first of these SMSs will be the entire far-future community of mathematical reasoners, denoted $\varphi_{1}$.
The second SMS will represent a finite (possibly singleton) group of mathematical reasoners, denoted $\varphi_{2}$.
We will typically be interested in the specific case where the probability distribution over answers output by $\varphi_1$ is
defined in terms the limit distribution of $\varphi_1$, i.e., the distribution over answers to mathematical questions that will be generated in the infinite long run of mathematical practice.
A present-day mathematician $\varphi_{2}$ is calibrated with that far-future
community $\varphi_{1}$ (for a particular question and current claim set of the current mathematician)
if with high probability, that current mathematician assigns approximately the same probability
to the possible answers to the question as does
the community of mathematicians under its limit distribution.

\subsection{Formal definition of calibration}
\label{sec:calibration_1}

To make the concept of calibration fully formal, first
let $\psi: \mathcal{Q}_{2} \rightarrow \mathcal{Q}_{1}^{*}$ be a partial function from the question set of the SMS $\varphi_{2}$ to the set of all finite vectors of questions that can be asked by the SMS $\varphi_{1}$.
%The function $\psi$ has the same basic interpretation as the function $\eta$ defined above;
Intuitively, we will interpret
any question $q \in\dom\psi \subseteq \mathcal{Q}_{2}$ as asking `what is the joint distribution assigned to the possible answer(s) to the question(s) specified in $\psi(q)$, under the distribution of $\varphi_{1}$?'.
We indicate the number of components of $\psi(q)$ for any specific $q \in \rm{Dom} \, \psi$ by $|\psi(q)|$.
%By using $\psi$ we can have the SMS $\varphi_{2}$ pose questions about the joint probability distributions over possible answers to questions that can be asked by the different SMS, $\varphi_{1}$.\par

Next, let $\Psi$ be a partial function mapping $\mathcal{Q}^m_1 \times \mathcal{V}_2$ for all
positive integers $m$ into the associated
%$\Delta_{\mathcal{V}_{1}^*}$ the
space of distributions over
$\mathcal{V}_{1}^m$, i.e., the space of distributions over the set of
all vectors of $m$ answers that might be output by $\varphi_{1}$.\footnote{As a technical note, we allow
there to be values $m$ such that
the partial function $\Psi$ is undefined for \textit{all} arguments that involve $m$ elements of $\mathcal{Q}_1$.}
Intuitively, we interpret any $\Psi(\psi(q),v)$ as the specification of a particular joint probability distribution over possible
answers to the question(s)
specified in $\psi(q)$. So, when a given SMS $\varphi_{2}$ considers the question `what is the probability distribution over possible answer(s) to the question(s) specified in $\psi(q)$?' and outputs the answer $v\in\mathcal{V}_2$,
the function $\Psi$ translates this answer into a probability distribution over the set of possible answers to the question(s) in $\psi(q)$ that might be given by the \textit{other} SMS, $\varphi_{1}$.\par

Note that in general, for any question $q\in\mathcal{Q}_{2}$ that we want to interpret as asking for a distribution over answers
to questions $\psi(q)\in\mathcal{V}^{*}_{1}$, there may be a non-delta function probability distribution over the possible associated
answers of the SMS $\varphi_{2}$. In other words, there is more than one possible distribution over answers to $\psi(q)$
that the SMS $\varphi_{2}$ might specify.
%\dhwc{Integrate the following several paragraphs with one another.}
To illustrate, an SMS may be a science journalist,
able to generate questions of the form, `what is the probability distribution over possible answers by
a randomly chosen mathematician from the current community to the question `does $\texttt{P}=\texttt{NP}$?'?'.
The journalist might produce two answers to this question, with equal probability. One
answer could be, ``the probability that a randomly chosen mathematician believes $\texttt{P}=\texttt{NP}$ is $.95$'', and the other one could
be  ``the probability that a randomly chosen mathematician believes $\texttt{P}=\texttt{NP}$ is $.5$''.
In this case, the journalist would be just as likely to
become convinced by their reading of the literature that $95\%$ of mathematicians
believe $\texttt{P}=\texttt{NP}$
%(assigning that proposition the probability of $.95$),
or instead become convinced that an equal number of mathematicians
believe  $\texttt{P}=\texttt{NP}$ and believe $\texttt{P} \ne \texttt{NP}$.
However, our definition of calibration below would penalize the journalist for having
such uncertainty concerning the probability that they think (that a randomly chosen mathematician in) the community will assign to
the possible answers, in that there is in fact only one probability that the community
will indeed assign to the possible answers.

A current mathematician (or journalist) will not make predictions concerning the answers of the far-future community of mathematicians
``in a vacuum'', without any knowledge of present-day mathematics, as though they lived tens of thousands of years ago.
Rather they will make those
predictions based on a claim set $\hat{C}$ that includes many of the results that the current community of
mathematicians has come to agree on. In this paper, we are interested in how accurate that
mathematician's predictions for the answers of
the far-future community are supposing that that far-future community \textit{also} agrees with the claims in $\hat{C}$. Therefore
all of our results will be formulated where both distributions involving $\varphi_1$ and those involving $\varphi_2$ are
conditioned on some (identical) claim set $\hat{C}$. (Note that this will require that $\mathcal{\hat{C}}_{2}  \subseteq \mathcal{\hat{C}}_{1} $.)

Note also that there are $q \in \mathcal{Q}_2$ that we do \textit{not} want to interpret as asking for a distribution
over answers by $\varphi_1$. This is why $\psi$ is a partial function. Similarly, there are $v \in \mathcal{V}_2$ that we do not want to interpret as specifying a distribution over possible answers to some vector of question $\psi(q)$. This requires that $\Psi$ too be a partial function.
It will be convenient to restrict attention to $q$ and $\hat{{C}}$ such that the arguments of
both $\psi$ and $\Psi$ are within their respective domains of definition. We do that with the following shorthand:
\begin{Def}
Suppose we are given two SMSs $\varphi_1$ and $\varphi_2$ where $\mathcal{\hat{C}}_{2}  \subseteq \mathcal{\hat{C}}_{1} $.
Suppose we are also given a partial
function $\psi: \mathcal{Q}_{2} \rightarrow \mathcal{Q}_{1}^{*}$, and a partial function $\Psi$ that maps
$\mathcal{Q}^m_1 \times \mathcal{V}_2$ into the associated
space of distributions over $\mathcal{V}_{1}^m$ for any positive integer $m$. Let $n^1, n^2$ be two positive integers.
A pair of a question $q \in \mathcal{Q}_2$, and a claim set
$\hat{C} \in \hat{\mathcal{C}}_2$, is a \textbf{prediction pair} (for $\varphi_1, \varphi_2, \psi, \Psi$, $n^1$ and $n^2$) if
\begin{enumerate}
\item $P^{n^2}_{2}(v|q,{\hat{C}})$ is well-defined;
%\item  $\psi(q)$ is a well-defined vector of questions in $\mathcal{Q}_{1}$;
\item  $q \in \dom \psi$;
\item   For all $v \in \mathcal{V}_2$ such that $P^{n^2}_{2} (v \;|\; q,{\hat{C}})$
is nonzero,
%$\Psi(\psi(q), v)$ is a well-defined distribution over $\mathcal{V}_1^{|\psi(q)|}$
$(\psi(q), v) \in \dom \Psi(., .)$;
\item $P^{n^1}_1\left( \mathcal{V}_1^{|\psi(q)|} \;|\; \psi(q), {{\hat{C}}} \right)$
is a well-defined distribution over $\mathcal{V}_1^{|\psi(q)|}$;
\end{enumerate}
where we use subscripts to distinguish the distributions of the two SMSs.
\label{def:prediction_pair}
\end{Def}
\noindent
Note that there are some implicit requirements in
this definition, e.g., for $P^{n^2}_{2} (v|q,{\hat{C}})$ to
be well-defined, it must be that $P^{n^2}_2(q,{\hat{C}}) \ne 0$.
%Note that we use subscripts to distinguish the distributions of the two SMSs in this most
%general definition. However, b

The requirement in \cref{def:prediction_pair} that
$\mathcal{\hat{C}}_{2}  \subseteq \mathcal{\hat{C}}_{1} $ means that any question that can be asked by $\varphi_{2}$ can be asked by $\varphi_{1}$ and every answer that can be given by $\varphi_{2}$ can be given by $\varphi_{1}$. So for example,
any question a current mathematician might
ask can also be asked by the far-future community of mathematicians, and similarly for the answers that can be given
by the current mathematician and those of
the far-future community of mathematicians. The reason for this requirement is to ensure that both $\varphi_1$ and $\varphi_2$
can be conditioned on the same claim set $\hat{C}$, e.g., so that the set of claims accepted by the far-future community of mathematicians
includes those accepted by the current mathematician.

% This

%Recall that a distribution of an SMS whose argument
%list contains a claim set $\hat{{C}}$ and one or more questions  is the probability distribution that at
%the specified step the SMS produces a claim set that contains both $\hat{{C}}$ and those question(s).
%In addition, f
For any two probability distributions $p, r$ both defined over the same (sigma algebra with) event set $Y$,
we write ${{D}}[p(Y), \; r(Y)]$ for some convex divergence measure quantifying the difference between the distributions $p$ and $r$.
(The canonical example is the Kullback-Leibler divergence; see \citet{cover_elements_2012}.)
%~\cite{{cover_elements_2012}.)
Furthermore, for any $ m \in \Z^+$,
we write $Y^m$ for the event space of all $m$-tuples of elements of that set $Y$.

Given this notation, we define the calibration relation between an SMS $\varphi_{2}$ and an SMS $\varphi_{1}$ as follows:
\begin{Def}
Let $\varphi_1, \varphi_2$ be two SMSs where $\mathcal{\hat{C}}_{2}  \subseteq \mathcal{\hat{C}}_{1} $,
with associated partial functions $\psi$ and $\Psi$ that have domains and ranges as described in \cref{def:prediction_pair}.
Let $n^2$ be an integer, $\hat{C}$ a claim set in $\hat{\mathcal{C}}_2$,
and $q$ a question in $\mathcal{Q}_2$.
%such that $P^{n^2}_2(\hat{C}) \ne 0$.
%$\hat{C}^{\prime}$ such that an element of $U^{-1}(\hat{C}^{\prime})$ is
%which is produced by $\varphi_2$ at step $n_2$ with
%nonzero probability.
We say that $\varphi_2$ is \textbf{calibrated} (with $\varphi_1$) at steps $n^1, n^2$ for the pair
$(q, \hat{C})$ and for $\epsilon \ge 0$ iff
% \dk{should this be if and only if?}
\begin{enumerate}
%\item $\mathcal{Q}_{2}\subseteq\mathcal{Q}_{1}$ and $\mathcal{V}_{2}\subseteq\mathcal{V}_{1}$;

\item   $(q, {\hat{C}})$ is a {prediction pair} (for $n^1, n^2$);

\item $\sum_{v \in \mathcal{V}_2} P^{n^2}_2(v | q, \hat{C}) D\left[ \Psi(\psi(q), v)(\mathcal{V}_1^m),
	{P}^{n^1}_1\left( \mathcal{V}_1^m \;|\; \psi(q), {{\hat{C}}} \right) \right]
                \le \epsilon  $
\end{enumerate}
where $m$ is shorthand for $|\psi(q)|$, and we use subscripts to distinguish whether quantities are defined with respect to SMS $\varphi_1$ or $\varphi_2$.
\label{def:calibrated}
\end{Def}
\noindent
In practice, we will usually leave $\epsilon$ implicit when discussing
calibration. Note that calibration is not a symmetric relation; for any $\epsilon$, if $\varphi_{2}$ is calibrated with
$\varphi_{1}$, then it is not necessarily the case that $\varphi_{1}$ is calibrated with $\varphi_{2}$.
Note also that for $(q, {\hat{C}})$ to be a prediction pair, it must be that both $P^{n^2}_2(q, \hat{{C}}) \ne 0$
and $P^{n^1}_1(\psi(q), \hat{{C}}) \ne 0$.

%From now on we restrict attention to the special case of  \cref{def:calibrated}
%where $\varphi_1$ is backward-consistent and $n^1 = \infty$
%In this special case, the three conditions in \cref{def:calibrated} for $\varphi_{2}$ to be calibrated with $\varphi_{1}$ at step $n$
%for the pair $(\hat{C},q)$ are, in sequence:
%\begin{enumerate}
%\item [1)]

The first condition in  \cref{def:calibrated}
%means that any question that can be asked by $\varphi_{2}$ can be asked by $\varphi_{1}$ and every answer that can be given by $\varphi_{2}$ can be given by $\varphi_{1}$. The reason for this requirement is to ensure that both $\varphi_1$ and $\varphi_2$ can be conditioned on the same  claim set $\hat{C}$. The second condition in  \cref{def:calibrated}
means that for
%
%%There is nonzero probability that in its $n$'th step, the SMS $\varphi_2$
%%will generate the question $\psi(q)$ concerning the answers of the SMS $\varphi_1$,
%%after having generated (a claim set which has as a subset) $\hat{C}$.
%%So $P^n_{2} (v \;|\; q, {{\hat{C}}})$ is a well-defined distribution;
%
%\item [2)]  For
any answer $v$ that $\varphi_{2}$ assigns with positive probability to question $q$ given the
claim set $\hat{C}$, $\Psi(\psi(q), v)$ is a probability distribution over the set of all vectors of $|\psi(q)|$
answers that $\varphi_{1}$ can generate.\footnote{Note that $\varphi_1$ and / or $\varphi_2$ might consider a particular
%(a claim vector with associated
%thinned claim set given by)
${{\hat{C}}}$ very unlikely \textit{a priori}.
This suggests an alternative definition of calibration to \cref{def:calibrated},
under which $P_1^{n^1}$ and / or $P_2^{n^2}$ are not conditioned on ${{\hat{C}}}$. We do not investigate this
alternative here, for reasons for space.}
%(Note that for $(q, \hat{C})$ to be a prediction
%pair implicitly requires that $P^n_2(q, {\hat{C}}) \ne 0$.);
%
%%\item [3)] The average of the predictions that $\varphi_2$ makes for what the response distribution
%%of $\varphi_1$ will be for the question $\psi(q)$, given
%%that $\varphi_2$ has already made the claims in the  claim set ${{\hat{C}}}$,
%%is well-defined;
%
%\item [3)] I
From now on we restrict attention to the special case of  \cref{def:calibrated}
where $\varphi_1$ is backward-consistent and $n^1 = \infty$. In this case we write $\overline{P}_1$ instead of $P^{n^1}_1$ in the fourth requirement of
\cref{def:prediction_pair}, and write $n$ instead of $n^2$.
Under these conditions, the second condition
%the three conditions
in \cref{def:calibrated}
%for $\varphi_{2}$ to be calibrated with $\varphi_{1}$ at step $n$
%for the pair $(\hat{C},q)$
%are, in sequence:
says that, in expectation conditioned on $\hat{C}$, the SMS $\varphi_{2}$ generates
answers to questions of the form `what is the probability distribution over vectors of $m$
answers to the associated $m$ questions in $\psi(q)$?' that diverge by less than $\epsilon$ from
$\varphi_{1}$'s limit distribution over possible vectors of $m$ answers to the $m$ questions in $\psi(q)$,
conditioned on the same claim set $\hat{C}$.\par

This special case is crucial to several applications of our framework described below. In many instances, we will take $\varphi_{1}$ to be a mathematical ``oracle'' that acts as the arbiter of what counts as the ``correct'' response distribution for a given question and claim set.
In particular, as described above, one can adopt the interpretation where
$\varphi_{1}$ represents the process by which all mathematicians in the universe ask and answer questions, from the
past on to the future. Under this interpretation, in keeping with the Peircean idea of truth as what wins out in the long run of intelligent inquiry, we can model an individual mathematician $\varphi_{2}$ as answering questions correctly to the extent that they are calibrated, at some present step $n$, with the SMS $\varphi_{1}$ in its infinite limit. That is, the mathematician $\varphi_{2}$ answers questions correctly to the extent that it answers them by sampling from a response distribution that is minimally divergent from the limit response distribution of $\varphi_{1}$, conditioned on the same claim set $\hat{C}$.
%\end{enumerate}

%Note that calibration is not a symmetric relation; for any $\epsilon$, if $\varphi_{2}$ is calibrated with
%$\varphi_{1}$, then it is not necessarily the case that $\varphi_{1}$ is calibrated with $\varphi_{2}$.
%Note as well that
%if $\varphi_{1}$ were non-repeating as well as backward-consistent, then there would a step after which it never revises its answers to any question $q$.
%Thus, we take the convergence assumption to characterize a mathematical oracle for which we a primarily interested in its limit distribution, where no questions are taken to be ``open problems.''\par

One might object to the real world relevance of calibration,
arguing that in practice, few human mathematicians ask themselves questions of
the form, ``what is \textit{the probability distribution over answers} to the following question (that would be given
by some oracle)?''. Instead of those kinds of questions, human mathematicians typically
just directly ask themselves, ``what is \textit{the answer} to the
following question?'' (which we can interpret as them asking themselves, ``what is the answer
that would be given by an oracle to the following question?'').
It might seem that the concept of calibration has nothing to say about this second,
more common kind of question. However, in our approach,
the second kind of question is just a special case of the first.\par

To see this formally, suppose that the two SMSs in \cref{def:prediction_pair} are actually identical.
In this case we do not use subscripts, and so for example we simply have a single integer
$n$ and have $\psi : \mathcal{Q} \rightarrow \mathcal{Q}^*$. We also do not make the
first requirement in \cref{def:prediction_pair} for such a \textbf{single-SMS prediction pair}. Choose
$\psi(.)$ to be the identity function (so $|\psi(q)| = m = 1$), choose $\Psi$ so that for all
$v_{1},v_2 \in \mathcal{V}_{1}=\mathcal{V}_2$, $\Psi(q, v_2)(v_1) = \delta(v_1, v_2)$
(the Kronecker delta function).
For those choices, and for any prediction pair $(q,\hat{C})$, \cref{def:calibrated} reduces to the condition
that
%\begin{linenomath}
\begin{equation}\label{eq:4.2}
    D[P^n_{2} (\mathcal{V}_{1} \;|\; q, {{\hat{C}}}), \overline{P}_1(\mathcal{V}_{1} \;|\; q, {{\hat{C}}})]\leq\epsilon.
\end{equation}
%\end{linenomath}
So for those choices, a current mathematician-SMS $\varphi_2$ is calibrated with the
oracle $\varphi_1$
if and only if their response distribution is sufficiently non-divergent from that of $\varphi_1$,
% (to within $\epsilon$),
for the question $q$.\footnote{To derive \cref{eq:4.2}, note that according to the assumptions made immediately above it,
%it, the distribution $P^{n}_{2}(\mathcal{V}_{1}|q,\hat{C})$ is such that, for any $v_{1}\in\mathcal{V}_{1}$,
$P^{n}_{2}(v_{1}|q,\hat{C})=\sum_{v_{2}}P^{n}_{2}(v_{2}|q,\hat{C})\delta(v_{1},v_{2})$
for any $v_{1}\in\mathcal{V}_{1}$. With abuse of notation, for any $v_{2}\in\mathcal{V}_{2}$, let $\delta(\mathcal{V}_{1},v_{2})$ be a vector such that each entry is the value of the delta function $\delta(v_{1},v_{2})$ for each $v_{1}\in\mathcal{V}_{1}$. Similarly, let $P^{n}_{2}(\mathcal{V}_{1}|q,\hat{C})$ be a vector such that each entry is the value of $P^{n}_{2}(\mathcal{V}_{1}|q,\hat{C})$ for each $v_{i}\in\mathcal{V}_{1}$. Thus, $P^{n}_{2}(\mathcal{V}_{1}|q,\hat{C})=\sum_{v_{2}}P^{n}_{2}(v_{2}|q,\hat{C})\delta(\mathcal{V}_{1},v_{2})$. Since $D$ is convex in its first argument, we have $D[P^{n}_{2}(\mathcal{V}_{1}|q,\hat{C}),\overline{P}_{1}(\mathcal{V}_{1}|q,\hat{C})]=D[\sum_{v_{2}}P^{n}_{2}(v_{2}|q,\hat{C})\delta(\mathcal{V}_{1},v_{2}),\overline{P}_{1}(\mathcal{V}_{1}|q,\hat{C})]\leq\sum_{v_{2}}P^{n}_{2}(v_{2}|q,\hat{C})D[\delta(\mathcal{V}_{1},v_{2}),\overline{P}_{1}(\mathcal{V}_{1}|q,\hat{C})]\leq\epsilon$.}

%\dhwc{Clean up / cannibalize the following paragraph.}
%This effectively states that if a mathematician-SMS \textit{cannot} ask questions about the limit distribution over possible answers to a given question, but can instead only output answers to that question, then they are calibrated with an arbiter of mathematical norms to the extent that their \textit{process} for answering those questions (represented by sampling from a probability distribution) is sufficiently similar (as measured by a convex divergence) to the process used by the arbiter. Thus, under these conditions the oracle or other arbiter of mathematical norms functions as an \textit{exemplar} of mathematical reasoning, with an individual mathematician or group of mathematicians calibrated with that exemplar to the extent that their reasoning mimics it.

There are several important reasons that the concept of calibration goes further than
\cref{eq:4.2}. First, it encompass sets of multiple question-answer pairs that come ``as a bundle'', so
to speak. Second, it allows a current mathematician to be explicitly uncertain about what they
think the oracle's answer to a question would be, crediting them (in a normative sense)
if that uncertainty matches the oracle's uncertainty.

There is also an  important distinction between calibration and simply requiring small divergence between the two response distributions
$P^n_{2}(\mathcal{V}| q, {{\hat{C}}})$ and $\overline{P}_1(\mathcal{V}|q, {{\hat{C}}})$.
Even if that divergence is identically zero, an answer given by randomly sampling $P^n_{2}(\mathcal{V}| q, {{\hat{C}}})$
may differ greatly from an answer given by sampling $\overline{P}_1(\mathcal{V}|q, {{\hat{C}}})$, however
we measure distance between answers. In other words,
%
%The SMSs $\varphi_1$ and $\varphi_2$ are independent stochastic processes. Accordingly, there cannot
%be any \textit{a priori} statistical coupling between the answers they give to the same question.
%This means that even if the divergence between $P^n_{2}(\mathcal{V}| q, {{\hat{C}}})$ and $\overline{P}_1(\mathcal{V}|q, {{\hat{C}}})$
%equals zero, if those two (identical) answer distributions had high entropy, there would be high probability that
the answer $v$ which the mathematician $\varphi_2$ generates in response to the question $q$ may differ greatly
from the answer $v$ which the oracle $\varphi_1$ generates in response to that question.
Due to this,
there is no normative basis to requiring that that divergence be small.

On the other hand, there \textit{is}
a normative basis to requiring that $\varphi_2$ be calibrated with $\varphi_1$; when $\varphi_2$ is calibrated with $\varphi_1$,
the answers of $\varphi_2$ to any question in a prediction pair are (on average) very \textit{close} to
the answers  of $\varphi_1$ (as quantified by the
divergence between the two associated distributions specified by those answers). This requirement amounts to requiring that (on average) $\varphi_2$ answers the question `what is the probability distribution over possible answers to $\psi(q)$?' by outputting a probability distribution with low divergence from the probability distribution that $\varphi_1$ outputs in response to the same question. So in this case,
the answer $v$ which the mathematician $\varphi_2$ generates in response to the question $q$ \textit{cannot} differ greatly
from the answer $v$ which the oracle $\varphi_1$ generates in response to that question (on average).
So when calibration is small, $\varphi_2$ answers mathematical questions about probability distributions over possible
mathematical facts in the same way that the oracle will (in the infinite long-run).

At first, it may seem strange to suppose that a mathematical \textit{oracle} would ever be unsure of the answer to any mathematical question. However, it is clearly possible under the far-future-community interpretation of such an oracle. Suppose that there is simply eternal disagreement on the answer to a mathematical question $\psi(q)$: some say that the answer is $v$, while others say that the answer is $v^{\prime}$ Eventually, the mathematical community converges on a stable distribution over possible answers:\ with probability $\rho$, $v$ is output as the answer, while with probability $(1-\rho)$, $v^{\prime}$ is output as the answer. Under these circumstances, we hold that the best way for a present-day mathematician to respond to the question `what is the distribution over possible answers to $\psi(q)$?' is to (in expectation) produce a distribution with low divergence from that which the mathematical community converges on. Thus, while our approach does not \textit{require} that such answers are the best way to answer mathematical questions, it allows for this possibility. We take this to be a virtue of our account, one that proves all the more fruitful when we turn our attention to scientific reasoning.

\subsection{Honesty}

Note that the requirement in \cref{def:prediction_pair} and \cref{def:calibrated} that $\mathcal{\hat{C}}_{2}  \subseteq \mathcal{\hat{C}}_{1} $
implies that there can be  pairs
$q \in {\mathcal{Q}}_1, v \in \mathcal{V}_1$ such that both $q \in \dom(\psi)$ and $(\psi(q), v) \in \dom(\Psi)$.
Given how we wish to interpret $\psi$ and $\Psi$, this means that an SMS $\varphi$ can ask the question,
 `what is the joint probability distribution that \textit{I assign} to the possible answer(s) to the question(s) specified in $\psi(q)$?' ---
and then provide an answer. This provides what can be seen as a necessary condition for $\psi$ and
$\Psi$ to have the interpretation we wish to assign them:

\begin{Def}
Let $\varphi$ be an SMS with associated single-SMS partial functions $\psi$
and $\Psi$.
% as in the special case described below \cref{def:prediction_pair} where there is a single SMS.
Let $q$ be an associated question, $v$ an associated answer, $n$ a positive integer, and $\hat{C}$
an associated claim set. Suppose that $(q, \hat{C} \cup \{(q, v)\})$ is a (single SMS) {prediction pair}.
Then we say that $\varphi$ is \textbf{honest} (for $(q, v, \hat{C}, n, \psi, \Psi)$)) if for all $\tilde{v} \in \mathcal{V}^{|\psi(q)|}$, $P^n(\tilde{v} \,|\, \psi(q), \hat{C} \cup {\{(q, v)\}}) = \Psi(\psi(q), v)(\tilde{v})$
\label{def:honest}
\end{Def}
\noindent
Note that in \cref{def:honest}, $\hat{C} \cup {\{(q, v)\}}$ is
the claim set given by the union of $\hat{C}$ and the single pair $(q, v)$.
We make the obvious extension of removing $n$ from the definition of honesty
to allow consideration of the limit distribution of a backward-consistent SMS rather than a step-$n$ distribution.
We will sometimes simply say that $\varphi$ is honest, without specifying the precise $q$ and $v$, if it is honest
for all $(q, v)$ such that $q \in \dom \psi, (\psi(q), v) \in \dom \Psi$.

Intuitively, an SMS is honest if it is perfectly calibrated (i.e., calibrated for $\epsilon = 0$) with itself.
Viewed alternatively, an SMS is honest if it does indeed ``mean'' $\psi(q)$ to be a set
of questions, and ``means'' $\Psi(\psi(q), v)$ to be an associated joint distribution over $|\psi(q)|$ answers.
As an example, suppose that we interpret  $\varphi_1$ in \cref{def:calibrated}
as the far-future community of mathematicians.
%
%, it may help if they
%presume that the limit distribution of that far-future community is honest for the
%provided $\hat{C}$ and $(q, v)$.
Then requiring that the far-future community
is an honest SMS (for the provided $\hat{C}$ and $(q, v)$) essentially amounts to requiring that they
use what \citet{lewis1971immodest} calls ``immodest'' inductive methods. Once they accept a given response distribution $\Psi(\psi(q), v)$ as the correct one for assigning answers to a vector of questions,
they make predictions about the answers to those questions by sampling from that distribution $\Psi(\psi(q), v)$.
%since it takes that distribution to provide, in expectation, more accurate predictions than any other.

Strictly speaking, a current mathematician can be perfectly calibrated with the far-future community
even if they are not honest. If this were the case, then (the distribution over) their personal answers to the
questions in $\psi(q)$ would not be the same as (the distribution over) the answers that they predict the far-future community
would provide. For some reason or other, when asked, ``do you think $\texttt{P}=\texttt{NP}$?'' they would respond by sampling from a different distribution over
answers than the one they would sample if instead asked, ``do you predict that the far-future community of mathematicians
would think $\texttt{P}=\texttt{NP}$?''.

While our analysis allows this, typically one would presume that those two
distributions must be identical, essentially as a normative principle. Under this presumption, \cref{def:honest}
would be a \textit{definition} of $\Psi(., .)$, a definition that one would make before
defining prediction pairs or calibration. However, many of our results below do not require that $\varphi_2$ be
honest, allowing the current mathematician (for whatever reasons) to give a different set of
answers to a given set of questions from the set of answers they predict that the far-future community
would give. Accordingly, here we do not impose honesty, formulating it as a property that an SMS
may or may not have, rather than as a requirement.

%For similar reasons, it may help the reader to ground their thinking if they presume that
%a present-day mathematician who is calibrated with that far-future community for
%partial functions $\psi$ and $\Psi$ is also honest for the same partial functions.

%Note though that it is not meaningful to even ask whether a given SMS is honest unless we have also defined the partial
%functions $\psi$ and $\Psi$. Moreover, \textit{a priori}, even if $\psi$ and $\Psi$ are defined for some SMS $\varphi$, and we have a $q \in \dom(\psi)$ and a $v \in \dom(\Psi(\psi(q), v)$, it may be that $\varphi$ is not honest for
%any $n, \hat{C}$ for that pair $(q, v)$. Also, none of our results  below rely on an SMS being honest (although
%some possible problems with interpreting our results are obviated if the associated
%SMSs are honest).

\section{SMSs, Calibration and Empirically Grounded science}\label{sec:science}

While the details are still unclear, the deepest understanding of the nature of physical reality
that modern physics provides us is quantum cosmology, and various possible elaborations of
quantum cosmology like string theory~(see \cite{barrow1991theories,greene1999elegant,carroll2021quantum}). Crucially,
%in the quantum cosmology community
the consensus view among a large subset of quantum cosmologists is that quantum cosmology is a \textit{complete} description of physical reality,
in the sense that every experiment that a scientist could ever run (up to and including the ``experiment'' of posing questions to
themselves internally), and every observation they could ever make, is completely described either directly
by the laws of quantum cosmology or indirectly via formal implications of those laws. (Examples of those implications
being the rules of quantum chemistry, which in turn ultimately provide the statistical regularities of terrestrial biology, etc.)
Specifically, the modern view  is that the physical universe can be formulated \textit{in toto} as some
appropriate formal system~(see \cite{tegmark1998theory,tegmark2014our}). The outcomes of all experiments
ever conducted by scientists, all observations, either in the past or the future, are theorems in that formal system, all following from a particular
set of axioms. The crucial point is not that we have all those axioms in hand already --- we do not --- never mind that we could
even imagine calculating those theorems from those axioms. Rather, the crucial point is that physical reality can
be formulated as \textit{some} formal system, in such a way that almost tautologically, there
is no conceivable role played by some additional ``reality'' not captured by that formal system, no
sense in which such a reality could be physically accessible to us humans. In modern physics,
there is no unavoidable need for the concept of a ``physical reality'' outside of the mathematical laws of quantum cosmology.

The idea that the physical universe is ultimately a mathematical object and nothing more is explicitly promoted by, among others, Barrow (\citeyear{barrow1991theories}, \citeyear{barrow2011godel}), \citet{dipert1997mathematical}, \citet{schmidhuber1997computer}, in addition to Tegmark. The view that the universe is fundamentally mathematical is also \textit{close} to, and arguably entailed by, \textit{ontic structural realist} views in philosophy of science, such as those put forward by \citet{Ladyman1998-LADWIS-2}, \citet{ladyman2007scientific}, and \citet{pittphilsci20048}. (However, note that some ontic structural realists want to say that, although the physical world is a structure, it is a distinctly \textit{physical}, rather than mathematical, structure.)
For our part, though we offer no argument for it here, we are happy to embrace the literal interpretation of the general claim that the physical universe is a mathematical object.
%, as an entailment of the literal interpretation of the specific claim that it is a stochastic mathematical system.\par
We then extend the perspective of modern physics that the physical universe as a mathematical object is a formal
system, in the obvious way, to formulate the physical universe as an SMS.
%Just as representing mathematical truth as a stochastic process allows us to represent the possibility that mathematics might be inherently stochastic, r
%The decision to represent the physical universe as an SMS
%has two important consequences. First, it means that
%we wholly adopt the idea that the universe is essentially \textit{mathematical}. After all, an SMS is ultimately just a mathematical object:\ it is a pair consisting of a probability space $\mathfrak{C}$ and a function from $\mathbbm{Z}^{+}\times\Omega_{\mathfrak{C}}$ into $\mathcal{C}^{*}$, and of course $\mathfrak{C}$, $\mathbbm{Z}^{+}$, $\Omega_{\mathfrak{C}}$ and $\mathcal{C}^{*}$ are each, themselves, mathematical objects.
%
%Second, the claim that the physical universe is a stochastic mathematical system

Note that representing the physical universe as an SMS
entails that, as a matter of physical reality, there may be some questions that do not have a single prescribed answer, but are instead only answered by sampling from a non-delta-function distribution over possible answers. This is one way of saying that the physical universe is \textit{fundamentally stochastic}. %Representing the physical universe as an SMS allows for the possibility that the universe itself is inherently stochastic.
%That is, the claims or facts that obtain in a particular universe may be determined by sampling from a non-degenerate response distribution over possible answers to questions, given some claim set $\hat{C}$. Moreover, there may be facts of the matter as to the non-degenerate response distribution over possible outcomes of some vector of experiments.
%
There are many ways in which one might think that there are physical questions whose answers are determined in
such a fundamentally stochastic way. On the Copenhagen interpretation of quantum mechanics, the result of any measurement of a quantum system is determined stochastically. Alternatively, one might endorse a ``many-worlds'' interpretation of quantum mechanics (\cite{everett1957relative}) in which all possible quantum measurement outcomes occur in some branch of a multiverse structure, but still follow \citet{sebens2018self} in holding that it is an objective fact that observers in a given branch of that structure ought to assign specific non-extreme probabilities to their being in that branch, and to their being in some other possible branches. Or, one might adopt the \textit{Mentaculus view} of cosmology articulated by \citet{Albert2000-ALBTAC} and Loewer (\citeyear{Loewer2007-LOECAT}, \citeyear{Loewer2008-LOEWTI}, \citeyear{Loewer2009-CRAWIT}). On their view, which follows in the Humean tradition of \citet{lewis1986plurality}, our actual physical universe is composed of a spatiotemporal array of perfectly intrinsic properties. However, this actual world is generated by sampling from a probability distribution over the set of all possible physical universes defined as such.\par

All of these theories can be articulated within the SMS framework, where they differ with respect to their specification of the set of questions that are answered by sampling from a non-delta-function probability distribution. For instance, on the Copenhagen interpretation, the universe-SMS might consider a question of the form `will tritium atom $x$ decay by time $t$?', and arrive at the answer `Yes' or `No' by sampling from an unsure response distribution over those two possible answers. On a Sebens-and-Carroll style multiverse view, the universe-SMS considers questions of the form `where am I located in a branching multiverse structure?', and responds by sampling from an unsure response distribution. Finally, on a Mentaculus view, the universe-SMS considers questions of the form, `is the physical universe such that $x$ occurs?', where $x$ is some event, and answers them by sampling from an unsure response distribution over possible universes.
%For instance, the universe-SMS might consider a question of the form `will tritium atom $x$ decay by time $t$?', and arrive at the answer `Yes' or `No' by sampling from a non-sure response distribution over those two possible answers. Moreover, the universe-SMS might consider the question `what is the probability distribution over possible answers to the question ``will tritium atom $x$ decay by time $t$?'?', and respond with a non-degenerate probability distribution over the answers `Yes' and `No'.
%This probability distribution may itself be obtained by sampling from a non-degenerate mixture distribution, or not. Indeed,
However,
representing the universe with an SMS is \textit{also} perfectly consistent with an entirely deterministic universe; for
that to be the case we simply need to model
the universe as an SMS $\varphi_1$ such that the response distribution over all answers to any question $q$ and prior  claim
set $\hat{C}$ of that universe-SMS is sure. %So the SMS approach can embody
%the many-worlds interpretation of quantum
%mechanics, just as it can embody the Copenhagen interpretation.

% As described above, the SMS framework allows us to formalize mathematical reasoning and the conditions under which it can be described as ``correct.''
%Below we use calibration to formalize the validity of such reasoning, in terms of whether it converges
%to correct probability distributions.
%%the SMS framework, and prove that when calibration holds, some of the heuristics used in
%%both types of reasoning are formally justified.
% In addition to providing a model of reasoning in mathematics,
% the SMS framework can also provide a model of reasoning in the empirical sciences. Reasoning in the empirical sciences
% is inherently more complicated than reasoning in mathematics though, since the scientist exists in the
% same mathematical system (namely a physical universe) as the question-answer pairs they are interested in.
% %(See~\cite{bibid} for an investigation of some of the consequences )
% In this section we introduce extra definitions that
% allow us to capture that reasoning. Then in the following sections we prove that when calibration holds, some of the heuristics used in
% both mathematically and scientific reasoning are formally justified.

\subsection{What it means for a scientist to be embedded in a universe}

When formulating a model of reasoning in mathematics,
we used one SMS $\varphi_{2}$ to represent a mathematician and a second SMS $\varphi_{1}$ to represent
an arbiter of ``correct'' answers to mathematical questions, intuitively interpreting that arbiter as the far-future
community of mathematicians. Note that $\varphi_1$ and
$\varphi_2$ had no \textit{a priori} relation --- if the mathematician $\varphi_1$ was able to give
approximately correct answers, then they were calibrated with the far-future
community of mathematicians, but nothing more was formalized concerning their relationship.

%\dhwc{Work in much of pp. 23-25 in the 4.30.21 version of the document.} \dk{See the two paragraphs after the one below. I want to try to do this briefly, without taking a full detour into all of the metaphysics.}

We can model scientists making predictions about the outcomes
of physical experiments similarly --- but with one major difference. Choose $\varphi_2$ to be a
particular scientist, or more precisely an SMS
generating cognitive events of that scientist concerning physical phenomena (i.e.,
observations, memories, etc., formulated as question-answer pairs), much like the
cognitive events of a mathematician that concern mathematical phenomena.
Also choose $\varphi_1$ to be the arbiter of correct answers, just as in the case of mathematical
reasoning. Of course, just as real human mathematicians do
far more than predict how the far-future community of mathematicians will answer specific questions, real
scientists do far more than predict what answers will be provided by physical experiments to specific questions. For instance, they often build models of phenomena that contain explicitly unobservable phenomena. Moreover, scientists' interpretations of the outputs of experiments may be fundamentally shaped by the theories that they have an antecedent commitment to (i.e., their observations may be ``theory-laden'').
In this paper we are not concerned with those other aspects of scientific reasoning, and so abstract away from them for the sake of building our model, though we are optimistic that these other aspects of scientific reasoning could be formalized in a manner consistent with our broader framework. Similarly, we do not consider here the formal models of scientific (dis)agreement developed in the social epistemology of science (see \cite{sep-computational-philosophy} Sec.\ 3.2.2; \cite{Seselja2022-EELAMO-3}), though here too we believe future work could see overlap between formal modeling of the social dynamics of scientific inquiry and the model we provide here.\par

In the case of science, we take the arbiter of correct answers
to be the physical universe generating outcomes of experiments, not some far-future community of
scientists. Moreover, scientists are physically embedded in the very physical universe about which they wish to make predictions.
(After all, scientists are physical beings (or so we assume) and so must be accounted for as part of the
physical universe-SMS). This means that a scientist
$\varphi_2$ is \textit{itself} a set of (question-answer) pairs, posed to the
embedding physical universe $\varphi_1$. There is no analogous property relating a mathematician
and the far-future community of mathematicians. In essence, in the case of mathematical reasoning,
$\varphi_1$ and $\varphi_2$ are independent SMSs, while in the case of science, $\varphi_2$ is
a sub-SMS, embedded in $\varphi_1$.

We formalize what it means for a scientist-SMS to be physically
embedded in a universe-SMS by requiring that the distribution of claim sets that a scientist might make, each of which
is just some pattern in their physical brain, is given by a distribution over states of the universe
that ``project down'' to those patterns in the brain of the scientist. We formalize this by
requiring that there be a measurable (partial) function from the space of claim sets of the universe into
the space of claim sets of the scientist:

%As shorthand, we write
%$[{c}_\varphi; \hat{C}_\varphi]$ for the claim set $\{c_\varphi\} \cup \hat{C}_\varphi$.
% is defined to be the collection of all claim sets
%$\hat{C}'_\varphi \in \hat{\mathcal{C}}_\varphi$ such that
%$\{c_\varphi\} \cup \hat{C}_\varphi \, \subseteq \, \hat{C}'_\varphi$.

\begin{Def}\label{def:embedded}
Let $\varphi_{2}$ be an SMS and $\varphi_{1}$  an SMS that is backward-consistent.
We say that $\varphi_{1}$ \textbf{embeds} $\varphi_{2}$ for step $n$
if and only if there is a partial function $E$ from ${\hat{\mathcal{C}}_{1}}$ into
$ {\hat{\mathcal{C}}}_2$ such that for any non-empty claim set $\hat{C}'_2 \in \hat{\mathcal{C}}_2$,
$\overline{P}_{1}(E^{-1}[\hat{C}^{\prime}_2]) = P^{n}_{2}(\hat{C}^{\prime}_2)$.
\end{Def}
\noindent
%(Note that given some claim set $\hat{C}_\varphi \in \hat{\mathcal{C}}_\varphi$
%and claim $c_\varphi \not \in \hat{C}_\varphi$,  $\{c_\varphi\} \cup \hat{C}_\varphi$ is another claim set.)
We will sometimes say that ``$\varphi_2$ is embedded in $\varphi_1$'' rather than ``$\varphi_1$ embeds $\varphi_2$''.
We will also sometimes refer to $E$ as an \textbf{embedding function} from $\mathcal{\hat{C}}_1$ to $\mathcal{\hat{C}}_2$.
Often below we will leave out specification of the step $n$, assuming that $\varphi_2$ is embedded in $\varphi_1$
for all steps after some threshold step $m$, and that $n > m$.
We refer to $\varphi_2$ in \cref{def:embedded} as the ``scientist'' SMS, and refer to $\varphi_1$
as the ``universe'' SMS.

In the sequel we will be interested in the
case where $\hat{C}'_2$ can be written as $\{c_2\} \cup \hat{C}_2$ for some claim $c_2 \in \mathcal{C}_2$ and
(possibly empty) claim set $\hat{C}_2 \in \hat{\mathcal{C}}_2$.
% the embedding function $E$ in \cref{def:embedded}
%maps any set of sets of facts about the physical universe (i.e., any  claim set ${\hat{C}}_{1}$
%that might be generated by the universe-SMS $\varphi_{1}$) into a claim set given by the union of
%$\hat{C}'_2$ and $c_2$ that could be output by a scientist-SMS.
In these cases the claim $c_2$ will be identified with a question by the scientist of the form, ``what is the
distribution over possible outcomes for this specific experiment / observation?'', together with an
answer to that question. We will then identify the associated
claim set  $\hat{C}_2$ as (a subset of)
all that the scientist has previously observed that is relevant to their answer given in claim $c_2$.

As an example,
% of case (2), suppose  ${\vec{C}}_{2}$ is a trajectory
suppose that $\hat{C}_2$ specifies the (brain) state of a scientist
in which they remember some parts of a textbook concerning radioactive decay, and $c_2$ is
%containing only
the claim (`what is the probability that tritium atom $x$ decays by time $t$?', `$.9$'). In this case, we interpret
$E^{-1}[\{c_2\} \cup \hat{C}_2]$ as the collection of all physical
processes that might occur in the universe, within the brain of the scientist and (perhaps)
elsewhere in the universe, each of which is sufficient for the scientist to
have the (brain state in which) they think they have read some parts of a textbook
concerning radioactive decay and also have the thought that the tritium atom
$x$ decays by time $t$ with probability $.9$.

Note that in general there will
be many possible facts concerning the universe that are all consistent with some set of fixed claims made
by a scientist in a  claim set $\hat{C}_2$, physical facts that have nothing to do with that scientist or their claims. (As an illustration, there are
many possible numbers of stars in the Andromeda galaxy, and  all of those numbers are consistent with
the claims made by the typical present-day human scientist.) This is captured in \cref{def:embedded}
by the fact that the inverse image of the embedding function will in general be more than just one set
of physical facts --- one claim set of the universe-SMS ---
but rather a collection of multiple (perhaps infinite) such sets. Note as well that
the scientist's distribution in \cref{def:embedded} is of the form of  \cref{eq:second_dist},
not \cref{eq:first_dist}. This reflects the fact that like the state of the universe, the state of the scientist's
brain might have generated other claims besides those specified in $\{c_2\}\cup\hat{C}_2$.

As a formal point, suppose that any claim set produced by the scientist-SMS is
a subset of the claim set produced by the (infinite step limit of the) universe-SMS, i.e., some of the question-answer pairs
of the universe specify the physical properties of the scientist \textit{in toto}, as discussed above. Then \cref{def:embedded}
holds no matter what the universe-SMS is, i.e., no matter what $\overline{P}_1$ is. In this case, ontologically there is only one
SMS, that of the physical universe, with the scientist being in effect a projection the universe-SMS
onto a sub-SMS. In our discussion below
we will usually be informal and phrase the embedding relationship as though the scientist is indeed a sub-SMS of the
physical universe. It's worth pointing out though that
 \cref{def:embedded}  can hold even if $\varphi_1$ and $\varphi_2$ are completely
independent stochastic processes, without the claim set generated by $\varphi_2$ being a subset
of the one generated by $\varphi_1$, so long as the probability distributions of those two stochastic processes
happen to be related appropriately.

It is important to emphasize that $E^{-1}[\{c_2\} \cup \hat{C}_2]$ is a \textit{collection} of claim sets in $\hat{\mathcal{C}}_1$,
{not} a {union} of those claim sets. Formally,
 $E^{-1}[\{c_2\} \cup \hat{C}_2]$ is the union of a set of elements of $\hat{\mathcal{C}}_1$ --- the fact that those elements
happen to be sets is immaterial. So
in terms of the probability measure defining $\varphi_1$,  $E^{-1}[\{c_2\} \cup \hat{C}_2]$ is the event that
$\varphi_1$ generates one of the claim sets of $\varphi_1$ whose image under
$E$ is $\{c_2\} \cup \hat{C}_2$, not the event that $\varphi_1$ generates the union of all such claim sets.
Similarly, the probability $\overline{P}_{1}(E^{-1}[\{c_2\} \cup \hat{C}_2])$ is the
limit probability of  $\varphi_1$ generating (a superset of) any of the claim sets in $E^{-1}[\{c_2\} \cup \hat{C}_2]$.
It is not the probability of $\varphi_1$ generating
the claim set given by taking the union of all the claim sets whose image is $\{c_2\} \cup \hat{C}_2$.
Indeed, due to contradictions among claims in such a union, often the probability that $\varphi_1$
creates that union in its entirety will be zero.
%\footnote{Note
%that for the definition of embedding to have the meaning we wish, we implicitly presume that
%there is no claim set $\hat{C}_1$ which contains two distinct subsets, $\hat{C}'_1$ and $\hat{C}''_1$, both of which have the
%same image under the partial function $E$.}

As a final point, we note that many of the results below hold for any partial function $E$ from ${\hat{\mathcal{C}}_{1}}$ into
$ {\hat{\mathcal{C}}}_2$, even if it
violates the condition on probability distributions in \cref{def:embedded}. For simplicity though, in this paper we always assume that
that condition in \cref{def:embedded} is obeyed.

\subsection{Scientists who are calibrated with the universe that embeds them}

Given this formalization of what it means for a scientist-SMS to be embedded in a given universe-SMS,
we can define what it means for a scientist to accurately predict the outcomes of experiments or observations in the physical universe. This involves a slight modification to the definitions of prediction pair and calibration that were introduced in \cref{sec:calibration_1}.
% \cref{def:calibrated}.
In particular, we do not consider whether a scientist is making predictions consistent with the embedding universe
conditioned only on a set of their earlier experiments / observations, in analogy to the
case with calibration and mathematical reasoning. In addition to conditioning on a claim set, we also condition on the simple
\textit{physical} fact that the scientist is, in fact making a prediction. After all, the fact that they, a part
of the universe, are making a prediction, is a restriction on
what the precise state of the universe embedding them can be.

\begin{Def}
Let $\varphi_1$ be a backward-consistent SMS and $\varphi_2$ an SMS.
%where $\mathcal{\hat{C}}_{2}  \subseteq \mathcal{\hat{C}}_{1} $.
Let $\psi$ and $\Psi$
be associated partial functions that have domains and ranges as described in \cref{def:prediction_pair}.
Let $n$ be a positive integer
and let $E$ be an embedding function from $\varphi_1$ to $\varphi_2$.
A pair of a question $q \in \mathcal{Q}_2$, and a claim set
$\hat{C} \in \hat{\mathcal{C}}_2$, is an \textbf{embedded prediction pair} (for $E, \varphi_1, \varphi_2, \psi, \Psi$ and $n$) if
\begin{enumerate}
\item $P^{n}_{2} (v \;|\; q,{\hat{C}})$ is well-defined;
%\item  $\psi(q)$ is a well-defined vector of questions in $\mathcal{Q}_{1}$;
\item  $q \in \dom \psi$;
\end{enumerate}
For all $v \in \mathcal{V}_2$ such that $P^{n}_{2} (v \;|\; q,{\hat{C}})$
is nonzero:
\begin{enumerate}
\item  [(3)]
%$\Psi(\psi(q), v)$ is a well-defined distribution over $\mathcal{V}_1^{|\psi(q)|}$
$(\psi(q), v) \in \dom \Psi(., .)$;
\item [(4)]  $\overline{P}_1\left( \mathcal{V}_1^{|\psi(q)|} \;|\; \psi(q), E^{-1}[\{(q, v)\} \cup {\hat{C}}] \right)$
is a well-defined distribution over $\mathcal{V}_1^{|\psi(q)|}$;
\end{enumerate}
where as usual we use subscripts to distinguish the distributions of the two SMSs.
\label{def:embedded_prediction_pair}
\end{Def}
The definition of an embedded prediction pair in \cref{def:embedded_prediction_pair} is very similar to the
definition of a prediction pair in \cref{def:prediction_pair}. One difference is that
the requirement in \cref{def:prediction_pair} that $\mathcal{\hat{C}}_{2}  \subseteq \mathcal{\hat{C}}_{1} $ is replaced
in \cref{def:embedded_prediction_pair} with
the requirement that there be a specified embedding function from $\varphi_1$ to $\varphi_2$.
Note though that there is no requirement that $\varphi_1$ actually embeds $\varphi_2$ under that embedding
function (or any other) in \cref{def:embedded_prediction_pair}.

Another difference between the two definitions reflects the fact that the physical act of the
scientist asking question $q$ and giving answer $v$ is a restriction on the possible state of the physical universe
in which they are embedded. This requires us to add $(q, v)$ in the conditioning event in the
probability distribution in condition (4) of \cref{def:embedded_prediction_pair}.
That in turn requires us to stipulate that $v$ have nonzero probability under SMS $\varphi_2$. In contrast, no such stipulation is needed in
\cref{def:prediction_pair}. Sometimes we will wish to replace the requirement in \cref{def:embedded_prediction_pair} that
conditions (3, 4) hold for all $v$ such that $P^{n}_{2} (v \;|\; q,{\hat{C}})$ with the requirement that they hold
for some specific such $v$. In such cases we refer to $(q, \hat{C}, v)$ as an \textbf{embedded prediction triple}.
So in general every embedded prediction pair corresponds to many embedded prediction triples.

A small modification of the definition of calibration is also needed to define the
analogous concept for an empirical scientist:

\begin{Def}
Let $\varphi_1$ be a backward-consistent SMS and $\varphi_2$ an SMS,
where $\varphi_2$ is embedded in $\varphi_1$ with embedding function $E$ (for step $n$).
Suppose we are given associated partial functions $\psi$ and $\Psi$ that have domains and ranges
as described in \cref{def:prediction_pair}. Suppose further that
we are given a claim set $\hat{C} \in \hat{\mathcal{C}}_2$ and a question $q \in \mathcal{Q}_2$.
We say that $\varphi_2$ is \textbf{embed-calibrated} (with $\varphi_1$) at step $n$ for the pair
$(q, \hat{C})$ and for $\epsilon \ge 0$ iff
\begin{enumerate}
%\item $\varphi_2$ is embedded in $\varphi_1$ with embedding function $E$ (for step $n$);
\item   $(q, {\hat{C}})$ is an embedded {prediction pair} for embedding function $E$
\item $\sum_{v \in \mathcal{V}_2} P^n_2(v | q, \hat{C}) D\left[ \Psi(\psi(q), v)(\mathcal{V}_1^m),
	\overline{P}_1\left( \mathcal{V}_1^m \;|\; \psi(q), E^{-1}[\{(q, v)\} \cup {\hat{C}}]  \right) \right]
                \le \epsilon  $
\end{enumerate}
where $m$ is shorthand for $|\psi(q)|$.
\label{def:scientistcalibrated}
\end{Def}
\noindent
%(Recall the definition of an embedded prediction pair, in~\cref{def:embedded_prediction_pair}.)
% and prediction distribution $f$.
As before, the function $D$ in \cref{def:scientistcalibrated}
is some convex divergence measure between probability distributions, and $\epsilon$ is often implicit.

Intuitively, embed-calibration is the same as
calibration \textit{simpliciter}, just with three coupled changes.
First, in \cref{def:scientistcalibrated} we drop the assumption that $\mathcal{\hat{{C}}}_{2} \subseteq \mathcal{\hat{{C}}}_{1} $.
That assumption was needed in condition (3)
of \cref{def:calibrated}, to specify how conditioning on an already accepted claim set by a present-day
mathematician-SMS translates into conditioning on a claim set by the far-future-community SMS.
This assumption is replaced in
\cref{def:scientistcalibrated}
by the assumption that the scientist-SMS is embedded in the universe-SMS. Second, calibration requires that $(q, {\hat{C}})$ be a {prediction pair}. In contrast, embed-calibration requires
that they be an embedded prediction pair.
%, because the embedding function maps claim sets
%of the universe-SMS to claims made by the scientist-SMS.

Finally, in \cref{def:calibrated}
we required that, in expectation, the SMS $\varphi_2$'s answer to the question `what is the probability distribution over possible answers to the questions in $\psi(q)$?', conditional on the  claim set $\hat{C}$, diverges by less than $\epsilon$ from the SMS $\varphi_{1}$'s limit response distribution over possible answers to the equations in $\psi(q)$, conditional on the same claim set $\hat{C}$. By contrast, in the case of the calibration relation between a scientist-SMS and the universe-SMS in which it is embedded, we require that the scientist-SMS's answer to the question `what is the probability distribution over possible answers to the questions in $\psi(q)$?', conditional on the  claim set
$\hat{C}$, diverges by less than $\epsilon$ from the universe-SMS's limit response distribution over possible answers to the questions in $\psi(q)$,
\textit{conditional on the universe-SMS producing one of the claim sets in $E^{-1}[\{(q,v)\} \cup \hat{C}]$}.
%That is, whereas the scientist-SMS's response distribution is conditioned on the claims that it actually makes prior to considering the question $q$, the universe-SMS's response distribution is conditioned on the \textit{physical instantiations} of the scientist-SMS's claims prior to considering the question $q$.\par
%
%\dhwc{Nothing in this definition of embedding calibration requires the flexibility in the definition of embedding
%to handle both claim trajectories and
%pairs of claim trajectories and claim sets. Nor do I see it below. It appears though that
%the two associated ``technical assumptions'' concerning embedding are used, albeit implicitly, below?
%Does some of the material that was removed by DK exploit those assumptions? If so, we should
%discuss this issue here, right after the definition of embedding calibration.}

The goal of the scientist-SMS, broadly speaking, is to be embed-calibrated with the universe-SMS that embeds it.
To illustrate, consider a scientist-SMS who asks themselves
the question `what is the probability distribution over possible answers to the question `will the specific tritium
atom $x$ decay by specific time $t$?' that would be provided by a future experiment, as determined by the physical universe?'.
Suppose that with probability $1$ \textit{conditional on a body of nuclear theory the scientist knows,
which is identified with the  claim set $\hat{C}$}, the scientist generates the answer that the probability distribution is `Yes'
with probability $.35$ and `No' with probability $.65$. If the scientist is embedded in a universe-SMS with a backward-consistent response distribution for the question `will tritium atom $x$ decay by time $t$?'?',
and if \textit{conditional on the universe producing a claim set in $E^{-1}[\{(q,v)\} \cup \hat{C}]$},
that universe assigns probability $.35$ to `Yes' and $.65$ to `No', then that scientist is calibrated for the question `what is the probability distribution over possible answers to the question `will tritium atom $x$ decay by time $t$?'?',  claim set $\hat{C}$, and $m=1$, for any non-negative value of $\epsilon$.
%(Note that the s

Note that real world human scientists \textit{always} have statistical uncertainty about their predictions.
This uncertainty
%Statistical uncertainty accompanying
%predictions by
is automatically respected in our approach, due to the stochasticity of the universe-SMS.
Note also that in these examples $E^{-1}[\{(q,v)\} \cup \hat{C}]$ is \textit{not} a set of physical facts in the embedding universe
that comprise the nuclear theory adopted by the scientist-SMS. Rather, it is a set of physical facts that are jointly constitutive of the scientist-SMS having \textit{the beliefs that they do about nuclear theory}. Concretely, $E^{-1}[\{(q,v)\} \cup \hat{C}]$
specifies some spatio-temporal biochemical pattern in the brain of the scientist-SMS.
On the other hand, the distribution $\Psi(\psi(q), v)$ is a distribution over physical outcomes of an experiment, \textit{not} over
the state of the brain of the scientist after that experiment.
%It may be that the scientist-SMS accepts the claims in $\hat{C}$ only if those claims
%are also made by the embedding universe-SMS, i.e., if they are ``physically
%correct''. On the other hand, it

\section{Calibration and Truth}\label{sec:truth}
\label{sec:math}
%\dk{I've just made some minor touch-ups here.}

%\dhwc{I think it may make sense to expand this subsection into a full section, discussing both our (suggested) notions of mathematical truth and of scientific truth, and then to move it to after we introduce embedding and calibration in science, but before our two applications.}

The nature of
% truth in general, as well as both
mathematical and scientific truth
%articular,
is a notoriously thorny philosophical issue that we cannot hope to fully address here. For this reason we are careful not to use the term ``true'' in our formal definitions,
and to make no explicit claims concerning truth.
Nevertheless, there are some connections
% we offer here some comments on the connection
between the two calibration notions developed above and general, informal notions of ``truth'' which are worth commenting on.

In ordinary human reasoning, we often assert sentences, and moreover, we assert that some of these sentences are true. \citet{dummett1959viii} famously investigates \textit{why} a notion of truth should be part of our apparatus for reasoning, and concludes that the sole function of a notion of truth is to distinguish between those sentences that we can justifiably assert and those that we cannot. Indeed, he concludes that a sentence is true just in case it can be justifiably asserted. While this definition of truth does not yield a procedure for determining whether any given sentence is true (this would require a theory of justifiable assertion), it does provide a theory of the \textit{purpose} of the concept of truth; namely, truth has the normative function of distinguishing between those sentences that it is correct to assert and those that it is not correct to assert. One can easily extend this account to give an analysis of what it means to answer a question $q$ truthfully; $v$ is the true answer to $q$ just in case one is justified in answering $q$ with $v$. Note that the idea that truth provides normative governance for assertions extends to more technical notions of truth. For instance, when we define a semantics for a theory written in first-order logic by specifying a model of that theory, we implicitly stipulate which formulas of the theory are assertable and which ones are not.\par

%In our framework, reasoning is represented as a stochastic process of asking and answering questions. A calibration relation is then defined between SMSs.
On its own, calibration can be understood purely mathematically, without any normative significance.
%However, we can also ask \textit{why} it is important to define a notion of calibration between SMSs. One answer, which we have adopted i
However, in our discussion
of how to interpret the notion of calibration in the previous two sections we have interpreted it as also
having a normative function. This normative function is similar to the function of truth in Dummett's account, but is distinct in important ways. In Dummett's account, truth is a normative constraint on the claims that speakers of a language ought to make. Within the SMS framework, calibration can be interpreted as providing a normative constraint on \textit{how} an SMS ought to make claims.
%If we take calibration with some SMS $\varphi_{1}$ to be the normative goal of an SMS $\varphi_{2}$, then
%Most importantly, the definition of calibration
%given above tells us,
%provides a fully formal definition of what it means for a response distribution of $\varphi_2$ to
For a given question $q$ and claim set $\hat{C}$, one can understand an SMS $\varphi_{2}$ that generates answers to $q$ by sampling from a response distribution that allows $\varphi_{2}$ to achieve calibration with some other SMS $\varphi_{1}$ to be generating claims in a ``justifiable'' way.
%or correct choice of response distribution for $\varphi_{2}$ to sample from when generating answers to $q$.
%Namely, the SMS $\varphi_{2}$ should sample from those distributions that are not too divergent from the response distribution used by $\varphi_{1}$ for the question $q$ and claim set $\hat{C}$.

Note that defining ``justifiable'' in terms of calibration does not imply that any particular answer $v$ is the correct or justifiable one for $q$.
Defining justifiable (or more generally, ``truth'') in terms of
calibration allows for the possibility that the answers to some questions are objectively indeterminate
or chancy. In the context of calibration and mathematical reasoning this might seem peculiar,
 if the oracle is some nebulous version of ``Platonic truth''. However, it is not so unreasonable
if the oracle is interpreted as the far-future community of mathematicians.
(In any case, as discussed above, if in fact the infinite-step limit of
$\varphi_1$
% the mathematical community
\textit{does} output some
single answer to a particular question  --- as in conventional notions of mathematical truth ---
that simply means that the distribution over the answers they
provide is a delta function.)

In more detail, for any present-day real-world mathematician, there is a large set of claims made by earlier mathematicians
that they accept as ``true.'' In the real world, such earlier claims will be part of the claim set of that present-day mathematician. Moreover,
since the generation of claim sets is a stochastic process, the full set of claims accepted
by the present-day mathematician will modify the distribution of claims made by the far-future community of
mathematicians. Indeed, in the extreme case, one particular claim set accepted by the present-day mathematician
might completely rule out some claims of the far-future community, claims that would otherwise be
quite likely. This is the reason why the definition of calibration involves both
conditioning the present-day mathematician's answers on a given claim set and conditioning the
far-future community of mathematicians' answers on the same claim set. In addition, although we do not exploit
the fact here, backward-consistency allows us to require that the full claim set of
that far-future community of mathematicians ``develops out of''
the the claim set of the current community, as the SMS evolves.

The case of empirical science is slightly different.
In empirical science, we take the universe-SMS $\varphi_{1}$ as the arbiter of what counts as a correct or true distribution from which to sample answers to a given question about the outcome of an experiment. A scientist-SMS outputs true or correct answers to questions about an experiment to the extent that they are embed-calibrated (for small $\epsilon$) with the universe-SMS in which they are embedded. As in the case of mathematical reasoning, it may
be that there is no unique correct answer to a well-posed question about the outcome of an experiment, but that there is instead a correct probability distribution over possible answers to a well-posed question about the outcome of an experiment. Indeed, as discussed in Sec.~\ref{sec:science}, our best theories in quantum cosmology imply that there are well-posed questions about the outcomes of experiments that do not admit of a uniquely correct answer, but do admit of a uniquely correct probability distribution over possible answers, (i.e., possible observations of the outcome of a particular experiment). In this sense, our approach is in keeping with Hacking's (\citeyear{hacking1964foundations}) notion of objective statistical facts as grounded in the propensity of a particular ``chance set-up'' to result in a particular outcome. Indeed, Hacking's definition of a chance set-up as ``a device or part of the world on which might be conducted one or more trials, experiments, or observations; each particular trial must have a unique result which is a member of a definite class of possible results'' is very much in the spirit of what counts as a well-posed question for both a universe-SMS and a scientist-SMS (p.\ 3).\par

So in the case of mathematical reasoning calibration
involves the relationship between the response distribution of the far-future
community of mathematicians and that of the current mathematician
to the same set of questions, given that both accept the same claim set $\hat{C}_2$, i.e., given that
both accept the same set of chains of proof-like
mathematical reasoning. In the case of scientists though,  $\hat{C}_2$ is a physically concrete property
of the scientist's own brain, corresponding to their mental perception of
the outcomes of chains of experiments or observations they have made. There is no sense
in which the universe-SMS ``accepts the same claim set''. Rather that claim set is a projection
of the full state of the universe down to a set of properties of the brain of the scientist.

Note also that although the universe-SMS can make claims about what happens at particular points in space and time,
%whether or not the embedding relation holds does not depend on the ordering of the claim sets output by the universe-SMS.
we do \textit{not} presume that the universe-SMS successively generates its claims according to any privileged temporal ordering. In other words, for the universe-SMS, the ordering of steps need have nothing to do with
the ordering of time in some particular inertial frame of the physical universe generated by that SMS.
In this sense,
we effectively adopt a ``block universe'' perspective on quantum cosmology. By contrast, the scientist-SMS
\textit{does} have a privileged temporal ordering to its steps, since it perceives reality to unfold in a particular temporal
direction, from past to present.\footnote{One way of accounting for this ``psychological arrow of time'' is by appealing to the second law of thermodynamics; see \citet{wolpert1992memory} and \citet{davies1977physics}.} This is why we define embed-calibration in terms the limit probability of the universe-SMS; we can only be assured
that all cognitive events in the mind of the scientist will be derived as claims of the universe-SMS if we allow
that universe-SMS to run arbitrarily many steps, without necessarily interpreting those steps as temporally ordered.\par

Finally, we note that in the case of empirical science, it is possible for a scientist to be systematically biased in their
observations of the outcomes of experimental set-ups (e.g., if they are hallucinating). Indeed,
it may be that they are able to make accurate predictions about the outcomes of
experiments but are unable to accurately observe the outcomes of such experiments, and so
they might falsely conclude that their initial predictions were inaccurate. Even more extremely, it may be that they are unable to discriminate
among the possible outcomes of the experiments, due to limitations of their observational
apparatus. In such a case there would be no empirical meaning to their predictions, in the sense that the accuracy of those predictions cannot be verified.\footnote{Note that
there is not an analogous issue
for the case of mathematician-SMSs. That is because we suppose that mathematicians are
interested in ``mathematical truth'', which they may never directly observe in any sense. For
example, this is the case if mathematical truth is interpreted as being the response distribution
of some far-future community of mathematicians that the mathematician will not live long enough to encounter.}
%
%One response is to simply accept this, giving the scientist-SMS credit for making an accurate prediction of the outcome of an experiment whatever (possibly erroneous) previous evidence their predictions are based on and even if they
%misperceive the outcome of the experiment. That is,
%one might say that
%%under this version of our approach
%a scientist reasons correctly insofar as they are embed-calibrated with the universe-SMS in which they are embedded, regardless of the scientist's own internal phenomenology.\footnote{Note that there is not an analogous issue
%for the case of mathematician-SMSs. That is because we suppose that mathematicians are
%interested in ``mathematical truth'', which they may never directly observe in any sense. For
%example, this is the case if mathematical truth is interpreted as being the response distribution
%of some far-future community of mathematicians that the mathematician will not live long enough to encounter.}

%At the other extreme, one
One might object that this is a flaw in the idea of using embed-calibration to assess the quality of a scientist's predictions,
%that this is a shortcoming of our approach,
since ultimately any scientist is only interested in
the distribution over their \textit{perceptions} of the outcome of an experiment, not over the ``actual'' outcomes in some sense independent of their mind.
%We can render this argument moot though. First,
%
In order to address this objection we must first assure
% note that it implicitly assumes
that the perceptions of the scientist \textit{can} discriminate
between the actual outcomes of the experiment. To see how to do this, first note that
depending on the precise form of $E$,
there may be a set $S \in \dom E$ such that $S \subset E^{-1}E(S)$. In this case, having $E(S)$ be the state of the scientist's
brain would not fix the state of the universe to be an element in $S$. So a necessary condition for the observation of the scientist to
discriminate between the actual outcomes of the experiment
is that if we apply the embedding function $E$ to each of those possible outcomes of the experiment, we get
a state of the scientist's brain that uniquely distinguishes among those possible outcomes, via the inverse function $E^{-1}$.

We formalize this as follows:
\begin{Def}
Let $\varphi_1$ be a backward-consistent SMS and $\varphi_2$ an SMS, where $\psi$ and $\Psi$
are the associated partial functions that have domains and ranges as described in \cref{def:prediction_pair}.
Let $n$ be a positive integer and let $E$ be an embedding function from $\mathcal{\hat{C}}_1$ to $\mathcal{\hat{C}}_2$.
We say that an associated embedded prediction pair $(\hat{{C}}, q)$ is \textbf{discriminating} if for all associated embedded prediction
triples $(\hat{{C}}, q, v)$, for all associated $v^m_1 \in {\rm{supp}} \,
\overline{P}_1\left( \mathcal{V}_1^m \;|\; \psi(q), E^{-1}[\{(q, v)\} \cup {\hat{C}}]  \right) $,
\eq{
  E^{-1}E \left[\{(\psi(q), v^m_1)\} \cup E^{-1}[\{(q, v)\} \cup {\hat{C}}]\right]  =  \{(\psi(q), v^m_1)\} \cup E^{-1}[\{(q, v)\} \cup {\hat{C}}] \nonumber.
}
\label{def:distinguisher}
\end{Def}
\noindent
It is important to be clear about what \cref{def:distinguisher} does \textit{not} say. It does not say that $E$
must be an injective function from its domain to its image. Nor does it say that, for any $v\in\mathcal{V}_{2}$, $E^{-1}[\{(q, v)\} \cup {\hat{C}}]=\{(q, v)\} \cup {\hat{C}}$. Thus, $(\hat{{C}}, q)$ being a discriminating pair is still consistent with mental states of the scientist being multiply realizable in terms of outputs of the universe-SMS. Rather, \cref{def:distinguisher} says only that when a scientist poses a question $q$ of the form `what is the probability distribution over answers to the question(s) in $\psi(q)$?', and provides some answer $v\in\mathcal{V}_{2}$, and that answer is interpreted by the function $\Psi$ as specifying some probability distribution $\Psi(q,v)(\mathcal{V}^{m}_{1})$, the scientist-SMS can observationally distinguish between all elements in the support of $\Psi(q,v)(\mathcal{V}^{m}_{1})$, which are answers given by the \textit{universe}-SMS. Note that the scientist-SMS need not even ``speak universe'' for $(\hat{{C}}, q)$ to be a discriminating pair. That is, elements of the support of $\Psi(q,v)(\mathcal{V}^{m}_{1})$ need not be in the scientist-SMS's ``cognitive vocabulary'' of answers $\mathcal{V}_{2}$. All that is needed is that the scientist-SMS has \textit{some} way of individuating between elements of the support of $\Psi(q,v)(\mathcal{V}^{m}_{1})$, even if they do so using different ``terminology'' (i.e., they use elements of $\mathcal{V}_{2}$ to distinguish between elements of $\mathcal{V}^{m}_{1}$ via a translation process).\par

This definition of a discriminating pair is in keeping with an informal understanding of the
neurobiological process of observation. To illustrate, note that images in the human retina are displayed upside-down relative to how they are observed. So the scientist might be such that their observations systematically lack fidelity with the environment that they aim to represent. Nevertheless, the scientist is able to observe scenes of their environment because each of their mental representations --- upside-down though they may be ---  maps to a distinct set
of possible actual scenes. On the other hand, for a scientist to observe whether there is a particular pattern in ultraviolet on a flower does not require that the image be right-side up in their brain. Rather it requires that they be able to discriminate among different patterns in ultraviolet.
If they can only observe in the visible part of the spectrum, then they cannot make such an observation. In the same way, if a pair $(\hat{C},q)$ is  discriminating, then it is not necessarily the case that a scientist makes observes possible answers to the questions in $\psi(q)$ in exactly the same form that the universe outputs them. But it is necessarily the case that they can make observations that allow them to discriminate between these possible answers.\par

Combining \cref{def:distinguisher} with \cref{def:embedded} and \cref{def:scientistcalibrated} establishes the following:
\begin{prop}
Let $\varphi_1$ be a backward-consistent SMS and $\varphi_2$ an SMS, where $\psi$ and $\Psi$
are the associated partial functions that have domains and ranges as described in \cref{def:prediction_pair}.
Let $n$ be a positive integer.
%and let $E$ be an embedding function from $\varphi_1$ to $\varphi_2$.
Suppose that  $\varphi_2$ is {embed-calibrated} with $\varphi_1$ at step $n$ for the embedding function $E$, the pair
$(q, \hat{C})$, and for $\epsilon = 0$. Finally, suppose that $(q, \hat{C})$ is a discriminating pair. Then
for all $v^m_1 \in \mathcal{V}^m_1$,
\eq{
\Psi(\psi(q), v)(v_1^m) &=
	\frac{{P}^n_2\left( E\left[\{(\psi(q),  v_1^m)\} \cup E^{-1}[\{(q, v)\} \cup {\hat{C}}]\right]  \right)}
		{{P}^n_2\left( E\left[\{(\psi(q))\} \cup E^{-1}[\{(q, v)\} \cup {\hat{C}}]\right]  \right)}   \nonumber
}
\label{prop:projection}
\end{prop}
\noindent
To understand \cref{prop:projection} intuitively, suppose that
(like in the analysis of calibration for mathematical reasoning) $\mathcal{C}_2$ were a subset
of  $\mathcal{C}_1$. Suppose further that $\psi(q) \in \mathcal{Q}_2$, and so the values $v^m_1$ in \cref{prop:projection}
were all members of the product space $\mathcal{V}^m_2$. Then \cref{prop:projection} would imply
\eq{
\Psi(\psi(q), v)(v_1^m) &= P^n_2 \left(v^m_1 \,|\, \psi(q), E^{-1}[\{(q, v)\} \cup {\hat{C}}] \right)
}
This would essentially
amount to establishing that $\varphi_2$ is honest, in the sense of \cref{def:honest}.

In the remainder of this paper we turn to two concrete applications of our framework,
by showing how it vindicates two heuristics that are common in both mathematics and the empirical sciences. As stated in the introduction,
%in these two applications we vindicate
the first of these heuristics is
giving more weight to hypotheses that are established via multiple lines of reasoning or pieces of evidence, and the
second heuristic is the common practice of abductive inference.\par

\section{Application I:\ Multiple Lines of Evidence}\label{sec:multiplelines}

%\dhwc{Make sure to use the appropriate notation for  sets in all spots below where that's what
%you're actually using.}

%In this section, we use our formalization of the calibration relationship between a mathematician-SMS and a
%far-future-community-of-mathematicians SMS to derive the benefit to a mathematician of a
A common
heuristic among real world mathematicians
is to ascribe greater
epistemic value to a proposition if there are multiple lines of evidence that favor it.
%and the benefit of using abductive reasoning.
%We then use the calibration relation between empirical
%scientists and physical universes to similarly derive the benefit to a scientist of
What is essentially the same heuristic appears among empirical scientists when they greater
epistemic value to a prediction about the outcome of an experiment if there are multiple lines of evidence that
support that prediction.

%\subsection{Multiple Lines of Evidence in General}
%In the practice of empirical science, if one has \textit{multiple} lines of evidence for some claim, i.e.,
%the claim is supported by multiple, independent experiments or observations,
%then one is typically thought to have \textit{better} evidence for that claim, compared to the case of only having one of those lines of evidence for the claim. Similarly, in mathematics, one might have {multiple} lines of evidence for some claim,
%i.e., the claim is supported by multiple, independent lines of formal reasoning.
%In this instance, one is also typically thought to have {better} evidence for that claim, compared to the case of only having one of those lines of evidence for the claim.

In formal epistemology, this heuristic is sometimes called the ``variety of evidence thesis,'' (see, among others, \citet{Wayne1995-WAYBAD}, \citet{Myrvold1996-MYRBAD}, \citet{Fitelson1996-FITWHA}, \citet{Bovens2002-BOVBNA}, \citet{Claveau2019-CLATVT-2}, and Landes (\citeyear{Landes2020-LANTVO-6}, \citeyear{Landes2020-LANVOE})).
While it has decided intuitive appeal as a principle of reasoning in both empirical
science and mathematics, it has resisted straightforward formal vindication.
In this section we show how our approach to modeling of mathematical and scientific reasoning using SMSs vindicates a version of this heuristic,
assuming that the appropriate calibration relationship holds.\par

\subsection{Multiple Lines of Evidence in Mathematics}
Consider the probability that a given mathematician-SMS $\varphi_{2}$ will produce some particular answer $v^\ast$ in response to a given question $q$. We view some particular set of claims, $B$, as being a ``line of evidence'' or ``evidence path'' for that answer to $q$ if the probability that the SMS will in fact produce $v^\ast$ in response to $q$ (rather than some alternative $v' \ne v^\ast$) is higher if
the probability is conditioned
on the SMS having already generated $B$, as compared to the probability of that answer without conditioning on $B$.
In practice, such an evidence path
will typically take the form of a ``proof'' (e.g., as published in the mathematical literature).
That is, it will consist of written material that expresses, with varying degrees of formal precision (and varying degrees of confidence by
author of the proof), a series of reasoning steps intended to establish a given conclusion.
Such a proof $B$ is a line of evidence that $v^*$ is in fact the correct answer to $q$ if the posterior probability that the far-future community of mathematicians will answer $v^*$ in response to question $q$ conditioned on $B$ is larger than the prior, unconditioned probability of that answer by the far-future community of mathematicians.

Here we formally justify this heuristic by proving that if an SMS $\varphi_{2}$ representing a mathematician is calibrated with an SMS $\varphi_{1}$ representing an oracle that determines the correct response distribution for each prediction pair $(q,\hat{C})$ (e.g., the far-future community of mathematicians), then the more evidence paths there are for a given claim $(q,v)$ according to $\varphi_{2}$, the more likely it is that
the oracle also responds to $q$ with the answer $v$.
This provides a normative justification for the practice of assigning higher probability to a given answer to a
mathematical question the more evidence paths there are for that answer.
% Here we
% %formally justify this heuristic by
% %proving
% prove that if we increase the number of such evidence paths $B$ in favor of $v^\ast$
% that we condition on, then we increase the probability of the answer $v^\ast$ in response to the question $q$.
% \dhwc{Uh, that synopsis is close to meaningless. Where does calibration come in, for example?}
% Note that in general, for a claim set to be an evidence path will require some ``prior knowledge'' set of claims, capturing all the foundational understanding of mathematics that are necessary for $B$ to even be meaningful. In the sequel we will incorporate that claim set into the analysis explicitly.\par

Define some arbitrary question and answer space $\mathcal{Q}, \mathcal{V}$, and let $\hat{\mathcal{C}}$ be the associated collection of all finite sets of claims $(q \in \mathcal{Q}, v \in \mathcal{V})$.
%(The terminology is only for convenience; \textit{a priori}, $\mathcal{Q}, \mathcal{V}$ and $\hat{\mathcal{C}}$ need not have anything to do with an SMS.)
Let $P$ be any distribution defined over $\hat{\mathcal{C}}$. (While we will mostly have in mind the case where
$P$ is the distribution of a mathematician-SMS, that is not required below.)
Define ${\mathcal{B}} = \{B(i) : 1,\ldots, n\}$ to be some indexed collection of $n$ claim sets
in $\hat{\mathcal{C}}$. (Note that we allow the different $B(i)$ to have nonzero intersection, i.e., to share claims.)
Let $\beta$ also be a claim set
in $\hat{\mathcal{C}}$, and let $q \in \mathcal{Q}$ be any question. Abusing notation, we will implicitly interpret a list of claim sets as their union, e.g., implicitly equating $P(B(1), \beta) = P(B(1) \cup \beta)$.
Note that in the usual way, since any claim $(q, v)$ is a claim set (a claim set having a single element), expressions like $P(v \,|\, q, \beta, B(1))$ are well-defined, as ratios of probabilities of claim sets. (See \cref{def:responsedist}.)

We need to formalize the idea that each $B(i)$ is a line of evidence for a particular answer $v^*$ (conditioned on $q$ and $\beta$). We also need to formalize the notion that the separate $B(i)$ do not ``confound one another'', i.e., that it is not the case that taken on their own, the separate $B(i)$ each increase the probability of $v^*$, but taken as group, they do not have that effect. Finally, we need to formalize the notion that $\beta$ is a set of ``foundational'' claims, which underlie the ability of each $B(i)$ to be an evidence path.

There are many ways to translate these desiderata
concerning ${\mathcal{B}}, \beta, q$ and $v^*$ into formal requirements. Perhaps most straight-forward is captured in the
following definition:
%would be to simply require that for all $i : n \ge i \ge 1$,
%\eq{
%P(v^\ast \,|\, q, \beta, B(1), \ldots, B(i)) &> P(v^\ast \,|\, q, \beta, B(1), \ldots, B(i-1))
%\label{eq:used_to_be_}
%}
%(where we interpret the $i=1$ version of this statement to mean $P(v^\ast \,|\, q, \beta, B(1)) \ge P(v^\ast \,|\, q, \beta)$).
\begin{Def}
Let  ${\mathcal{B}} = \{B(i) : 1,\ldots, N\}$ be some indexed collection of $N$ claim sets
in $\hat{\mathcal{C}}$,  $\beta$ another claim set
in $\hat{\mathcal{C}}$, and $q \in \mathcal{Q}$ be any question.
Suppose that for all $i : 1 \le i \le N$,
\eq{
\label{eq:3.2}
P(v^\ast \,|\, q, \beta, B(i)) > P(v^\ast \,|\, q, \beta)
}
Suppose further that for all those $i$,
\eq{
P(v^\ast \,|\, q, \beta, B(1), \ldots, B(i)) &> P(v^\ast \,|\, q, \beta, B(1), \ldots, B(i-1))
\label{eq:used_to_be_}
}
(where we interpret the $i=1$ version of this statement to mean $P(v^\ast \,|\, q, \beta, B(1)) \ge P(v^\ast \,|\, q, \beta)$).
Finally, suppose that all of the conditional probabilities in those equations are nonzero.
Then we refer to each of the sets $B(i)$
 as an \textbf{evidence path} (for the claim $(q, v^\ast)$ under the distribution $P$ and conditioned on the claim set $\beta$),
with ${\mathcal{B}}$ being an \textbf{evidence collection}.\footnote{
Below we will assume that a given collection of claim sets is an evidence collection,
i.e., that \cref{eq:used_to_be_} holds. It turns out that we can \textit{derive}  \cref{eq:used_to_be_} instead,
from another assumption, one that might seem
more innocuous.
That alternative is presented, along with the proof that it results in \cref{eq:used_to_be_},  in Appendix \ref{sec:proofs}.
}
\label{def:3.4}
\end{Def}
\noindent
(The precise indexing of the $N$ claim sets typically does not matter.)

At this point, despite our language,
 \cref{def:3.4} does not justify the heuristic of ascribing greater credence
to a claim if it is supported by multiple lines of evidence.
The problem is that a mathematician-SMS might be delusional, generating arbitrary claims  like `$1 + 1 = 3$', or
even claims which are not well-formed formulas, like `$1 + \int dx = sin(d)$', or other claims which have
no connection to ``mathematical truth'', however that concept is conceptualized. So it may be that evidence
paths combine synergistically as in \cref{eq:used_to_be_}
for that mathematician, but that this just means their delusional reasoning process has a certain property,
with no implication for mathematical truth. To put this another way, even though a mathematician's reasoning may be such that $\mathcal{B}$ is an evidence collection for some claim, that reasoning lacks any positive normative valence; there is nothing to ground the idea that the mathematician \textit{ought} to regard $\mathcal{B}$ as a set of lines of evidence for that claim.

%In order to secure this kind of positive normative valence for the mathematician's treating some set $\mathcal{B}$ as an evidence path for a claim, we follow our strategy of identifying
To address this issue, as discussed in \cref{sec:truth}, we identify mathematical ``truth''
with the limit distribution of a mathematical oracle,
%which as usual can be taken to be the far-future community of mathematicians,
and suppose
that the mathematician is calibrated with that oracle for some small $\epsilon$. We now show that this supposition
implies that if a given
collection of claim sets are evidence paths of the step-$n$ mathematician for a given question-answer
pair, then taken as a whole, they increase the probability of that answer to that question \textit{by
the oracle}.

%Recall \cref{eq:def_f}, and the associated definition of $m$.
First, note that calibration is defined only for a specific $q$ and $\hat{C}$.
However, to investigate the implications of \cref{eq:used_to_be_} for a mathematician-SMS who
has multiple evidence paths for a given answer to a question, and is making
predictions about the responses of an oracle, we need to
consider multiple pairs of $q$ and ${\hat{ \mathcal{C}}}$. This requires
an additional definition:

\begin{Def}
\label{def:prediction_distribution}
Suppose we are given a backward-consistent SMS $\varphi_1$, an SMS
$\varphi_2$, associated partial functions $\psi, \Psi$ as in \cref{def:prediction_pair}, and some $n \in \Z^+$.
Assume that $\psi$ is invertible over its domain of definition.
Define the associated event space of triples,
%\begin{linenomath}
\eq{
%\Lambda := \{(v^{m} \in \mathcal{V}^{m}_1, q \in \image(\psi), \hat{C}_2 \in \hat{\mathcal{C}}_2)\}
\Lambda := \{(v^{b} \in \mathcal{V}^{b}_1, q^b \in \mathcal{Q}_1^b, \hat{C}_2 \in \hat{\mathcal{C}}_2) : b \in \mathbb{Z}^+\}.
\nonumber
}
%\end{linenomath}
The \textbf{prediction distribution} $F^n_2$ is a distribution over $\Lambda$ defined by
\begin{enumerate}
\item $
F^n_2(v^{m} \,|\, \psi(q), {{\hat{C}_2}}) := \sum_{v' \in \mathcal{V}_2}  P^n_{2} (v' \;|\; q, {{\hat{C}_2}} ) \;
	 [\Psi(\psi(q), v')(v^m)]$
\noindent for all  $(q, \hat{C}_2)$ that are prediction pairs for $\varphi_2$ and $n$;
\item $F^n_2((\psi(q), \hat{C}_2) \propto P^n_{2} (q, {{\hat{C}_2}} )$ for
all  $(q, \hat{C}_2)$ that are prediction pairs for $\varphi_2$ and $n$;
%all $(q, \hat{C}_2)$ that are prediction pairs for $n$;
\item
$F^n_2((\psi(q), \hat{C}_2) = 0$ for all $(q, \hat{C}_2)$
that are not a prediction pair for $\varphi_2$ and $n$;
%for $n$;
\end{enumerate}
\label{def:4.3a}
where as usual $m$ is shorthand for $|\psi(q)|$, and it is assumed that there is at least one
$(q, \hat{C}_2)$ that is a prediction pair for $\varphi_2$ and $n$.
\label{def:F}
\end{Def}
\noindent The only role of $\varphi_1$ in \cref{def:prediction_distribution}
is indirect, via the requirement that $(q, \hat{C}_2)$
be a prediction pair. Accordingly we will often shorten $F^n_2$ to simply $F^n$.
%\noindent (Recall the definition of prediction pair, \cref{def:prediction_pair}.)

Note that $F^n$ is a fully specified distribution. For example, (3) tells us that if
$(q, \hat{C}_2)$ is not a prediction pair, then
$F^n(\{(\psi(q),v^{m})\} \cup \hat{C}_2) = 0$ for any $v^{m}$.
As another example, we can define conditional distributions
like $F^n(B(1) \,|\, \psi(q), \beta)$ for any two claim sets $B(1), \beta$ such that
$(q, \hat{C}_2 = B(1) \cup \beta)$ is a prediction pair.
% and such that
%$P^k_{2} (q, {{\hat{C}^k_2}} ) \ne 0$.
Note also that although defined {in terms of} an SMS, $F^n$ is not an SMS, nor a marginal distribution of an SMS.
Rather, loosely speaking, $F^{n}$ is the average prediction made by the mathematician for the answer of the
oracle to the questions specified in $\psi(q)$, i.e., it is the average prediction of the mathematician
for the answers to the questions that will be output by the ``true'' distribution.\footnote{One must be careful with
this informal interpretation of the prediction distribution though. For example, the proportionality constant
in \cref{def:4.3a}(2) does not equal $1$ in general, since $q$ and
$\psi(q)$ live in different dimensional spaces.}$^,$\footnote{We require $\psi$ to be invertible
to avoid the problem of how to define $F^n(\psi(q), \hat{C}_2)$ if there is some $q' \ne q$ such
that $\psi(q') = \psi(q)$, but while $(q, \hat{C}_2)$ is a prediction pair, $(q', \hat{C}_2)$ is not a prediction pair,
e.g., because it has zero probability of being generated by $\varphi_1$.}

These definitions allow us to state the following result:

\begin{prop}
Suppose we are given some $N, n \in \Z^+,  \epsilon \in \RRR^+$, backward-consistent SMS $\varphi_1$ and SMS
$\varphi_2$  where $\mathcal{\hat{C}}_{2}  \subseteq \mathcal{\hat{C}}_{1} $.
Suppose further that we are given associated partial functions $\psi, \Psi$ as in \cref{def:prediction_pair}.
Let $q$ be a question such that $|\psi(q)| = 1$. Fix some $v^* \in \mathcal{V}_{1}$ and
suppose that
\begin{enumerate}
\item $\mathcal{B} = \{B(i)\}$ is  a set of $N$ evidence paths under prediction distribution $F^n$ for $(\psi(q), v^*)$
conditioned on a claim set $\beta \in \hat{\mathcal{C}}^2$;

\item $(q, \hat{C}_2)$ is a prediction pair for $\varphi_2$ at step $n$, for any claim set
$\hat{C}_2  \subseteq \beta \cup \bigcup_i B(i)$;
\label{item:7.3.3}

\item $\varphi_2$ is calibrated with  $\varphi_1$ at step $n$ for $q, \epsilon$ and any
$\hat{C}_2  \subseteq \beta \cup \bigcup_i B(i)$, with divergence measure $D[., .]$;

\item ${{D}}[., .]$ is a locally Lipschitz continuous function of its
arguments.\footnote{Note that depending on the choice of the divergence ${{D}}[., .]$,
this Lipschitz continuity assumption might require that those distributions have full support.}
\end{enumerate}
Then for all $1 \le i \le N$, for small enough $\epsilon$,
%\begin{linenomath}
\eq{
\overline{P}_1(v^\ast \,|\, \psi(q), \beta, B(1), \ldots, B(i)) \;>\; \overline{P}_1(v^\ast \,|\, \psi(q), \beta, B(1), \ldots, B(i-1)).
\nonumber
}
%\end{linenomath}
\label{prop:3}
\end{prop}
\noindent (The restriction that $|\psi(q)| = 1$ in \cref{prop:3} is simply for convenience, to avoid
introducing extra notation.)

Intuitively, \cref{prop:3} means that if the mathematician is calibrated with an oracle,
then the greater the number of evidence paths for the claim $(\psi(q),v^{*})$ under the distribution $F^n$
--- the average prediction made by that mathematician for the response distribution of the oracle  ---
the more likely it is that the oracle actually responds with the answer $v^{*}$
to the question $\psi(q)$.
% \dhwc{Clean up this paragraph. E.g., the mentioned discussion cannot
% refer to $F^n($ --- it hadn't been introduced yet. Also, is the ending sentence simply the assumption
% that the mathematician and the far-future community are both honest, for the same partial functions,
% and that the mathematician is perfectly calibrated with that community (or something like that)?}
%To illustrate this, recall from the interpretation the limit distribution of the SMS $\varphi_{1}$ as the limit output of the far-future community of mathematicians. Under this interpretation, $F^{n}$ is the average prediction made by the mathematician for the answer of the
%oracle to the question specified in $\psi(q)$, i.e., for the average prediction of the mathematician
%for what the ``true'' answer to the question is. %Moreover, assuming the mathematician
%% responds to any question they consider with the answer that they think the oracle
%% would give to that question, $F^n($ (evaluated for any prediction pair $(q, \hat{C})$)
%% is just the response distribution of that mathematician themself.
%%Finally, as emphasized above, we could in fact
%%use \textit{any} oracle $\varphi_1$ to define mathematical truth --- our use of the far-future community
%%of mathematicians is simply an intuitive convenience.
In this sense, if a calibrated mathematician-SMS follows the
multiple lines of evidence heuristic in their reasoning, assigning greater probability to a proposition if there
are more evidence paths supporting it, then the far-future community of mathematicians
(or, more generally, any arbiter of mathematical truth) also assigns greater probability to that proposition.

\subsection{Multiple Lines of Evidence in Science}

We can analyze the potential benefit of multiple lines of evidence to a scientist using
essentially the same approach we used for analyzing the potential benefit
of multiple lines of evidence to a mathematician. We simply need to consider a scientist-SMS
embedded in a universe-SMS rather than a mathematician-SMS predicting a far-future
community of mathematicians SMS.

The formal definition of an evidence path for the case of scientist-SMSs embedded in a universe-SMS
is exactly the same as in the case of mathematics, given in \cref{def:3.4}. However, the interpretation of an evidence path is different.
In a scientific context, we conceive of evidence paths not as chains of proof-like mathematical
reasoning, but instead as chains of (states of the brain of the scientist that correspond
to the) outcomes of experiments or observations the scientist has made. Similarly, the formal
definitions of prediction pairs and of the distribution $F^n$ are unchanged, but
the interpretation of $F^n$ is different. In the case of mathematicians,
it is the (expected value of) the prediction of the mathematician for the
outcomes generated by the far-future community (or whatever the oracle is) to the mathematical questions specified in $q$.
%but our interpretations of the two associated SMSs $\varphi_1$ and  $\varphi_2$ are different.
In the case of a scientist, it is instead the (expected value of) the prediction of the scientist for the
outcomes generated by the universe to the experiments specified in $q$.
%but our interpretations of the two associated SMSs $\varphi_1$ and  $\varphi_2$ are different.

There is also one concrete, non-interpretational distinction between the case of mathematicians
predicting the answers of their far-future community to mathematical questions and scientists predicting
the outcomes of physical experiments. In
%the first case, calibration means that the far-future
%community of mathematicians shares all of the conclusions of the current mathematician
%relevant to the question at hand, given by $\hat{C}_2$, i.e., given by the chains of proof-like
%mathematical reasoning. In the case of scientists though,  $\hat{C}_2$ is a physically concrete property
%of the scientist's own brain, corresponding to their mental perception of
%the outcomes of chains of experiments or observations they have made. Moreover,
in
the case of mathematicians and calibration, the quantity that we require to have been
considered by the far-future community of mathematicians, the  claim set,
is identical to the chains of proof-like mathematical reasoning considered by the mathematician. In contrast,
the associated quantity arising in the assumption of embedding
calibration is $E^{-1}([\{(q,v)\} \cup \hat{C}_2])$,
the collection of all universe-SMS claim sets that are consistent both with the scientist's mental perception of
the outcomes of chains of experiments or observations they have made, and
with their subsequent prediction about the probability distribution over the outcomes of the experiments in $\psi(q)$.

These distinctions yield a result that is similar to \cref{prop:3}, but not quite identical. To state this new result,
define an \textbf{embedded prediction distribution} exactly as a prediction distribution is defined, in \cref{def:F},
with two changes. First, an embedding function $E$ is specified.
Second, wherever the pair $(q, \hat{{C}})$ arises in \cref{def:F}, it is now taken to be an
embedded prediction pair for embedding function $E$, rather than a prediction pair.

\begin{prop}\label{prop:3s}
Suppose we are given some $N, n \in \Z^+, \epsilon \ge 0$, backward-consistent SMS $\varphi_1$ and SMS
$\varphi_2$ with embedding function $E$.
Suppose further that we are given associated partial functions $\psi, \Psi$ as in \cref{def:prediction_pair}, and
a question $q$ such that $|\psi(q)| = 1$. Fix some $v^* \in \mathcal{V}_{1}$ and
suppose that
\begin{enumerate}
\item $\mathcal{B} = \{B(i)\}$ is  an $N$-element evidence collection under embedded prediction distribution $F^n($ for $(\psi(q), v^*)$ and
the embedding function $E$, conditioned on a claim set $\beta \in \hat{\mathcal{C}}^2$;

\item $(q, \hat{C}_2)$ is an embedded prediction pair for $\varphi_2$ at step $n$ for embedding function $E$, for any claim set
$\hat{C}_2  \subseteq \beta \cup \bigcup_i B(i)$;

\item Given the $\hat{C}_2$ in condition (2), $\varphi_2$ is embed-calibrated with the SMS $\varphi_1$ at step $n$ for $\epsilon$ and
$(q, \hat{C}_2)$, with divergence measure $D[., .]$ and embedding function $E$;

\item ${{D}}[., .]$ is a locally Lipschitz continuous function of its arguments.
%\footnote{Note that depending on the choice of the divergence ${{D}}[., .]$,
%this Lipschitz continuity assumption might require that those distributions have full support.}
\end{enumerate}
Then for all $1 \le i \le N$, for small enough $\epsilon$,
%\begin{linenomath}
\begin{multline}
    \sum_{v\in\mathcal{V}_{2}}P^{n}_{2}(v|q,\beta, B(1), \ldots, B(i))\overline{P}_1
	\left(v^\ast \,|\, \psi(q), E^{-1}[\{(q,v\}\cup\beta, B(1), \ldots, B(i))]\right) \\ > \ \sum_{v\in\mathcal{V}_{2}}P^{n}_{2}
	\left(v|q,\beta, B(1), \ldots, B(i-1))\overline{P}_1(v^\ast \,|\, \psi(q), E^{-1}[\{(q,v\}\cup\beta, B(1), \ldots, B(i-1)]\right).
\label{eq:prop.7.4.equation}
\end{multline}
%\end{linenomath}
\label{prop:4}
\end{prop}
\noindent (As with \cref{prop:3}, the restriction in \cref{prop:4} that $|\psi(q)| = 1$ is for convenience.)
%
%\dhwc{Troubling fact: We never (?) use the definition of ``embedded'', in \cref{def:embedded}. \textit{Any} function
%$E$ with the appropriate domain and range would work. If that's true, then ... I'm not sure what.}

Suppose that a scientist embedded in a universe-SMS
views $\mathcal{B}$ as an evidence collection for some claim $(\psi(q),v^\ast)$, and that
this scientist is embed-calibrated with the universe-SMS. Loosely speaking, \cref{prop:3s} says that
under these conditions,
%the universe-SMS is such that
the more evidence paths
% in $\mathcal{B}$ that this
the scientist observes, the higher (the expected value of) the probability that \textit{the universe-SMS itself} outputs $v^{*}$ as the answer to the question $\psi(q)$.\par

There is an important difference between the conclusion of \cref{prop:3} and the conclusion of \cref{prop:3s}. \cref{prop:3} establishes that
the probability that the mathematical oracle (e.g., the far-future community of mathematicians) $\varphi_1$
outputs $v^\ast$ as the answer to $\psi(q)$
grows with the number of evidence paths that the present-day mathematician $\varphi_2$ observes for the claim $(\psi(q),v^\ast)$.
%
%if a calibrated mathematician-SMS $\varphi_{2}$ treats $\mathcal{B}$ as an evidence collection for some claim $(\psi(q),v^\ast)$, then so does the oracle-SMS $\varphi_{1}$ with which $\varphi_{2}$ is calibrated.
By contrast, \cref{prop:3s} shows that the expected probability that the universe-SMS $\varphi_1$
outputs $v^\ast$ as the answer to $\psi(q)$, where that expectation is taken across all possible predictions that the scientist-SMS
$\varphi_2$ might make about the probability distribution over answers to $\psi(q)$,
grows with the number of evidence paths that the scientist-SMS observes for the claim $(\psi(q),v^\ast)$.

Ultimately, the reason for this difference is that
in the case of empirical science there is only a single SMS  (the physical universe) and
a projection of that SMS, which can be interpreted as another,
dependent SMS (the scientist). In contrast, in the case of mathematics, there are two SMSs
which are allowed \textit{a priori} to be completely independent (e.g., a current mathematician and
an oracle, like the far-future community of mathematicians).
Concretely, in the case of a scientist, the physical fact that that scientist's brain is in a state that
involves asking question $q$ and giving answer $v$ restricts the set
of possible physical universes, over and beyond the restriction given by the fact that the state of the scientist's brain also
involves a claim set $\beta, B(1), \ldots$.
(This is why the claim $(q, v)$ is part of the argument of the inverse embedding function in the definition of embed-calibration.)
That in turn means that in the definition of embed-calibration,
the average of $D(., .)$ over all answers $v$ by the scientist affects both of $D(., .)$'s arguments. In contrast, in the case of a
mathematician-SMS predicting the answers of an oracle-SMS, we don't assume that the oracle-SMS generates the claim $(q, v)$
made by the mathematician-SMS. So the average of $D(., .)$ only affects the first of $D(., .)$'s arguments.

%Focusing just on the precise form of \cref{prop:3s},
Another important point is that
one might think that the supposition that  $\mathcal{B}$ is an evidence collection already says
all that's important, from a normative perspective, with none of the other assumptions needed.
Recall though that as discussed above, $\mathcal{B}$ could be an evidence collection simply because the scientist
is prone to hallucinations about the results of experiments and
is self-consistent in their hallucinations.
To rule out this possibility of self-consistent hallucination we need to couple
what the scientist thinks they are likely to
observe in the future with the actual future states of the physical reality that they're embedded in.
The third supposition in \cref{prop:3s}, that the scientist is also embed-calibrated, is one way to enforce this coupling.
The result of this coupling is made explicit in \cref{eq:prop.7.4.equation}, which involves both the
distribution of the scientist and of the universe.

Even given this coupling though, so that the universe-SMS produces answer $v^*$ to the experiments specified in $\psi(q)$,
one might worry that the scientist-SMS observing the universe might not be able to observe the results of such experiments.
Or that even if they did observe those results, they might do so in a highly biased manner, being led to believe that the
the results of the experiments was \textit{not} $v^*$, even if it was. This possibility can be addressed by adding the
requirement that the combination of each of the evidence path with the question $q$ is a discriminating pair (see \cref{def:distinguisher}).
%
%is directly ruled out in \cref{prop:3s}
%though, by the supposition that  $\mathcal{B}$ is an evidence collection under distribution $F^n($.

There are other ways to couple what the scientist thinks they are likely to
observe in the future with the actual future states of the physical reality that they're embedded in.
%However, the
%fact that (by assumption) $\mathcal{B}$ is an evidence collection directly establishes
%
%
%For an evidence collection
%of a scientist (comprising a set of their own previous observations) to have a normative role, that collection
%must increase the probability that \textit{the results of their own future
%observations} will have the associated particular values.
%This is the precise property established by the supposition that
%%, then that implies that those observations will in fact have those values.
One way to do this is to flip the roles of the scientist's distribution and
the universe's distribution in the suppositions in \cref{prop:3s}, as follows:

\begin{prop}
Suppose we are given some $N, n \in \Z^+, \epsilon \ge 0$, backward-consistent SMS $\varphi_1$ and SMS
$\varphi_2$ with embedding function $E$.
Suppose further that we are given associated partial functions $\psi, \Psi$ as in \cref{def:prediction_pair}, and
a question $q$ such that $|\psi(q)| = 1$. Fix some $v^* \in \mathcal{V}_{1}$ and
suppose that
\begin{enumerate}
\item $\mathcal{B} = \{B(i)\}$ is an $N$-element evidence collection under distribution $\overline{P}_{1}$ for $(q, v^*)$
conditioned on a claim set $\beta$;

\item $(q, E[\hat{C}_1])$ is an embedding prediction pair for $\varphi_2$ at step $n$ and any claim set
$\hat{C}_1\subseteq\beta \cup \bigcup_i B(i)$;

\item Given the $\hat{{C}}_1$ in condition (2), $\varphi_2$ is embed-calibrated with the universe-SMS $\varphi_1$ at step $n$ for
$(q, E[\hat{C}_1])$, with divergence measure $D[., .]$;
% \subseteq \beta \cup \bigcup_i B(i)$;

\item ${{D}}[., .]$ is a locally Lipschitz continuous function of its
%probability distribution
arguments;
%(where those distributions are considered as vectors in a Euclidean metric
%space);

\item for all $1 \le i \le N$ and all $v\in\mathcal{V}_{2}$, $$\overline{P}_{1}(v^{*}|\psi(q),E^{-1}[\{(q,v)\}\cup E[\beta,B(1),\dots,B(i)]])\propto \overline{P}_{1}(v^{*}|\psi(q),\beta,B(1),\dots,B(i));$$
\end{enumerate}

Then for all $1 \le i \le N$, $v\in\mathcal{V}_{1}$, for small enough $\epsilon$, the embedded prediction distribution  obeys
%\begin{linenomath}
\eq{F^{n}_2(v^{*}|q,E[\beta,B(1),\dots,B(i)])
	\;>\; F^{n}_2(v^{*}|q,E[\beta,B(1),\dots,B(i-1)]).
	}
%\end{linenomath}
%where $F^{n}$ is the embedded prediction distribution of $\varphi_{2}$.
\label{prop:3b}
\end{prop}

%embedding function $E$.

%Formally, expand the conditional distribution on the LHS of \cref{eq:4.8} and
%then use \cref{eq:6.1a} to get
%\eq{
%\dfrac{\overline{P}_1((v^\ast, \psi(q)), E^{-1}[\beta, B(1), \ldots, B(i))}
%	{\overline{P}_1(\psi(q), E^{-1}[\beta, B(1), \ldots, B(i)])}
%	&=
%\dfrac{\overline{P}_1((v^\ast, \psi(q)), E^{-1}[\beta, B(1), \ldots, B(i)])}
%	{\sum_{v_1 \in {\mathcal{V}}_1} \overline{P}_1((v_1, \psi(q)), E^{-1}[\beta, B(1), \ldots, B(i)])} \\
%	&=
%\dfrac{P_2^{n+1}((v^\ast, \psi(q)), \beta, B(1), \ldots, B(i))}
%		{\sum_{v_2 \in {\mathcal{V}}_2} P_2^{n+1}((v_2, \psi(q)), \beta, B(1), \ldots, B(i))} \\
%	&= P_2^{n+1}((v^\ast \;|\;  \psi(q)), \beta, B(1), \ldots, B(i))
%}
%Expanding the conditional distribution on the RHS of \cref{eq:4.8} the same way and plugging in
%both of our expansions to  \cref{eq:4.8} gives
%%a ratio of probabilities, and then apply
%% \cref{eq:6.1a} to both of those probabilties. This transforms that condidtional distribution to
%%applying \cref{eq:6.1a} twice in \cref{eq:6.8}, once to the distribution on each side, we get
%%%\begin{linenomath}
%%\eq{
%%P_2^{n+1}((v^\ast , \psi(q)), \beta, B(1), \ldots, B(i)) &> P_2^{n+1}((v^\ast, \psi(q)), \beta, B(1), \ldots, B(i-1)),
%%}
%%%\end{linenomath}
%%and therefore
%%\begin{linenomath}
%\eq{
%P_2^{n+1}(v^\ast \,|\, \psi(q), \beta, B(1), \ldots, B(i)) &> P_2^{n+1}(v^\ast \,|\, \psi(q), \beta, B(1), \ldots, B(i-1)).
%}
%%\end{linenomath}
\noindent
\cref{prop:3b} establishes that in a universe-SMS where $\mathcal{B}$ is an evidence collection for the claim $(q,v^{*})$, a scientist-SMS that is
both embedded in that universe and calibrated with it should, in expectation, adopt a greater degree of belief that $v^\ast$ is the answer to $q$. Note that the fifth condition, which has no analog in \cref{prop:3s}, states that the probability that the universe-SMS responds to question $\psi(q)$ with answer $v^{*}$, conditional on any evidence collection $\mathcal{B}$ and a claim set $\beta$, is proportional to the probability that the universe-SMS responds to question $\psi(q)$ with answer $v^{*}$ conditional on the claim set $E^{-1}[\{(q,v)\}\cup E[\beta,B(1),\dots,B(i)]]$ (where $E^{-1}[\{(q,v)\}\cup E[\beta,B(1),\dots,B(i)]]$ denotes the universe-SMS description of the scientist-SMS observing the events $\beta,\beta,B(1),\dots,B(i)$ and then outputting the answer $v$ to the question $q$). Concretely, we can say that probability that the universe-SMS outputs the answer $v^{*}$ to the question $\psi(q)$, conditional on $\beta,B(1),\dots,B(i)$ is not meaningfully effected by either the fact the the scientist-SMS observes the events $\beta,B(1),\dots,B(i)$ or the fact that scientist-SMS outputs an answer $v$ to the question $q$.\par

\section{Application II:\ Abduction}\label{sec:abduction}

\subsection{Abduction in General}
%Under Frankfurt's (\citeyear{frankfurt1958peirce}) exegesis of (\citet{peirce1960collected}), an `abductive inference'
%is any inference with the following pattern:
%\begin{enumerate}
%    \item A ``surprising fact'' $S$ is observed.
%
%    \item If $R$ were true, then $S$ would be true as a matter of course.
%
%    \item Therefore, there is reason to suspect that $R$ is true.
%\end{enumerate}
%\noindent
%Note that this differs slightly from more modern uses of `abduction' as a synonym for `inference to the best explanation' (see Douven \citeyear{sep-abduction}). In the classic Peircean definition of abduction, there is no implicit or explicit ranking of explanations $A$ for the surprising fact $C$, such that we can infer a reason to suspect that $A$ is true on the ground that $A$ best explains $C$. Rather, on this definition abduction is simply the inference to a reason to suspect that $A$ is true in virtue of $A$'s rendering $C$ less surprising, regardless of whatever other facts might also render $C$ less surprising. \citet{viteri2020explosive} argue that, descriptively speaking, Peircean abduction has played a significant role in mathematical progress throughout its history.\par
%
Following modern uses of `abduction' as a synonym for `inference to the best explanation' (see Douven \citeyear{sep-abduction}),
%here we define abduction so that it contains explicitly probabilistic language, and is put in terms of probabilistic relations between questions and answers, rather than facts. That is,
we take abduction to be the following inference pattern:
\begin{enumerate}
    \item The question $q^{*}$ has answer $v^{*}$ with probability $x^{*}$.

    \item The question $q^{\dagger}$ has answer $v^{\dagger}$ with probability $x^{\dagger}$.

    \item If $q^{\dagger}$ does indeed have answer $v^{\dagger}$, then $q^{*}$ has answer $v^{*}$ with
% probability would have
 probability $y^{*}>x^{*}$.

    \item It turns out that $q^{*}$ does have answer $v^{*}$.

    \item We conclude that $q^{\dagger}$ has answer $v^{\dagger}$ with a probability $y^{\dagger}>x^{\dagger}$.
\end{enumerate}
%While slightly more detailed, this inference pattern preserves the basic structure of Peircean abduction. The claim $(q^{*},v^{*})$ is taken to be surprising, to the extent that probability $x^{*}$ is low. However, despite being surprising, it is supposed that $q^{*}$ does have answer $v^{*}$. Moreover, it is supposed in addition that if the question $q^{\dagger}$ is assigned answer $v^{\dagger}$, then this renders the claim $(q^{*},v^{*})$ less surprising. From this, we infer that the claim $(q^{\dagger},v^{\dagger})$ is also made less surprising by the claim $(q^{*},v^{*})$. This mirrors the way in which, in the original presentation of Peircean abduction, the occurrence of $C$ licenses an inference to the more likely truth of $A$ because $A$ renders $C$ less surprising.\par
%
In particular, if we take $y^* = 1$, then (3) means that the answer $v^\dagger$ to the
question $q^\dagger$ implies that $v^*$ is the answer to the question $q^*$. The conclusion (5) then
says that in light of this premise, if in fact $v^*$ is the answer to the question $q^*$,
then it is more likely that $v^\dagger$ is the indeed the answer to the question $q^\dagger$.

In what follows we prove that
% in both mathematical and scientific contexts,
% the norms of reasoning that we have presented above lead to the conclusion that
both mathematicians and scientists ought to engage in abductive reasoning, under certain assumptions.
This is a partial justification of using abduction as a heuristic in human reasoning.\footnote{See \citet{viteri2022epistemic} for a descriptive study of the role of abduction in mathematical reasoning.}

\subsection{Abduction in Mathematics}
In this section, we show that if a mathematical oracle treats one claim as abductively supporting another, then an individual mathematician should follow the same inference pattern. To begin, we let $P(.)$ be an arbitrary joint probability distribution defined over $\mathcal{C} \times \mathcal{C}
\times \hat{\mathcal{C}}$ for some arbitrary claim space $\mathcal{C}$. (At this point, $P(.)$ has nothing to do with SMSs of any sort.) For simplicity, suppose that under $P$ there is probability $1$ that the
two claims have questions $q^\ast$ and $q^\dagger$, respectively. Suppose as well that if the answer $v^\dagger$  is generated in response to the question $q^\dagger$, that would make it more likely that the answer $v^\ast$ would be generated in response to the question $q^\ast$,  conditioned on a
particular claim set $\hat{{C}}$. Formally, this condition means that for some $\alpha > 1$,
%\begin{linenomath}
\eq{
P\left(v^\ast\;|\; q^\ast, (q^\dagger, v^\dagger),\hat{C}\right) = \alpha P\left(v^\ast\;|\; q^\ast,\hat{C} \right),
\label{eq:5}
}
%\end{linenomath}
or equivalently,
%\begin{linenomath}
\eq{
\dfrac{ P\left( (q^\ast, v^\ast), (q^\dagger, v^\dagger),\hat{C} \right)}
				{P\left( q^\ast, (q^\dagger, v^\dagger),\hat{C} \right)}
    &= \alpha \dfrac{ P\left((q^\ast, v^\ast\right),\hat{C})} {P(q^\ast,\hat{C})},
}
%\end{linenomath}
and so repeatedly using our assumption that both $q^\ast$ and $q^\dagger$ occur with probability $1$,
%\begin{linenomath}
\eq{
\dfrac{P\left((q^\ast, v^\ast), (q^\dagger, v^\dagger),\hat{C}\right)}
				{P\left((q^\dagger, v^\dagger),\hat{C} \right)}
    &= \alpha  \dfrac{ P\left((q^\ast, v^\ast),\hat{C}  \right)}{P(\hat{C})} \\
\dfrac{P\left((q^\ast, v^\ast), (q^\dagger, v^\dagger),\hat{C}\right)}
						{P\left((q^\ast, v^\ast),\hat{C} \right)}
    &= \alpha  \dfrac{ P\left((q^\dagger, v^\dagger),\hat{C}  \right)}{P(\hat{C})} \\
\dfrac{P\left( (q^\ast, v^\ast), (q^\dagger, v^\dagger),\hat{C}\right)}
			{P\left(q^\dagger, (q^\ast, v^\ast),\hat{C} \right)}
    &= \alpha \dfrac{P\left((q^\dagger, v^\dagger),\hat{C}\right)} {P(q^\dagger,\hat{C})}
}
%\end{linenomath}
i.e.,
%\begin{linenomath}
\eq{
    P\left(v^\dagger \;|\; q^\dagger, (q^\ast, v^\ast),\hat{C} \right) = \alpha
					P(v^\dagger \;|\; q^\dagger,\hat{C}).
\label{eq:6.7}
}
%\end{linenomath}
In the sequel we will refer to \cref{eq:5} as the ``abduction premise'' and to \cref{eq:6.7} as the ``abduction implication'',
and simply say that ``abduction holds'' when the former implies the latter.
The abduction premise
says that if the answer $v^\dagger$ is generated in response to the question $q^\dagger$, that would make it more likely that the answer $v^\ast$ would be generated in response to the question $q^\ast$. We have just shown that this premise implies its converse:\ if the answer $v^\ast$ were generated in response to the question $q^\ast$,  that would make it more likely that the answer $v^\dagger$ would be generated in response to the question $q^\dagger$. So we have established that abduction holds for the distribution $P(.)$, for the two question-answer pairs $(q^\ast, v^\ast), (q^\dagger, v^\dagger)$.

In particular, consider the case where
the conditional distributions $P(. \,|\, .)$ in \cref{eq:5} are $P^k_2(. \,|\, ., \hat{C}^k_2)$ for the step $k$ of a mathematician-SMS $\varphi_2$, where as usual $\hat{C}^k_2$ is the  claim set associated with a step-$k$ claim set that they generate with nonzero probability, $\hat{C}^k_2$. Then we have established that abduction holds for that mathematician, at step $k$, simply by the laws of probability theory.
Similarly, we have established that abduction holds for the (limit distribution of the) oracle SMS, e.g., for
the far-future community of mathematicians, again simply by the laws of probability theory.

However, as in the case of establishing the evidential value of multiple proof paths, we have made no assumptions at this stage that the output of the mathematician-SMS is in any way connected to mathematical truth. So, while we have stated some conditions under which a mathematician-SMS will reason abductively, we have not yet established that the same mathematician-SMS \textit{ought} to reason abductively.

To address this issue, recall
the prediction distribution $F^{n}$ defined in \cref{def:4.3a}, which is the structure that
captures what an SMS mathematician thinks the distribution of answers of an oracle-universe
would be to a given set of questions. Specifically,
${\PPP}(v^\ast, v^\dagger \;|\; q^\ast, q^\dagger,\hat{C}_{2})$ is a distribution
over joint answers $(v^\ast, v^\dagger)$ by the oracle-universe to the pair of questions $(q^\ast, q^\dagger)$,
constructed by averaging the ``predictions'' for what that distribution is made by the mathematician
SMS $\varphi_2$ at step $n$ (via the function $\psi(q)$). (All of these distributions are conditioned on the
mathematician having also generated the claims listed in $\hat{C}_{2}$ at step $n$.)
Using $\PPP$ and our notational convention \cref{eq:3.44a}, we can establish the
following  sufficient conditions for an SMS $\varphi_{2}$ that is calibrated with an oracle-SMS $\varphi_{1}$ to be justified in abductive reasoning:

\begin{prop}\label{prop:mathabduction}
Suppose we are given some $N, n \in \Z^+$, backward-consistent SMS $\varphi_1$, SMS
$\varphi_2$, and associated partial functions $\psi, \Psi$ as in \cref{def:prediction_pair}, where
$\psi$ is an invertible function. Let  $v^*, v^\dagger \in \mathcal{V}_1$, and let $q^*, q^\dagger \in \mathcal{Q}_1$ be
two questions such that all three tuples $q^*$, $q^\dagger$ and $(q^*, q^\dagger)$ are in the codomain of $\psi$.
Suppose as well that
\begin{enumerate}
    \item
 $\PPP$ satisfies the abductive premise, \cref{eq:5}, for all $\hat{C}_{2}\in\mathcal{\hat{C}}_{2}$ and
for the specific questions $q^*$, $q^\dagger$ and answers $v^{*},v^\dagger$;

    \item $(\psi^{-1}(q^{*}),\hat{C}_{2})$, $(\psi^{-1}(q^{\dagger}),\hat{C}_{2})$, and $(\psi^{-1}(q^{*},q^{\dagger}),\hat{C}_{2})$ are all prediction pairs for $\varphi_{2}$ at step $n$ for any $\hat{C}_{2}\in\hat{\mathcal{C}}_{2}$;

    \item
%$F^{n} \left(v^\ast \,|\, q^\ast, q^\dagger, \hat{C}_{2}\right)
$\sum_{v^\dagger} F^{n} \left(v^\ast, v^\dagger \,|\, q^\ast, q^\dagger, \hat{C}_{2}\right)
= F^{n} \left(v^\ast \,|\, q^\ast, \hat{C}_{2}\right)$;

    \item $\varphi_2$ is calibrated with $\varphi_{1}$ at step $n$ for $\hat{C}_{2}$ for each of the three questions $\psi^{-1}(q^{*})$, $\psi^{-1}(q^{\dagger})$, and $\psi^{-1}(q^{*},q^{\dagger})$, and for $m=1,2$;
%, where $q^\ast, q^\dagger \in \mathcal{Q}_1$;

    \item ${{D}}[., .]$ is a locally Lipschitz continuous function of its probability distribution arguments (where those distributions are considered as vectors in a Euclidean metric space);

    \item
%$\overline{P}_1 \left(v^\ast \,|\, q^\ast, q^\dagger, \hat{C}_{2}\right)
$\sum_{v^\dagger} \overline{P}_1 \left(v^\ast, v^\dagger\,|\, q^\ast, q^\dagger, \hat{C}_{2}\right)
= \overline{P}_1 \left(v^\ast \,|\, q^\ast, \hat{C}_{2}\right)$.
\end{enumerate}
Then $\overline{P}_{1}(v^{\dagger}|q^{\dagger},(q^{*},v^{*}), \hat{C}_{2})>\overline{P}_{1}(v^{\dagger}|q^{\dagger}, \hat{C}_{2})$.
\end{prop}
\noindent The requirement that $\psi$ be invertible is for convenience. Note that
since by condition (2) there are questions $q_2 \in \mathcal{Q}_2$
such that  $(q_2,\hat{C}_{2})$ is a prediction pair, $\mathcal{\hat{C}}_{2}  \subseteq \mathcal{\hat{C}}_{1} $.
%In condition (3), we take $q^\dagger$ to be an extra question for $\varphi_1$ rather than the question part of
%a claim in the claim set of $\varphi_2$, i.e., we take $F^{n} \left(v^\ast \,|\, q^\ast, q^\dagger, \hat{C}_{2}\right)$
%as shorthand for $\sum_{v^\dagger} F^{n} \left(v^\ast, v^\dagger \,|\, q^\ast, q^\dagger, \hat{C}_{2}\right) $
Intuitively, condition (3) means that
the extra information that $\varphi_{1}$ produces the question $q^\dagger$ doesn't affect the expectation of the conditional probability that the mathematician assigns to the event that the oracle-SMS
produces the answer $v^\ast$ to the question $q^\ast$ (given the information already contained in $\hat{C}_{2}$).
Similarly, condition (5) means that the extra information that the limit distribution of $\varphi_{1}$ generates the question $q^\dagger$ doesn't affect the conditional probability that it produces the answer $v^\ast$ to the question $q^\ast$ (given the information already contained in $\hat{C}_{2}$).

\cref{prop:mathabduction} says that if the distribution $\PPP$ representing the \textit{beliefs} of a mathematician-SMS $\varphi_{2}$ about the output of an oracle-SMS $\varphi_{1}$ satisfies the abductive premise
for $q^*, q^\dagger, v^*$ and $v^\dagger$,
and if $\varphi_{2}$ is appropriately calibrated with $\varphi_{1}$,
then the oracle-SMS must obey the associated abduction implication.
More concretely, note that we have specifically assumed that the mathematician-SMS is calibrated with the oracle-SMS for the questions $\psi^{-1}(q^{*})$, $\psi^{-1}(q^{\dagger})$, and $\psi^{-1}(q^{*},q^{\dagger})$.
These are questions for the mathematician, which can be
interpreted as asking them,  `what is the probability distribution over possible answer(s) by the
oracle-SMS  to the questions $q^{*}$/$q^{\dagger}$/$(q^{*},q^{\dagger})$ (resp.)?'. So a mathematician who is calibrated with the oracle with respect to the probabilities that the oracle assigns to the answers to the questions $q^{*}$ and $q^{\dagger}$ will also follow the abductive inference pattern when answering these questions. To the extent that we take it to be a normative goal of a mathematician or finite group of mathematicians that their predictions about the answers to questions match those of the far-future mathematical community, it stands to reason that said mathematicians ought to make use of abduction in at least some circumstances.\par

\subsection{Abduction in Science}

% \dhwc{Don't we need to redefine a tweak to the definition of $G$, appropriate for the case of scientific reasoning,
% just like we needed to define a tweak to the definition of prediction distribution for the case of scientific reasoning?}
As when we considered the benefit of multiple evidence paths for empirical scientific reasoning rather than for
mathematical reasoning, we begin by re-defining $F^{n}$ as an embedded prediction distribution, rather than a prediction distribution \textit{simpliciter}. This allows us to state a result
analogous to \cref{prop:mathabduction}, for the case of scientific rather than mathematical reasoning.
% Our vindication of abduction in science proceeds similarly to our vindication of abduction in mathematics. As we did in our the case of multiple evidence paths, we begin with an analogous proposition to \cref{prop:mathabduction}:
% \dhwc{Do we just want to say, ``We begin with an analogous proposition to \cref{prop:mathabduction}?}
\begin{prop}\label{prop:scienceabduction2}
Suppose we are given some $N, n \in \Z^+$, backward-consistent SMS $\varphi_1$, SMS
$\varphi_2$, and associated partial functions $\psi, \Psi$ as in \cref{def:prediction_pair}.
Let  $v^*, v^\dagger \in \mathcal{V}_1$, and let $q^*, q^\dagger \in \mathcal{Q}_1$ be
two questions such that all three tuples $q^*$, $q^\dagger$ and $(q^*, q^\dagger)$ are in the codomain of $\psi$.
Let $E$ be an embedding function. Suppose as well that

\begin{enumerate}
\item  $\PPP$ satisfies the abductive premise, \cref{eq:5}, for all $\hat{C}_{2}\in\mathcal{\hat{C}}_{2}$ and
for the specific questions $q^*$, $q^\dagger$ and answers $v^{*},v^\dagger$;

%    \item The prediction distribution
% $\PPP$ satisfies the abductive premise, \cref{eq:5}, for all $\hat{C}_{2}\in\mathcal{\hat{C}}_{2}$ and
%for the specific questions $q^*$, $q^\dagger$ and answers $v^{*},v^\dagger$;

    \item $(\psi^{-1}(q^{*}),\hat{C}_{2})$, $(\psi^{-1}(q^{\dagger}),\hat{C}_{2})$, and $(\psi^{-1}(q^{*},q^{\dagger}),\hat{C}_{2})$ are all embedded prediction pairs for $\varphi_{2}$ at step $n$ for any $\hat{C}_{2}\in\hat{\mathcal{C}}_{2}$;

    \item $\sum_{v^\dagger} F^{n} \left(v^\ast, v^\dagger \,|\, q^\ast, q^\dagger, \hat{C}_{2}\right)
= F^{n} \left(v^\ast \,|\, q^\ast, \hat{C}_{2}\right)$;
%$F^{n} \left(v^\ast \,|\, q^\ast, q^\dagger, \hat{C}_{2}\right) = F^{n} \left(v^\ast \,|\, q^\ast, \hat{C}_{2}\right)$, i.e., the extra information that $\varphi_{2}$ produces the question $q^\dagger$ doesn't affect the expectation of the conditional probability that it produces the answer $v^\ast$ to the question $q^\ast$ (given the information already contained in $\hat{C}_{2}$.);

    \item $\varphi_2$ is embed-calibrated with $\varphi_{1}$ at step $n$ for $\hat{C}_{2}$ for each of the three questions $\psi^{-1}(q^{*})$, $\psi^{-1}(q^{\dagger})$, and $\psi^{-1}(q^{*},q^{\dagger})$, and for $m=1,2$;

    \item ${{D}}[., .]$ is a locally Lipschitz continuous function of its probability distribution arguments (where those distributions are considered as vectors in a Euclidean metric space);

    \item  For any $v\in\mathcal{V}_{2}$,
$$\sum_{v^\dagger}\overline{P}_1 \left(v^\ast, v^\dagger \,|\, q^\ast, q^{\dagger}, E^{-1}[\{(\psi^{-1}(q^\ast,q^\dagger),v)\}\cup \hat{C}_{2}]\right) = \overline{P}_1 \left(v^\ast \,|\, q^\ast, E^{-1}[\{(\psi^{-1}(q^\ast),v)\}\cup \hat{C}_{2}]\right);$$
%i.e., the extra information that the sample of the limit distribution of $\varphi_{1}$ produces the question $q^\dagger$ doesn't affect the conditional probability that it produces the answer $v^\ast$ to the question $q^\ast$ (given the information already contained in $E^{-1}[\{(\psi^{-1}(q^\ast,q^\dagger),v)\}\cup \hat{C}_{2}]$ or $E^{-1}[\{(\psi^{-1}(q^\ast),v)\}\cup \hat{C}_{2}]$ for any $v\in\mathcal{V}_{2}$).
\end{enumerate}
Then
\begin{multline*}
    \sum_{v\in\mathcal{V}_{2}}P^{n}_{2}(v|\psi^{-1}(q^{\dagger}),(q^{*},v^{*}),\hat{C}_{2})\overline{P}_{1}(v^{\dagger}|q^{\dagger},(q^{*},v^{*}), E^{-1}[\{(\psi^{-1}(q^{\dagger}),v)]\}\cup\hat{C}_{2}]) \\ > \ \sum_{v\in\mathcal{V}_{2}}P^{n}_{2}(v|\psi^{-1}(q^{\dagger}),\hat{C}_{2})\overline{P}_{1}(v^{\dagger}|q^{\dagger}, E^{-1}[\{\psi^{-1}(q^{\dagger}),v)\}\cup\hat{C}_{2}]).
\end{multline*}
\end{prop}
\noindent
This proposition shows that if a scientist-SMS believes that one claim $(q^{*},v^{*})$ abductively supports a second claim $(q^{\dagger},v^{\dagger})$, and they are embed-calibrated with the universe-SMS, then in expectation the universe-SMS itself is such that $(q^{*},v^{*})$ abductively supports a second claim $(q^{\dagger},v^{\dagger})$.

% in the scientific case, we note that the switch to the scientific context means that we change condition (2) in \cref{prop:mathabduction} to the requirement that $\varphi_{2}$ (i.e., the scientist-SMS) is embed-calibrated with $\varphi_{1}$ (i.e., the universe-SMS), and change the conclusion to read
% %\begin{linenomath}
% \begin{equation}\label{eq:scienceabductionconc}
%     \overline{P}_{1}(v^{\dagger}|q^{\dagger},(q^{*},v^{*}), E^{-1}[\hat{C}_{2}])> \overline{P}_{1}(v^{\dagger}|q^{\dagger}, E^{-1}[\hat{C}_{2}])
% \end{equation}
% %\end{linenomath}
% for any $v\in\mathcal{V}_2$. This states that the universe-SMS $\varphi_{1}$ follows the inductive inference pattern with respect to the claims $(q^{\dagger},v^{\dagger})$ and $(q^{*},v^{*})$ conditional on any claim set $E^{-1}[(\psi^{-1}(q^{\dagger}),v),\hat{C}]$ representing the physical instantiation (e.g., a neuro-biological instantiation) of the reasoning of a scientist-SMS that is embedded within it. Note that $(\psi^{-1}(q^{\dagger}),v)$ is a claim, made by the scientist-SMS, specifying a probability distribution over possible answers to the question $q^{\dagger}$.\par

Next, as we did in the case of multiple evidence paths, we use the assumption that $\varphi_{2}$ is embed-calibrated with $\varphi_{1}$ to ``project down'' from the fact that the outcomes of experiments in the universe conform with the abductive inference pattern, to establish the same property for the observations by the scientists of the outcomes of those experiments:
\begin{prop}\label{prop:scienceabduction}
%Suppose we are given some $N, n \in \Z^+$, backward-consistent SMS $\varphi_1$, SMS
%$\varphi_2$, and associated partial functions $\psi, \Psi$ as in \cref{def:prediction_pair}. Suppose as well that
%\begin{enumerate}
%    \item  The distribution $\overline{P}_{1}$  satisfies the abductive premise, \cref{eq:5}, for all $\hat{C}_{1}\in\mathcal{\hat{C}}_{1}$ and
%for the specific questions $q^*$, $q^\dagger$ and answers $v^{*},v^\dagger$ ;
%
%    \item $(\psi^{-1}(q^{*}),\hat{C}_{2})$, $(\psi^{-1}(q^{\dagger}),\hat{C}_{2})$, and $(\psi^{-1}(q^{*},q^{\dagger}),\hat{C}_{2})$ are all embedded prediction pairs for $\varphi_{2}$ at step $n$ for any $\hat{C}_{2}\in\hat{\mathcal{C}}_{2}$;
Suppose we are given some $N, n \in \Z^+$, backward-consistent SMS $\varphi_1$, SMS
$\varphi_2$, and associated partial functions $\psi, \Psi$ as in \cref{def:prediction_pair}.
Let  $v^*, v^\dagger \in \mathcal{V}_1$, and let $q^*, q^\dagger \in \mathcal{Q}_1$ be
two questions such that all three tuples $q^*$, $q^\dagger$ and $(q^*, q^\dagger)$ are in the codomain of $\psi$.
Let $E$ be an embedding function. Suppose as well that

\begin{enumerate}
\item  $\PPP$ satisfies the abductive premise, \cref{eq:5}, for all $\hat{C}_{2}\in\mathcal{\hat{C}}_{2}$ and
for the specific questions $q^*$, $q^\dagger$ and answers $v^{*},v^\dagger$;

%    \item The prediction distribution
% $\PPP$ satisfies the abductive premise, \cref{eq:5}, for all $\hat{C}_{2}\in\mathcal{\hat{C}}_{2}$ and
%for the specific questions $q^*$, $q^\dagger$ and answers $v^{*},v^\dagger$;

    \item $(\psi^{-1}(q^{*}),\hat{C}_{2})$, $(\psi^{-1}(q^{\dagger}),\hat{C}_{2})$, and $(\psi^{-1}(q^{*},q^{\dagger}),\hat{C}_{2})$ are all embedded prediction pairs for $\varphi_{2}$ at step $n$ for any $\hat{C}_{2}\in\hat{\mathcal{C}}_{2}$;

    \item For any $v\in\mathcal{V}_{2}$, $$\sum_{v^\dagger}\overline{P}_1 \left(v^\ast, v^\dagger \,|\, q^\ast, q^{\dagger}, E^{-1}[\{(\psi^{-1}(q^\ast,q^\dagger),v)\}\cup \hat{C}_{2}]\right) = \overline{P}_1 \left(v^\ast \,|\, q^\ast, E^{-1}[\{(\psi^{-1}(q^\ast),v)\}\cup \hat{C}_{2}]\right);$$
% i.e., the extra information that the sample of the limit distribution of $\varphi_{1}$ produces the question $q^\dagger$ doesn't affect the conditional probability that it produces the answer $v^\ast$ to the question $q^\ast$ (given the information already contained in $E^{-1}[\{(\psi^{-1}(q^\ast,q^\dagger),v)\}\cup \hat{C}_{2}]$ for any $v\in\mathcal{V}_{2}$ or $E^{-1}[\{(\psi^{-1}(q^\ast),v)\}\cup \hat{C}_{2}]$ for any $v\in\mathcal{V}_{2}$ ).

    \item $\varphi_2$ is embed-calibrated with $\varphi_{1}$ at step $n$ for $E[\hat{C}_{1}]$ for each of the three questions $\psi^{-1}(q^{*})$, $\psi^{-1}(q^{\dagger})$, and $\psi^{-1}(q^{*},q^{\dagger})$, and for $m=1,2$;

    \item ${{D}}[., .]$ is a locally Lipschitz continuous function of its probability distribution arguments (where those distributions are considered as vectors in a Euclidean metric space);

    \item $\sum_{v^\dagger}\PPP \left(v^\ast, v^\dagger \,|\, q^\ast, q^\dagger, E[\hat{C}_{1}]\right) = \PPP \left(v^\ast \,|\, q^\ast, E[\hat{C}_{1}]\right)$;
%, i.e., the extra information that the sample of the limit distribution of $\varphi_{2}$ produces the question $q^\dagger$ doesn't affect the expectation of the conditional probability that $\varphi_{2}$ produces the answer $v^\ast$ to the question $q^\ast$ (given the information already contained in $E[\hat{C}_{1}]$.)
\end{enumerate}
Then $\PPP(v^{\dagger}|q^{\dagger},(q^{*},v^{*}),E[\hat{C}_{1}])>\PPP(v^{\dagger}|q^{\dagger}, E[\hat{C}_{1}])$.
\end{prop}
Intuitively, \cref{prop:scienceabduction} shows that if the occurrence of the experimental outcome $(q^{*},v^{*})$ increases the probability that the universe-SMS assigns the question $q^{\dagger}$ the answer $v^{\dagger}$, then for a scientist who is embed-calibrated with the universe that embeds them, conditioning on $(q^{\dagger},v^{\dagger})$ ought to increase their expected degree of belief that the question $q^{*}$ has the answer $v^{*}$, in keeping with the abductive inference pattern. To the extent that being embed-calibrated with their physical universe is an appropriate normative goal for an empirical scientist, this result shows that such a scientist should engage in abductive reasoning, in expectation.\par

%\dhwc{Again, fix the following to not involve step-sensitivity. (In fact, despite what the
%text before the proposition says, is step-sensitivity actually used?)}

\section{Relation with previous work}
\label{sec:earlier_work}
%Before further developing our framework, we make
There are several points worth elaborating about how the SMS framework differs from frameworks
previously investigated in the literature.
First, it is important to contrast the SMS framework with concepts that have
been studied in \textit{probability logic} (e.g., Carnap \citeyear{carnap1962logical}, Burgess \citeyear{burgess1969probability}, Hoover \citeyear{hoover1978probability}, Leblanc \citeyear{leblanc1979probabilistic}, Hailperin et al. \citeyear{hailperin1984probability}, Nilsson \citeyear{nilsson1986probabilistic}, Fagin, Halpern, and Meggido \citeyear{fagin1990logic}, Leitgeb \citeyear{leitgeb2008probabilistic}, Haenni et al. \citeyear{haenni2010probabilistic}, Christiano et al. \citeyear{christiano2013definability}, and Campbell-Moore \citeyear{campbell2015express}).
When building a system of probability logic, one seeks ``to start with a classical (propositional/modal/etc.) system of logic and to `probabilify' it in one way or another, by adding probabilistic features to it'' (Demey, Kooi, and Sack \citeyear{sep-logic-probability}). That is, whereas ordinary system of logic aim to axiomatize truth-preserving inferences between sentences that do not contain probabilistic language, probability logics seek either to axiomatize the process of making \textit{probability}-preserving inferences from one sentence to another, or to axiomatize the process of making truth-preserving inferences between sentences that contain probabilistic language.
(See also work on \textit{Markov logic networks} due to \cite{richardson2006markov}.)
In contrast, the SMS framework does not require the axiomatization of either truth-preserving or probability-preserving inferences.
Rather, we aim to provide a framework in which \textit{any} pattern of inferences can be understood as the output of a stochastic process.
%That is, we want an approach that does not require one to first specify a particular ``base logic'' that is then used to axiomatize probabilistic relations between proposition in that logic.
%This is because we wish to remain
%We are agnostic as to what the underlying logic of formal reasoning is that is ``probabilified'',
%and leave room for the
%In contrast, we allow the possibility that no such underlying logic even exists.
%It is in part for this reason that we develop here a new framework for probabilistic inference, rather than adopting an existing system of probability logic. Having said this, it is worth noting that, because a one-step SMS can be used to represent any probability distribution over some domain, such an SMS can be designed so that it is equivalent to any given system of probabilistic logic, including a Markov logic network or logical inductor.\par

Although his proposal is less technical than the work in probability logic discussed immediately above, \citet{franklin1987non} argues for the descriptive claim that a non-deductive logic, formulated in the language of probability theory, plays a crucial role in the development of mathematical results. His proposal contains a germ of an idea that we have developed here, namely that mathematical reasoning in human agents can be represented as a stochastic process. However, he holds that non-extreme probabilities assigned to mathematical statements should be treated as solely representing the subjective degrees of belief of individual mathematicians, and that such degrees of belief are only justified prior to a statement's being proven. By contrast, we hold that in principle, a mathematician may always be justified in assigning non-extreme degree of belief to a mathematical claim.\par

%\dk{This next paragraph is new.}
Second, our framework shares some features with ``erotetic logic,'' which aims to provide a formal logic of questions and answers. \citet{Belnap1976-BELTLO-2} provide a textbook introduction to the subject, and subsequent developments include Hintikka's Interrogative Model of Inquiry (see \cite{Hintikka1991-HINWIT-8}) and Wi\'sniewski's Inferential Erotetic Logic (IEL) (see \cite {Wisniewski2013-WINQIA}). Our framework and erotetic logic are both primarily concerned with asking and answering of questions. However, research in erotetic logic follows the general methodology of formal logic, defining a syntax, semantics, and proof theory that allows one to prove that certain questions have certain answers. By contrast, we use a probabilistic formalism to descriptively represent the asking and answering of questions.
We also extend this probabilistic formalism to model conditions under which this process of question answering can be held to certain normative standards. In addition,
insofar as it provides a model for inquisitive reasoning under uncertainty, our framework shares some similarities with recent work by \citet{hoek2022questions}. However, the focus of Hoek's work is on the connections between an inquisitive model of belief and decision theory, rather than mathematical or scientific reasoning.\par

Third, our approach has significant similarities to work from \citet{garrabrant2016logical,garrabrant2017formal}, who present a computable algorithm that assigns probabilities to well-formed formulas in a given logical language, and is able to update those probabilities as it obtains more evidence about possible theorems of that language. We submit that their algorithm is a particular instance of the more general class of SMSs. The type of SMS Garrabrant et al.\ define is one that considers questions of the form `what is the probability distribution over possible answers to whether $\phi$ a theorem of the logical language?' and returns a probability distribution over the possible answers `Yes' and `No.' As a point of contrast between their approach and ours, we note that the normative standards that they impose on their algorithm, where that algorithm is understood as an SMS, are entirely \textit{internal} to the SMS. For their algorithm to reason ``correctly,'' they require only that it is impossible for a Turing machine to construct a series of bets that leads the algorithm to sure loss in time polynomial in the size of the input to that Turing machine. In other words, they require that their algorithm not be internally inconsistent in a way that can be efficiently exploited. By contrast, in our approach an SMS is subject to normative standards that are \textit{external} to the SMS itself, in virtue of its calibration with a second SMS. We note that one can draw a similar comparison between our approach and work by \citet{lample2019deep}, in which they use deep neural networks to approximate integrals and solve differential equations. These deep neural networks can be thought of a special case of an SMS, in which questions are considered and answers are output depending on the pattern of activation within the neural network, with the distribution over possible outputs updated via supervised learning.\par

We note also that the approaches due to Garrabrant et al.\ do not provide a mechanism for considering questions that are about the \textit{probability distribution} over the answers produced by a different
SMS, made by that SMS in response to other questions.
This flexibility is particularly important when using SMSs to model human scientists making predictions
about the outcome of a future experiment concerning the embedding SMS of the physical universe, since such
predictions invariably take the form of probability distributions. For example, when asked to
predict what the weather will be tomorrow over London, commonly a scientist would respond
by providing a probability distribution over the possible states of the weather, rather than a point prediction.\par

Fourth, we note the connections between our work and work in computational cognitive science on methods for defining probabilistic generative models, which are models that enable agents to generate random samples from a distribution. These include the programming language Church, first introduced by \citet{goodman2012church}, and the related probabilistic Turing machine $\texttt{QUERY}$, which is introduced by \citet{freer2014towards}, and uses conditional simulation to implement hierarchical Bayesian inference. As in the case of Garrabrant et al.'s algorithm, we understand $\texttt{QUERY}$ as a particular type of SMS that answers questions about its environment by randomly sampling a probability distribution over possible answers to a question, randomly sampling an answer from the sampled distribution, and comparing its predicted answer to observed answers, updating both of the distributions that it samples from over time. Notably, $\texttt{QUERY}$ is the sort of SMS that can answer questions about the probability distributions over answers to some question, since it uses this information as part of its hierarchical Bayesian architecture. However, the category of SMSs is more general that the category of inference techniques implementable in $\texttt{QUERY}$. Unlike $\texttt{QUERY}$, the distributions the define an SMS are not required to be computable, nor do Freer, Roy, and Tenenbaum provide any explicit notion of semantics or normative constraints for $\texttt{QUERY}$ that could be analogized to the notion of calibration for an SMS. That said, more carefully considering the connections between SMSs, $\texttt{QUERY}$, heirarchical Bayesian inference, and calibration is an intriguing avenue for future work.\par

Fifth, Icard (\citeyear{icard2020calibrating}) considered a probabilistic version of the Chomsky-Sch\"utzenberger (\citeyear{chomsky1959algebraic}) hierarchy of linguistic grammars. Each grammar in Icard's hierarchy is defined by a set of restrictions on the rules for generating finite strings from other finite strings by sampling from other finite strings. Icard derives results showing that as we move to less-restricted probabilistic grammars, we are able to use those grammars to define increasingly general probability distributions. (For example, the most restricted class of probabilistic grammars can define finite-support, rational-valued probability distributions, while the least restricted can define any probability distribution that can be implemented by a probabilistic Turing machine.) Any SMS in which the claim set is restricted to pairs of finite strings can be placed in Icard's hierarchy, with implications for the probability distribution(s) that define each SMS. However, SMSs also allow for even more general probabilistic grammars, if we allow for infinite strings in the claim set. An intriguing avenue for future work would be to consider how calibration between SMSs is effected by where those two SMSs stand in Icard's hierarchy. (See also work by \citet{lin2017critical} on criticality in probability distributions definable at different levels of a probabilistic Chomsky hierarchy.)\par

Sixth and finally, we note recent work in the foundations of physics by Gisin (\citeyear{gisin2019indeterminism}, \citeyear{gisin2021indeterminism}). He argues that even classical mechanics should be understood as indeterministic, on the grounds that real-valued quantities are physically unrealistic, and should be replaced with \textit{intuitionistic numbers}, which are finite descriptions of stochastic processes that produce real numbers in the infinite limit. By representing physical systems (even classical physical systems) as stochastic processes, Gisin's framework bears a significant resemblance to ours. However, there are also crucial differences. In particular, Gisin argues that indeterminism in physical systems is driven by the fact that physical quantities only exist in finite lime, and so their precise value is not allowed to tend towards a limit. This emphasis on the role of time in explaining key aspects of physical reality leads him to explicitly reject the block universe hypothesis in favor of an ontologically meaningful notion of time (\citeyear{gisin2021indeterminism}, p.\ 13364). By contrast, within the SMS framework time plays no explicit role in conceptualizing indeterminism in physical systems. As explicitly noted above, we intend for our framework to be consistent with a block-universe view.\par

\section{Discussion}\label{sec:conclusion}
%This completes our initial exegesis of the stochastic mathematical systems framework, and our application of the framework to two aspects of mathematical and scientific reasoning. To conclude, we begin by commenting on the fact that, in the examples above,

There are several aspects of the stochastic mathematical systems framework we introduced and
explored above that are worth exploring further. First, note that we needed to introduce additional
structure when formalizing scientific reasoning, beyond the formal structure needed to formalize
mathematical reasoning. In particular,
we had to introduce additional formal structure (namely, the concept of an embedding) to capture our normative constraints on scientific reasoning, beyond the formal structure needed to capture normative constraints on
mathematical reasoning. This reflected the fact that we take abstract mathematical reasoning by a
human reasoner not to be normatively constrained by the physical universe in which the mathematician finds themselves, whereas we do take scientific reasoning to be normatively constrained in this way. This amounts to a sharp distinction between the norms of mathematics and the norms of science that we hope to explore in future work.

Focusing just on the case of scientific reasoning, there has been previous work that considers the implications of a reasoner being embedded in the system they are reasoning about. In particular, see the ``inference devices'' framework of \citet{wolpert2017constraints}. Within that framework a ``monotheism'' theorem is proven, showing that no two reasoners embedded in the same physical system can both be perfectly accurate predictors. An intriguing avenue for future work
on the SMS framework is to see whether similar results hold for two scientist-SMSs both embedded within the same universe-SMS.\par

%Sticking with the question of the limits of reasoning in the SMS framework, a
Another potentially fruitful line of research might be consider whether there are analogs of G\"odel's second incompleteness theorems for the SMS framework. Specifically, potential future work would investigate whether, if an SMS assigns non-zero probability to both answers to all true/false questions, then it would in particular assign non-zero probability to the answer of `false' to self-referential questions about that same SMS's probabilistic pattern of answers to questions.\par

In a similar vein, we believe that our approach may be fruitful for extending the ``no free lunch theorems'' of \citet{wolpert1996lack,woma97,wolpert2021implications}, which provide a set of formal bounds on how well any machine learning algorithm or search algorithm can be guaranteed to perform if one does not make assumptions for the prior probability distribution of the underlying stochastic process, to mathematical and scientific learning and reasoning more generally. Fourth, we hope to study the relationships between multiple mathematician-SMSs or scientist-SMSs that are calibrated with the same community-SMS and/or embedded in the same universe-SMS. Our hope is that such an inquiry will be fruitful for understanding, in a general way, the properties of reasoners that, despite differing in crucial respects, share the same normative goals.\par

%There are still more areas of research wherein we believe that our framework might be directly applicable.
Another potential line of future research involves  the ``problem of (logical) omniscience'' in epistemic logic.
Briefly, the problem is that many forms of epistemic logic suppose that if a reasoner knows any proposition
$A$, and knows that $A$ implies $B$ according to some canonical set of axioms, then they know $B$. For example, according to many epistemic logics
a reasoner who knows the axioms of number theory automatically must also
know all theorems of number theory, which is clearly absurd. However, consider
defining what an SMS at a step $n$ ``knows'' to be the contents of a claim set it produces at that step.
In general, an SMS can with nonzero
probability produce a claim $A$, a claim  $A \rightarrow B$, and also a claim $\neg B$.
Thus, the problem of logical omniscience disappears in the SMS framework, if we adopt its natural
definition of what it means to ``know''.
These considerations can also be generalized into a Bayesian setting. There the problem of omniscience
%Traditional Bayesian approaches to reasoning suffer from the ``problem of omniscience.''
is that, because the sample space over which probabilities are defined is typically assumed to be closed under the entailment relation of some logic, a reasoner whose degrees of belief are modeled by that probability function must be represented as maximally confident in any theorems of that logic. See \citet{skipper2020bayesianism} for a recent and more detailed discussion of this problem in a Bayesian context. The SMS framework makes no assumptions about the closure of the sample space under any logic, and so might therefore be seen as immune to the problem of logical omniscience, though more work is needed to defend this thesis rigorously.\par

As another possible research direction, we note that in both the mathematical universe approach of Tegmark (\citeyear{tegmark1998theory}, \citeyear{tegmark2008mathematical}, \citeyear{tegmark2014our}), and the work of the ontic structural realists, physical reality is given a \textit{structural} representation as a collection of functions between sets. While this may seem to be a very different approach to our representation of the physical universe as an SMS, we hope in future work to show how one can model structure as being generated through the SMS's process of asking and answering questions. We suspect that in so doing, we can shed new light on the relationships between computability, indeterminism, and the structuralist approach to the representation of physical reality.\par

More speculatively, we are also interested in whether the framework presented here may be helpful in quantifying the extent to which foundational results in mathematics, like G\"odel's incompleteness theorem or the undecidability of the halting problem are robust to small perturbations in the oracle-SMS that is used to set the standards of mathematical correctness. This could in turn lead to improved understanding of the modal status of these results across alternative possibilities with respect to the norms of mathematical reasoning.

Having said all this, we hope that there are many more applications, both practical and theoretical, of the approach presented here that we have not yet anticipated. Indeed, the two applications of our SMS framework
presented in this paper--- establishing the
benefit of abduction and establishing the benefit of multiple lines of evidence --- do not
exploit the fact that the relevant response distributions arise via an infinitely iterated stochastic process. The
analysis would still go through if that were not the case. However, we expect that some of these future
directions of research on other applications of the SMS framework would exploit that fact.

%
%\printbibliography

\begin{appendix}

\section{Additional Facets of the SMS Framework}\label{sec:appendixfullframework}
\subsection{Claim Trajectories}
Although it was not needed for the results or discussion in this paper, we can also use the SMS framework to assign probabilities to trajectories of claims produced at different steps of an SMS.
\begin{Def}
A \textbf{claim vector trajectory} $\vec{C}^{n}$ is a map taking any $i \in 1, \ldots, n$ to a claim vector if $n$ is finite, and taking any $i \in \mathbbm{Z}^+$ to a claim vector if $n = \infty$.
\end{Def}
\noindent
As shorthand, if a claim vector
trajectory contains only a single claim vector, then we will sometimes write that trajectory as that single claim vector.
In the sequel the spaces $\mathcal{Q}$ and $\mathcal{V}$ will often be implicit. We will also cavalierly use
$|.|$ to indicate cardinality. So for example, $|C|$ is the number of components of a claim vector
$C$, $|\hat{C}|$ is the number of elements in the claim set $\hat{C}$, etc.\par

It will occasionally be useful to consider probability distributions concerning the event that the claim vector
generated at the $n$-th step of a particular SMS contains the answer $v$ in response to the question $q$.
We present some definitions that will be useful when
working with such distributions for the case
where $\mathcal{C}$ is countable; their extensions to other spaces is straight-forward.
\begin{Def}
\label{def:2.10}
Suppose we are given an SMS $\varphi=(\mathfrak{C},X)$ and a claim vector trajectory $\vec{C}^{n-1} =
(C^1, C^2, \ldots, C^{n-1})$ for some finite $n > 1$.
The associated \textbf{question semi-distribution} is the function mapping all $q\in{\mathcal{Q}}$ to the value
%\begin{linenomath}
\begin{equation*}
    P^n(q, \vec{C}^{n-1}) := \sum_{v^{\prime}\in{\mathcal{V}}}
        P_{\mathfrak{C}} \left(X(1)=C^{1},\dots, X(n-1) = C^{n-1}, (q, v') \in \hat{X}(n)\right).
\end{equation*}%\end{linenomath}
\end{Def}
\noindent
$P^n(q, \vec{C}^{n-1})$ is the joint
probability that the question $q$ occurs in a claim in the step-$n$ claim vector and that the earlier claim vectors were $C^1, C^2, \ldots, C^{n-1}$. Intuitively, one can think of this as the probability that a reasoner considers a particular question at step $n$ and also outputs a particular set of claims at previous steps. Note that in general, for any fixed $\vec{C}^{n-1}$,
the associated question semi-distribution is not normalized when considered as a function
of questions $q$.\par

We define claim semi-distributions $P^n((q,v),\vec{C}^{n-1})$
analogously to question semi-distributions. We combine these definitions as follows:
\begin{Def}
For an SMS $\varphi=(\mathfrak{C},X)$, question $q$, and claim vector trajectory $\vec{C}^{n-1}$
where $P^n(q, \vec{C}^{n-1}) \ne 0$, the associated
\textbf{response distribution} is the map from all $v\in{\mathcal{V}}$ to the value
%\begin{linenomath}
\begin{equation*}
    P^n(v | q,\vec{C}^{n-1})=\frac{P^n((q, v), \vec{C}^{n-1})}{P^n(q, \vec{C}^{n-1})}.
\end{equation*}%\end{linenomath}
\end{Def}
\noindent
The response distribution for a given SMS $\varphi$, question $q$, and claim vector trajectory $\vec{C}^{n-1}$ specifies the probability that $q$ is given an answer $v$ in the $n$-th step of $\varphi$, conditional on $\varphi$ having produced the earlier claims specified in the partial claim vector trajectory $\vec{C}^{n-1}$ and having produced \textit{some} claim that includes the question $q$ in the $n$-th claim vector.\footnote{The terms `earlier' and `later' implicitly presuppose a temporal interpretation of the ordering of the SMS. While such an interpretation is not required, it is a natural one, especially when our framework is taken to represent mathematical reasoning by actual humans.} This conditional character of the response probability distribution will be exploited
below to represent path-dependent mathematical reasoning, in which the answers given by an SMS
to questions earlier in the process of generating mathematical claims
affect the answers given to questions posed later in the process.\par

\subsection{An Alternative Distribution Over Claim Sets}
\label{app:alt_distribution}
In the body of this paper, we considered the following distribution over claim sets:
%\begin{linenomath}
\begin{equation}
    P^{n}(\hat{C}):=\sum_{\{C\in\mathcal{C}^{*}:\hat{C}\subseteq U(C)\}}P_{\mathfrak{C}}(X(n)=C).
\end{equation}
%\end{linenomath}
This is the probability that an SMS outputs, at step $n$, a claim vector whose unordering is a superset of $\hat{C}$. But we may also want to consider the probability that an SMS outputs, at a step $n$, a vector whose unordering simply \textit{is} the claim set $\hat{C}$. Such a distribution, which we denote $\underline{P}^{n}(\hat{C})$, can be defined in the obvious way:
%\begin{linenomath}
\begin{equation}
    \underline{P}^{n}(\hat{C}):=\sum_{\{C\in\mathcal{C}^{*}:\hat{C}= U(C)\}}P_{\mathfrak{C}}(X(n)=C).
\end{equation}
%\end{linenomath}
We can then use this distribution to define a related question semi-distribution and response distribution:
%\begin{linenomath}
\begin{equation}
        \underline{P}^{n}((q,v),\hat{C}):=\sum_{\{C\in\mathcal{C}^{*}:\hat{C}\cup\{(q,v)\}= U(C)\}}P_{\mathfrak{C}}(X(n)=C)
\end{equation}
%\end{linenomath}
%\begin{linenomath}
\begin{equation}
        \underline{P}^{n}(q,\hat{C}):=\sum_{v\in\mathcal{V}}\underline{P}^{n}((q,v),\hat{C})
\end{equation}
%\end{linenomath}
%\begin{linenomath}
\begin{equation}
        \underline{P}^{n}(v|q,\hat{C}):=\frac{\underline{P}^{n}((q,v),\hat{C})}{\underline{P}^{n}(q,\hat{C})}
\end{equation}
%\end{linenomath}
Finally, we can extend this distribution to apply to trajectories of claim \textit{sets}, rather than claim vectors, e.g.:
%\begin{linenomath}
    \eq{
\underline{P}^n(q,\vec{\hat{C}}^{n-1}) &:= \sum_{\vec{C}^{n-1} : U(C^i) = \hat{C}^i \, \forall 1 \le i \le n-1} P^n(q,\vec{C}^{n-1}),
}
%\end{linenomath}
where $\vec{\hat{C}}^{n-1}$ is a trajectory of unordered claim sets $(\hat{C}^{1},\dots,\hat{C}^{n-1})$.

\section{Proofs of Propositions}
\label{sec:proofs}
\subsection{Proof of Lemma \ref{lem:limit}}
\begin{proof}
Choose some SMS that is backward-consistent after step $\kappa$, and some claim set
$\hat{C}$. For all steps $j > \kappa$ and $j > i > \kappa$,
\eq{
P^j(\hat{C}) &=  \sum_{\hat{C}' : \hat{C} \subseteq \hat{C}'} \PP^j(\hat{C}')   \nonumber \\
	&=  \sum_{{C}' : \hat{C} \subseteq U({C}')} P_{\mathfrak{C}} (X(j, .) = {C}')  \nonumber \\
	&=  \sum_{{C}' : \hat{C} \subseteq U(C')} \left[  \sum_{{C}'' : \hat{C} \subseteq U(C'')} P_{\mathfrak{C}} (X(j, .) = {C}', X(i, .) = {C}'')
				+ \sum_{{C}'' : \hat{C} \not \subseteq U(C'')} P_{\mathfrak{C}} (X(j, .) = {C}', X(i, .) = {C}'')  \right]  \nonumber \\
	&=  \sum_{{C}' : \hat{C} \subseteq U(C')} \bigg[  \sum_{{C}'' : \hat{C} \subseteq U(C'')}
			P_{\mathfrak{C}} (X(j, .) = {C}' \, \vert \, X(i, .) = {C}'') P_{\mathfrak{C}} (X(i, .) = {C}'') 	\nonumber \\
				&\qquad\qquad \qquad\qquad\qquad\qquad \qquad\qquad
									+ \sum_{{C}'' : \hat{C} \not \subseteq U(C'')} P_{\mathfrak{C}} (X(j, .) = {C}', X(i, .) = {C}'')  \bigg]
}
Due to backward-consistency, $U(C^{\prime\prime})\subseteq U(C^{\prime})$ and so if $\hat{C}\subseteq U(C^{\prime\prime})$, then $\hat{C}\subseteq U(C^{\prime})$. This means that for all $C^{\prime\prime}$ such that $\hat{C}\subseteq U(C^{\prime\prime})$,
\begin{equation}
    P_{\mathfrak{C}}(X(j,.)\in\{C^{\prime}: \hat{C}\subseteq U(C^{\prime})\}|X(i,.)=U(C^{\prime\prime}))=1.
\end{equation}
Thus, the sum over $C'$ of the conditional distribution in the summand in the last line equals $1$,
for all $C''$ in the associated sum. Therefore we get
\eq{
P^j(\hat{C}) &=    \sum_{{C}'' : \hat{C} \subseteq U(C'')}
		P_{\mathfrak{C}} (X(i, .) = {C}'')
				+  \sum_{{C}' : \hat{C} \subseteq U(C')}
						\sum_{{C}'' : \hat{C} \not \subseteq U(C'')} P_{\mathfrak{C}} (X(j, .) = {C}', X(i, .) = {C}'')    \nonumber \\
	&= P^i(\hat{C}) +  \sum_{{C}' : \hat{C} \subseteq U(C')}
				 \sum_{{C}'' : \hat{C} \not \subseteq U(C'')} P_{\mathfrak{C}} (X(j, .) = {C}', X(i, .) = {C}'')    \nonumber \\
	&\ge P^i(\hat{C})
}
So for any claim set $\hat{C}$, $P^j(\hat{C})$ is a monotonically increasing function of $j$. In addition,
that probability is upper-bounded as one varies over all $j$, by $1$. By the completeness axiom of the
reals, that means it has a least upper bound. Therefore the limit of that probability goes to infinity is
well-defined.  Replacing $j$ with $m$ establishes the claim.
\end{proof}

\subsection{Proof of \cref{prop:projection}}
%\begin{linenomath}
\begin{proof}
The divergence function $D$ in the definition of embed-calibrated can only equal $0$ if
its two arguments are identical. Therefore for all embedded prediction triples $(v,  q, \hat{C})$, and $v^m_1$,
\eq{
\Psi(\psi(q), v)(v_1^m) =
	\overline{P}_1\left( v_1^m \;|\; \psi(q), E^{-1}[\{(q, v)\} \cup {\hat{C}}]  \right).
}
Define
\eq{
\hat{C}_2(v^m_1) &:= E \left[\{(\psi(q), v^m_1)\} \cup E^{-1}[\{(q, v)\} \cup {\hat{C}}]\right].
}
Since $(q, \hat{C})$ is a discriminating pair,
\eq{
E^{-1}(\hat{{C}}_2(v^m_1)) = \{(\psi(q), v^m_1)\} \cup E^{-1}[\{(q, v)\} \cup {\hat{C}}].
}
By the definition of embedding function, this means that
\eq{
\overline{P}_{1} \left(\{(\psi(q), v^m_1)\} \cup E^{-1}[\{(q, v)\} \cup {\hat{C}}]\right)
	&= P^{n}_{2} \left(\hat{C}_2(v^m_1)\right).
}
Since this is true for all $v^m_1$, it also is true when we sum over such $v^m_1$. Therefore
\eq{
\overline{P}_{1} \left( v^m_1 \,|\, \psi(q), E^{-1}[\{(q, v)\} \cup \hat{C}] \right)
		&= \frac{P^{n}_{2} (\hat{C}_2(v^m_1))} {\sum_{v^m_1} P^{n}_{2} (\hat{C}_2(v^m_1))}.
}
(See discussion just below \cref{def:responsedist}.) This completes the proof.
\end{proof}
%\end{linenomath}

\subsection{Derivation of \cref{def:3.4}}
We begin with an alternative definition, that a collection of evidence paths
 $\mathcal{B}$ is a set of claim sets $\{B(i)\}$ such that
the premise:
%The first of the requirements is that
%\eq{
%\label{eq:3.22}
%&	\forall \, i \in \{1, \ldots, n\} \qquad P(v^\ast \,|\, q, \beta, B(i)) > P(v^\ast \,|\, q, \beta)
%}
%Intuitively speaking, \cref{eq:3.22} tells us that each $B(i)$ is a line of evidence for the claim $(q, v^\ast)$, when considered in isolation of all the others.
%
%The second requirement is more subtle, providing one formalization of what it means for the
%different evidence paths not to ``thwart one another'':
%\begin{linenomath}
\eq{
\label{eq:3.4}
& \forall \, i \in \{2, \ldots, n\} \qquad
 \dfrac{P\left(B(1), \ldots, B(i) \,|\,  (q, v^\ast), \beta\right)} {P\left(B(1), \ldots, B(i-1) \,|\, \beta, (q, v^*)\right) \times P\left(B(i) \,|\, \beta, (q, v^*)\right)}  \\
			&\qquad\qquad\qquad \qquad\qquad \qquad \ge
				\dfrac{P\left(B(1), \ldots, B(i) \,|\, q, \beta\right)} {P\left(B(1), \ldots, B(i-1) \,|\, q, \beta) \times P(B(i) \,|\, q, \beta \right)}  \nonumber
}
%\end{linenomath}
%\begin{Def}
%If  \cref{eq:3.2,eq:3.4} holds for some collection $\mathcal{B}$, claim $(q, v^*)$, and claim set $\beta$, where all of the conditional probabilities in those equations are nonzero, then we refer to each of the sets $B(i)$ in that collection as an \textbf{evidence path} (for the claim $(q, v^\ast)$ under the distribution $P$ and conditioned on the claim set $\beta$).
%\label{def:3.4}
%\end{Def}
\noindent
Intuitively speaking, \cref{eq:3.2} tells us that each $B(i)$ is a line of evidence for the claim $(q, v^\ast)$, when considered in isolation of all the others.  \cref{eq:3.4} goes on to tell us that those different lines of evidence do not ``work at cross-purposes'', thwarting one another. (The precise indexing of the $n$ sets in $\mathcal{B}$ is not important, so long
as we can find such an indexing for which \cref{eq:3.4} holds.)

%
% for all $i : N \ge i \ge 1$,
%\eq{
%P(v^\ast \,|\, q, \beta, B(1), \ldots, B(i)) &> P(v^\ast \,|\, q, \beta, B(1), \ldots, B(i-1))

To illustrate \cref{eq:3.4}, suppose that $B(1)$ and $B(2)$ are two lines of reasoning that each establish
for a given mathematician-SMS that (conditional on $\beta$) the answer to $q$ is $v^{*}$.
However, suppose that both lines of reasoning depend on steps such that, if the answer to $q$
did indeed turn out to be $v^{*}$ rather than something else, it would imply to
the mathematician-SMS that $B(1)$ and $B(2)$ are very unlikely to both be true.
%such that they will be output in the long run by the commmunity-SMS.
Thus, from the perspective of the mathematician-SMS, $P\left(B(1), B(2) \,|\,  (q, v^\ast), \beta\right) \ll P\left(B(1), B(2) \,|\, q, \beta\right)$.
This is one way of formalizing the notion that $B(1)$ and $B(2)$ thwart
one another, {conditioned on the answer $v^\ast$
in response to the question $q$ (in addition to the foundational claim set $\beta$)}, but not if that answer is
unspecified.\par

We still need to quantify the strength of that thwarting though. To see how to do this, suppose that
although the supposition that $q$ has the answer $v^{*}$ might increase the mathematician-SMS's degree of belief that $B(1)$ is true, the overall effect is mild:\ $P(B(1)|(q,v^{*}), \beta)$ is not much greater than $P(B(1)|q,\beta)$. Suppose that similarly,
$P(B(2)|(q,v^{*}), \beta)$ is not
much greater than $P(B(2)|q,\beta)$. Under these conditions, \cref{eq:3.4} may not hold for $B(1)$ and $B(2)$, so that they are not both evidence paths for the claim $(q,v^{*})$.
Intuitively, they thwart each other too much. Thus, enforcing \cref{eq:3.4} ensures that such thwarting does not occur.\par

If $\mathcal{B}$ is a collection of evidence paths, then not only is each $B(i)$ a line of evidence for the answer $v^\ast$ to the question $q$ when considered in isolation, but furthermore, as we iteratively combine more and more of those evidence paths, we strictly increase the probability of the response $v^\ast$ to the question $q$:
\begin{prop}
\label{prop:used_to_be_}
Let $\mathcal{B}$ be a collection of $N$ evidence paths
for the claim $(q, v^*)$. Then for all $i : N \ge i \ge 1$,
\eq{
P(v^\ast \,|\, q, \beta, B(1), \ldots, B(i)) &> P(v^\ast \,|\, q, \beta, B(1), \ldots, B(i-1))
}
(where we interpret the $i=1$ version of this statement to mean $P(v^\ast \,|\, q, \beta, B(1)) \ge P(v^\ast \,|\, q, \beta)$).
\end{prop}

\begin{proof}
Evaluating \cref{eq:3.2} for $i = 1$ directly establishes the $i=1$ version of \cref{eq:used_to_be_}:
%\begin{linenomath}
\eq{
P(v^\ast \,|\, q, \beta, B(1)) \;>\; P(v^\ast \,|\, q, \beta)  \nonumber
}
%\end{linenomath}
Next, for any $i : 2 \le i \le n$, we can use Bayes' theorem to expand
%\begin{linenomath}
\eq{
\label{eq:3.6}
\dfrac{P(v^\ast \,|\, q, \beta, B(1), \ldots, B(i))}{P(v^\ast \,|\, q, \beta, B(1), \ldots, B(i-1))} &=
			\dfrac{P(v^\ast \,|\, q, {\beta, B(i))}} {P(v^\ast \,|\, q, {\beta})}    \\
			& \qquad \times \dfrac{F^n({B(i)} \,|\, q, {\beta}) \times P({B(1), \ldots, B(i-1)} \,|\, q, {\beta})} {P({B(1), \ldots, B(i)} \,|\, q, {\beta})} \nonumber \\
			& \qquad \times \dfrac{P({B(1), \ldots, B(i)} \,|\,  (q, v^\ast), {\beta})} {P({B(i)} \,|\,  (q, v^\ast), {\beta})  \times
									P({B(1), \ldots, B(i-1)} \,|\,  (q, v^\ast), {\beta})}  \nonumber
}
%\end{linenomath}
Using \cref{eq:3.2} again establishes that the first term on the RHS of \cref{eq:3.6} is $> 1$.
\cref{eq:3.4} establishes that the product of the second and third terms on the RHS is also $> 1$.
Therefore the entire RHS of \cref{eq:3.6} is $> 1$. This establishes the claim.
\end{proof}

So this alternative definition of an evidence path implies the one we adopted in the main text.

\subsection{Proof of Prop.~\ref{prop:3}}
\begin{proof}

By hypothesis,
%$n \in \Z^+$ and $|\psi(q)| = 1$.
%Also by hypothesis, $(q, \hat{C}_2)$ is a prediction pair for step $n$ for any  claim set $\hat{C}_{2}\subseteq \beta \cup \bigcup_i B(i)$.
%By hypothesis it is also true that
$\varphi_2$ is calibrated with the community-SMS $\varphi_1$ at step $n$ for any $\hat{C}_{2}\subseteq \beta \cup \bigcup_i B(i)$.
So for all $1 \le i \le n$,
%\begin{linenomath}
\eq{
		\sum_{v\in\mathcal{V}_{2}}P^{n}_{2}(v|q,\beta, B(1), \ldots, B(i))D\left[\Psi(\psi(q),v)(\mathcal{V}_{1}), \;
			\overline{P}_1\left( \mathcal{V}_1 \;|\; \psi(q),  \beta, B(1), \ldots, B(i) \right) \right]
                &\le \epsilon.
\label{eq:3.85}
}
%\end{linenomath}
The convexity of $D[.,.]$ over its first argument then establishes
%\begin{linenomath}
\eq{
	D\left[\sum_{v\in\mathcal{V}_{2}}P^{n}_{2}(v|q,\beta, B(1), \ldots, B(i))\Psi(\psi(q),v)(\mathcal{V}_{1}), \;
			\overline{P}_1\left( \mathcal{V}_1 \;|\; \psi(q),  \beta, B(1), \ldots, B(i) \right) \right]
                &\le \epsilon.
\label{eq:3.851}
}
%\end{linenomath}
i.e.,
%\begin{linenomath}
\eq{
		D\left[ F^n(\mathcal{V}_1 \,|\, \psi(q), {\beta, B(1), \ldots, B(i)}), \;
			\overline{P}_1\left( \mathcal{V}_1 \;|\; \psi(q), {\beta, B(1), \ldots, B(i)} \right) \right]
                &\le \epsilon.
\label{eq:3.9}
}
%\end{linenomath}
Next, since all of the $B(i)$ are evidence paths for $(\psi(q), v^*)$ under the distribution $F^n$ and
conditioned on $\beta$, by \cref{eq:used_to_be_}, for all $1 \le i \le N$,
%\begin{linenomath}
\eq{
F^n(v^\ast \,|\, \psi(q), {\beta, B(1), \ldots, B(i)}) &> F^n(v^\ast \,|\, \psi(q), {\beta, B(1), \ldots, B(i-1)})
\label{eq:3.10}
}
%\end{linenomath}
Finally, by hypothesis
${{D}}[., .]$ is a locally Lipschitz continuous function of its probability distribution arguments
(where those distributions are considered as vectors in a Euclidean metric space) when evaluated
for the distributions specified in ~\cref{eq:3.9,eq:3.10}.
Then since a divergence equals zero only if its arguments are identical, for all $1 \le i \le n$,
for small enough $\epsilon$,
%\begin{linenomath}
\eq{
\overline{P}_1(v^\ast \,|\, \psi(q), {\beta, B(1), \ldots, B(i)}) &> \overline{P}_1(v^\ast \,|\, \psi(q), {\beta, B(1), \ldots, B(i-1)}).
}
%\end{linenomath}
as claimed.
\end{proof}

\subsection{Proof of Prop.~\ref{prop:3s}}
\begin{proof}
%By hypothesis, $n \in \Z^+$ and $|\psi(q)| = 1$.
%Also by hypothesis, $(q, \hat{C}_2)$ is a prediction pair for step $n$ for any claim set $\hat{C}_{2}\subseteq \beta \cup \bigcup_i B(i)$.
By hypothesis
% it is also true that
$\varphi_2$ is embed-calibrated with the universe-SMS $\varphi_1$ at step $n$ for any $\hat{C}_{2}\subseteq \beta \cup \bigcup_i B(i)$.
So for all $1 \le i \le n$,
%\begin{linenomath}
\begin{multline}\label{eq:3.851a}
\sum_{v\in\mathcal{V}_{2}}P^{n}_{2}(v|q,\beta, B(1), \ldots, B(i))D\Big[\Psi(\psi(q),v)(\mathcal{V}_{1}), \;
			\\ \overline{P}_1\left( \mathcal{V}_1 \;|\; \psi(q),  E^{-1}[\{(q,v)\}\cup\beta, B(1), \ldots, B(i)] \right) \Big]
                \le \epsilon.
\end{multline}
%\end{linenomath}
The convexity of $D[.,.]$ then establishes
%\begin{linenomath}
\begin{multline}\label{eq:3.852}
D\Big[\sum_{v\in\mathcal{V}_{2}}P^{n}_{2}(v|q,\beta, B(1), \ldots, B(i))\Psi(\psi(q),v)(\mathcal{V}_{1}), \;
			\\ \sum_{v\in\mathcal{V}_{2}}P^{n}_{2}(v|q,\beta, B(1), \ldots, B(i))\overline{P}_1\left( \mathcal{V}_1 \;|\; \psi(q),  E^{-1}[\{(q,v)\}\cup\beta, B(1), \ldots, B(i)] \right) \Big]
                \le \epsilon.
\end{multline}
%\end{linenomath}
i.e.,
%\begin{linenomath}
\begin{multline}\label{eq:3.853}
D\Big[F^{n}_{2}(v|q,\beta, B(1), \ldots, B(i)), \;
			\\ \sum_{v\in\mathcal{V}_{2}}P^{n}_{2}(v|q,\beta, B(1), \ldots, B(i))\overline{P}_1\left( \mathcal{V}_1 \;|\; \psi(q),  E^{-1}[\{(q,v)\}\cup\beta, B(1), \ldots, B(i)] \right) \Big]
                \le \epsilon.
\end{multline}
%\end{linenomath}
Next, since all of the $B(i)$ are evidence paths for $(\psi(q), v^*)$ under the distribution $F^n$ and
conditioned on $\beta$, by \cref{eq:used_to_be_}, for all $1 \le i \le N$,
%\begin{linenomath}
\eq{
F^n(v^\ast \,|\, \psi(q), {\beta, B(1), \ldots, B(i)}) &> F^n(v^\ast \,|\, \psi(q), {\beta, B(1), \ldots, B(i-1)})
\label{eq:3.101}
}
%\end{linenomath}
Finally, by hypothesis
${{D}}[., .]$ is a locally Lipschitz continuous function of its probability distribution arguments
(where those distributions are considered as vectors in a Euclidean metric space) when evaluated
for the distributions specified in ~\cref{eq:3.853,eq:3.101}.
Then since a divergence equals zero only if its arguments are identical, for all $1 \le i \le n$,
for small enough $\epsilon$,
%\begin{linenomath}
\begin{multline}
    \sum_{v\in\mathcal{V}_{2}}P^{n}_{2}(v|q,\beta, B(1), \ldots, B(i))\overline{P}_1(v^\ast \,|\, \psi(q), E^{-1}[\{(q,v)\}\cup\beta, B(1), \ldots, B(i)) \\ > \ \sum_{v\in\mathcal{V}_{2}}P^{n}_{2}(v|q,\beta, B(1), \ldots, B(i-1))\overline{P}_1(v^\ast \,|\, \psi(q), E^{-1}[\{(q,v)\}\cup\beta, B(1), \ldots, B(i-1)]).
\end{multline}
%\end{linenomath}
as claimed.
\end{proof}

\subsection{Proof of Prop.~\ref{prop:3b}}
\begin{proof}
By hypothesis,
%$n \in \Z^+$ and $|\psi(q)| = 1$.
%Also by hypothesis,
$(q, E[\hat{C}_1])$ is an embedding prediction pair for step $n$ for any  claim set $\hat{C}_{1}\subseteq \beta \cup \bigcup_i B(i)$.
By hypothesis it is also true that
$\varphi_2$ is embed-calibrated with SMS $\varphi_1$ at step $n$ for the prediction pair $(q, E[\hat{C}_1])$ for any such $\hat{C}_1$.
So
%Since  $|\psi(q)| = 1$, this means that
for all $1 \le i \le n$,
\begin{multline}\label{eq:helper0}
    \sum_{v\in\mathcal{V}_{2}}P^{n}_{2}(v|q,E[\beta, B(1), \ldots, B(i)])D\Big[\Psi(\psi(q),v)(\mathcal{V}_{1}), \\
			\overline{P}_1\left( \mathcal{V}_1 \;|\; \psi(q),  E^{-1}[\{(q,v)\}\cup E[\beta, B(1), \ldots, B(i)]] \right) \Big]
                \le \epsilon.
\end{multline}
Due to the convexity of $D[.,.]$, this means that
%\begin{linenomath}
\begin{multline}\label{eq:helper0.5}
    D\Big[\sum_{v\in\mathcal{V}_{2}}P^{n}_{2}(v|q,E[\beta, B(1), \ldots, B(i)])\Psi(\psi(q),v)(\mathcal{V}_{1}), \\
			\sum_{v\in\mathcal{V}_{2}}P^{n}_{2}(v|q,E[\beta, B(1), \ldots, B(i)])\overline{P}_1\left( \mathcal{V}_1 \;|\; \psi(q),  E^{-1}[\{(q,v)\}\cup E[\beta, B(1), \ldots, B(i)]] \right) \Big]\le \epsilon.
\end{multline}
%\end{linenomath}
or
%\begin{linenomath}
\begin{multline}\label{eq:helper1}
    D\Big[F^n(\mathcal{V}_1 \,|\, \psi(q), E[\beta, B(1), \ldots, B(i)]), \\
			\sum_{v\in\mathcal{V}_{2}}P^{n}_{2}(v|q,E[\beta, B(1), \ldots, B(i)])\overline{P}_1\left( \mathcal{V}_1 \;|\; \psi(q),  E^{-1}[\{(q,v)\}\cup E[\beta, B(1), \ldots, B(i)]] \right) \Big]\le \epsilon.
\end{multline}
%\end{linenomath}
The same line of reasoning for $E[\beta, B(1), \ldots, B(i-1)]$ yields:
%\begin{linenomath}
\begin{multline}\label{eq:helper2}
    D\Big[F^n(\mathcal{V}_1 \,|\, \psi(q), E[\beta, B(1), \ldots, B(i-1)]), \\
			\sum_{v\in\mathcal{V}_{2}}P^{n}_{2}(v|q,E[\beta, B(1), \ldots, B(i-1)])\overline{P}_1\left( \mathcal{V}_1 \;|\; \psi(q),  E^{-1}[\{(q,v)\}\cup E[\beta, B(1), \ldots, B(i-1)]] \right) \Big]
                \le \epsilon
\end{multline}
%\end{linenomath}
Since, by hypothesis, all of the $B(i)$ are evidence paths for $(\psi(q), v^*)$ under the distribution $\overline{P}_{1}$ and
conditioned on $\beta$, for all $1 \le i \le N$,~\cref{eq:used_to_be_} allows us to write
%\begin{linenomath}
\eq{
\overline{P}_1(v^\ast \,|\, \psi(q), \beta, B(1), \ldots, B(i))
	\;>\; \overline{P}_1(v^\ast \,|\, \psi(q), \beta, B(1), \ldots, B(i-1)).
\label{eq:helper3}
}
%\end{linenomath}
Via the fifth condition of
%the antecedent of
the proposition, this implies that
%\begin{linenomath}
\begin{multline}\label{eq:helper3.5}
    \overline{P}_1\left( v^{*} \;|\; \psi(q),  E^{-1}[\{(q,v)\}\cup E[\beta, B(1), \ldots, B(i)]] \right)
	\\ > \overline{P}_1\left( v^{*} \;|\; \psi(q),  E^{-1}[\{(q,v)\}\cup E[\beta, B(1), \ldots, B(i-1)]] \right)
\end{multline}
%\end{linenomath}
for all $v\in\mathcal{V}_{2}$. By hypothesis ${{D}}[., .]$ is a locally Lipschitz continuous function of its probability distribution arguments
(where those distributions are considered as vectors in a Euclidean metric space) when evaluated
for the distributions specified in ~\cref{eq:helper1,eq:helper2,eq:helper3}.
Then since a divergence equals zero only if its arguments are identical, for all $1 \le i \le n$,
for small enough $\epsilon$,
%\begin{linenomath}
\begin{equation}\label{eq:helper4}
    F^n(\mathcal{V}_1 \,|\, \psi(q), E[\beta, B(1), \ldots, B(i)])>F^n(\mathcal{V}_1 \,|\, \psi(q), E[\beta, B(1), \ldots, B(i-1)]),
\end{equation}
%\end{linenomath}
as claimed.
\end{proof}

\subsection{Proof of Prop.~\ref{prop:mathabduction}}
\begin{proof}
By hypothesis, $\PPP$ satisfies the abductive premise,
%\begin{linenomath}
\eq{
\PPP(v^\ast\;|\; q^\ast, (q^\dagger, v^\dagger),\hat{C}_{2}) = \alpha \PPP(v^\ast\;|\; q^\ast, \hat{C}_{2}),
\label{eq:abductionproof1}
}
%\end{linenomath}
which means the abduction implication must hold,
%\begin{linenomath}
\eq{
    \PPP(v^\dagger \;|\; q^\dagger, (q^\ast, v^\ast),\hat{C}_{2}) = \alpha
					\PPP(v^\dagger \;|\; q^\dagger,\hat{C}_{2}).
\label{eq:abductionproof2}
}
%\end{linenomath}
Since hypothesis $
F^{n}_{2}(v^\ast \,|\, q^\ast, q^\dagger, \hat{C}_{2}) :=
\sum_{v^\dagger} F^{n}_{2}(v^\ast, v^\dagger \,|\, q^\ast, q^\dagger, \hat{C}_{2}) =
 F^{n}_{2}(v^\ast \,|\, q^\ast, \hat{C}_{2})$, \cref{eq:abductionproof2} implies
%\begin{linenomath}
\eq{
 \PPP(v^\ast, v^\dagger \;|\; q^\ast, q^\dagger, \hat{C}_{2} ) = \alpha
					\PPP(v^\dagger \;|\; q^\dagger, \hat{C}_{2})  \PPP(v^\ast \;|\; q^\ast, \hat{C}_{2}),
\label{eq:abductionproof2.5}
}
%\end{linenomath}
Also by hypothesis, $\varphi_2$ is calibrated with $\varphi_{1}$ at step $n$ for $\hat{C}_{2}$, for each of the three questions $\psi^{-1}(q^{*})$, $\psi^{-1}(q^{\dagger})$, and $\psi^{-1}(q^{*},q^{\dagger})$, and for $m=1,2$. Then choosing $m=2$, we get
%\begin{linenomath}
\eq{
		\sum_{v\in\mathcal{V}_2}P^{n}_{2}(v|\psi^{-1}(q^{*},q^{\dagger}),\hat{C}_{2})D\left[ \Psi(\psi(\psi^{-1}(q^{*},q^{\dagger})),v)(\mathcal{V}_{1},\mathcal{V}_{1}),  \;
		 \overline{P}_1\left( \mathcal{V}_{1}, \mathcal{V}_{1} \,|\, q^\ast, q^\dagger, \hat{C}_{2}\right) \right]
  				              \le \epsilon.
\label{eq:abductionproof4}
}
%\end{linenomath}
The convexity of $D[.,.]$ establishes
%\begin{linenomath}
\eq{
		D\left[\sum_{v\in\mathcal{V}_2}P^{n}_{2}(v|\psi^{-1}(q^{*},q^{\dagger}),\hat{C})\Psi(\psi(\psi^{-1}(q^{*},q^{\dagger})),v)(\mathcal{V}_{1},\mathcal{V}_{1}),  \;
		 \overline{P}_1\left( \mathcal{V}_{1}, \mathcal{V}_{1} \,|\, q^\ast, q^\dagger, \hat{C}_{2}\right) \right]
  				              \le \epsilon,
\label{eq:abductionproof4.1}
}
%\end{linenomath}
or
%\begin{linenomath}
\eq{
		D\left[\PPP(\mathcal{V}_{1},\mathcal{V}_{1} \,|\, q^\ast, q^\dagger,\hat{C}_{2}),  \;
		 \overline{P}_1\left( \mathcal{V}_{1}, \mathcal{V}_{1} \,|\, q^\ast, q^\dagger, \hat{C}_{2}\right) \right]
  				              \le \epsilon
\label{eq:abductionproof5}
}
%\end{linenomath}
where with abuse of notation we write $\overline{P}_1\left( \mathcal{V}_{1}, \mathcal{V}_{1} \,|\, q^\ast, q^\dagger, \hat{C}_{2}\right)$ for a fixed value of the triple $(q^\ast, q^\dagger, \hat{C}_{2})$ for the associated distribution over an event space $\mathcal{V}_{1} \times \mathcal{V}_{1}$.\par

Go through the analogous reasoning for $m=1$ twice, once for each of the two distinct questions $q' \ne q, q'' \ne q$ that we assume exist, questions which (via $\psi(.)$) specify the single question $q^\ast$ and the single question $q^\dagger$, respectively. In these two cases calibration means that:
%\begin{linenomath}
\eq{
\label{eq:abductionproof6}
		D\left[	\PPP(\mathcal{V}_{1} \,|\, q^\ast,\hat{C}_{2}),  \;  \overline{P}_1 \left(\mathcal{V}_{1} \,|\, q^\ast, \hat{C}_{2}\right) \right]
  				              &\le \epsilon}
\eq{
\label{eq:abductionproof6.5}
		D\left[	\PPP(\mathcal{V}_{1} \,|\, q^\dagger,\hat{C}_{2}),  \;
				\overline{P}_1 \left(\mathcal{V}_{1} \,|\, q^\dagger, \hat{C}_{2}\right) \right]
  				              &\le \epsilon.
}
%\end{linenomath}
(where we have extended the definition of $\PPP$ in the obvious way to the case where it has one claim as an argument rather than two).\par

By hypothesis, ${{D}}[., .]$ is a locally Lipschitz continuous function of its probability distribution arguments (where those distributions are considered as vectors in a Euclidean metric space) when evaluated for the distributions specified in ~\cref{eq:abductionproof5} to \cref{eq:abductionproof6.5}. Then since a divergence equals zero only if its arguments are identical, for small $\epsilon$ \cref{eq:abductionproof5} implies:
%\begin{linenomath}
\eq{
 \overline{P}_1\left( v^\ast, v^\dagger \,|\, q^\ast, q^\dagger, \hat{C}_{2}\right) &\simeq \PPP(v^\ast, v^\dagger \,|\, q^\ast, q^\dagger,\hat{C}_{2}) = \alpha
					\PPP(v^\dagger \;|\; q^\dagger,\hat{C}_{2})  P(v^\ast \;|\; q^\ast, \hat{C}_{2}).
\label{eq:abductionproof6.6}
}
%\end{linenomath}
\cref{eq:abductionproof6} implies:
%\begin{linenomath}
\eq{
 \overline{P}_1\left( v^\ast \,|\, q^\ast,\hat{C}_{2}\right) &\simeq \PPP(v^\ast \,|\, q^\ast,\hat{C}_{2}),
\label{eq:abductionproof6.7}
}
%\end{linenomath}
and \cref{eq:abductionproof6.5} implies:
%\begin{linenomath}
\eq{
 \overline{P}_1\left( v^\dagger \,|\, q^\dagger,\hat{C}_{2}\right) &\simeq \PPP(v^\dagger \,|\, q^\dagger, \hat{C}_{2}).
\label{eq:abductionproof6.8}
}
%\end{linenomath}
Together, \cref{eq:abductionproof6.6} through \cref{eq:abductionproof6.8} imply:
%\begin{linenomath}
\eq{
 \overline{P}_1\left( v^\ast, v^\dagger \,|\, q^\ast, q^\dagger, \hat{C}_{2}\right) &\simeq \alpha
				\overline{P}_1 \left(v^\ast \,|\, q^\ast, \hat{C}_{2}\right)
				\overline{P}_1 \left(v^\dagger \,|\, q^\dagger, \hat{C}_{2}\right)
\label{eq:abductionproof7}
}
%\end{linenomath}
By hypothesis,
%\begin{linenomath}
\eq{
\overline{P}_1 \left(v^\ast \,|\, q^\ast, q^\dagger, \hat{C}_{2}\right) &= \overline{P}_1 \left(v^\ast \,|\, q^\ast, \hat{C}_{2}\right)
}
%\end{linenomath}
and so \cref{eq:abductionproof7} implies that
%\begin{linenomath}
\begin{equation}
    \overline{P}_{1}(v^{\dagger}|q^{\dagger},(q^{*},v^{*}), \hat{C}_{2})\simeq\alpha\overline{P}_{1}(v^{\dagger}|q^{\dagger}, \hat{C}_{2})>\overline{P}_{1}(v^{\dagger}|q^{\dagger}, \hat{C}_{2}),
\end{equation}
%\end{linenomath}
as claimed.
\end{proof}

\subsection{Proof of Prop.~\ref{prop:scienceabduction2}}
\begin{proof}
By hypothesis, $\PPP$ satisfies the abductive premise. Following the same steps as at the beginning of the proof of \cref{prop:mathabduction}, we get
%\begin{linenomath}
\eq{
 \PPP(v^\ast, v^\dagger \;|\; q^\ast, q^\dagger, \hat{C}_{2} ) = \alpha
					\PPP(v^\dagger \;|\; q^\dagger, \hat{C}_{2})  \PPP(v^\ast \;|\; q^\ast, \hat{C}_{2}),
\label{eq:abductionproof2.5s'}
}
%\end{linenomath}
By hypothesis, $\varphi_{2}$ is embed-calibrated with $\varphi_{1}$ for $\hat{C}_{2}$ for each of the three questions $\psi^{-1}(q^\ast)$, $\psi^{-1}(q^\dagger)$, and $\psi^{-1}(q^\ast,q^\dagger)$, and for $m=1,2$. Then choosing $m=2$, we get
%\begin{linenomath}
\begin{multline}\label{eq:abductionproof4s'}
		\sum_{v\in\mathcal{V}_2}P^{n}_{2}(v|\psi^{-1}(q^{*},q^{\dagger}),\hat{C}_{2}])D\Big[ \Psi(\psi(\psi^{-1}(q^{*},q^{\dagger})),v)(\mathcal{V}_{1},\mathcal{V}_{1}),  \\
		 \overline{P}_1\left( \mathcal{V}_{1}, \mathcal{V}_{1} \,|\, q^\ast, q^\dagger, E^{-1}[\{(\psi^{-1}(q^{*},q^{\dagger}),v)\}\cup\hat{C}_{2}]]\right) \Big]\le \epsilon.
\end{multline}
%\end{linenomath}
The convexity of $D[.,.]$ establishes
%\begin{linenomath}
\begin{multline}\label{eq:abductionproof5s'}
		D\Big[\PPP(\mathcal{V}_{1}|\psi^{-1}(q^{*},q^{\dagger}),\hat{C}_{2}]),  \\
		 \sum_{v\in\mathcal{V}_2}P^{n}_{2}(v|\psi^{-1}(q^{*},q^{\dagger}),\hat{C}_{2}])\overline{P}_1\left( \mathcal{V}_{1}, \mathcal{V}_{1} \,|\, q^\ast, q^\dagger, E^{-1}[\{(\psi^{-1}(q^{*},q^{\dagger}),v)\}\cup\hat{C}_{2}]]\right) \Big]\le \epsilon.
\end{multline}
%\end{linenomath}
Go through the analogous reasoning for $m=1$ twice, once for each of the two distinct questions $q' \ne q, q'' \ne q$ that we assume exist, questions which (via $\psi(.)$) specify the single question $q^\ast$ and the single question $q^\dagger$, respectively. In these two cases calibration means that:
%\begin{linenomath}
\eq{
\label{eq:abductionproof6s'}
		D\left[	\PPP(\mathcal{V}_{1} \,|\, q^\ast,\hat{C}_{2}),  \;  \sum_{v\in\mathcal{V}_2}P^{n}_{2}(v|\psi^{-1}(q^{*}),\hat{C}_{2})\overline{P}_1 \left(\mathcal{V}_{1} \,|\, q^\ast, E^{-1}[\{(\psi^{-1}(q^\ast),v)\}\cup \hat{C}_{2}]\right) \right]
  				              &\le \epsilon,}
\eq{
\label{eq:abductionproof6.5s'}
		D\left[	\PPP(\mathcal{V}_{1} \,|\, q^\dagger,\hat{C}_{2}),  \;
				\sum_{v\in\mathcal{V}_2}P^{n}_{2}(v|\psi^{-1}(q^{\dagger}),\hat{C}_{2})\overline{P}_1 \left(\mathcal{V}_{1} \,|\, q^\dagger, E^{-1}[\{(\psi^{-1}(q^\dagger),v)\}\cup \hat{C}_{2}]\right) \right]
  				              &\le \epsilon.
}
%\end{linenomath}
(where we have extended the definition of $\PPP$ in the obvious way to the case where it has one claim as an argument rather than two).\par

By hypothesis, ${{D}}[., .]$ is a locally Lipschitz continuous function of its probability distribution arguments (where those distributions are considered as vectors in a Euclidean metric space) when evaluated for the distributions specified in ~\cref{eq:abductionproof5s'} to \cref{eq:abductionproof6.5s'}. Then since a divergence equals zero only if its arguments are identical, for small $\epsilon$ \cref{eq:abductionproof5s'} implies:
%\begin{linenomath}
\begin{multline}\label{eq:abductionproof6.6s'}
    \PPP(v^\ast,v^\dagger|\psi^{-1}(q^{*},q^{\dagger}),\hat{C}_{2}) = \PPP(v^\dagger \;|\; q^\dagger, \hat{C}_{2})  \PPP(v^\ast \;|\; q^\ast, \hat{C}_{2}) \\ \simeq \ \sum_{v\in\mathcal{V}_2}P^{n}_{2}(v|\psi^{-1}(q^{*},q^{\dagger}),\hat{C}_{2})\overline{P}_1\left(v^\ast,v^\dagger \,|\, q^\ast, q^\dagger, E^{-1}[\{(\psi^{-1}(q^{*},q^{\dagger}),v)\}\cup \hat{C}_{2}]\right).
\end{multline}
%\end{linenomath}
\cref{eq:abductionproof6s'} implies:
%\begin{linenomath}
\eq{
 \PPP(v^\ast \,|\, q^\ast,\hat{C}_{2}) &\simeq \sum_{v\in\mathcal{V}_2}P^{n}_{2}(v|\psi^{-1}(q^{*}),\hat{C}_{2})\overline{P}_1 \left(v^\ast \,|\, q^\ast, E^{-1}[\{(\psi^{-1}(q^\ast),v)\}\cup \hat{C}_{2}]\right),
\label{eq:abductionproof6.7s'}
}
%\end{linenomath}
and \cref{eq:abductionproof6.5s'} implies:
%\begin{linenomath}
\eq{
 \PPP(v^\dagger \,|\, q^\dagger,\hat{C}_{2}) &\simeq \sum_{v\in\mathcal{V}_2}P^{n}_{2}(v|\psi^{-1}(q^{\dagger}),\hat{C}_{2})\overline{P}_1 \left(v^\dagger \,|\, q^\dagger, E^{-1}[\{(\psi^{-1}(q^\dagger),v)\}\cup \hat{C}_{2}]\right).
\label{eq:abductionproof6.8s'}
}
%\end{linenomath}
Together, \cref{eq:abductionproof6.6s'} through \cref{eq:abductionproof6.8s'} imply:
\begin{multline}\label{eq:abductionproof7s'}
    \sum_{v\in\mathcal{V}_2}P^{n}_{2}(v|\psi^{-1}(q^{*},q^{\dagger}),\hat{C}_{2})\overline{P}_1\left(v^\ast,v^\dagger \,|\, q^\ast, q^\dagger, E^{-1}[\{(\psi^{-1}(q^{*},q^{\dagger}),v)\}\cup \hat{C}_{2}]\right) \\ \simeq \ \alpha \sum_{v\in\mathcal{V}_2}P^{n}_{2}(v|\psi^{-1}(q^{*}),\hat{C}_{2})\overline{P}_1 \left(v^\ast \,|\, q^\ast, E^{-1}[\{(\psi^{-1}(q^\ast),v)\}\cup \hat{C}_{2}\right) \\ \times \ \sum_{v\in\mathcal{V}_2}P^{n}_{2}(v|\psi^{-1}(q^{\dagger}),\hat{C}_{2})\overline{P}_1 \left(v^\dagger \,|\, q^\dagger, E^{-1}[\{(\psi^{-1}(q^\dagger),v)\}\cup \hat{C}_{2}]\right).
\end{multline}
By hypothesis, for all $v\in\mathcal{V}_{2}$,
\begin{equation}
    \overline{P}_1 \left(v^\ast \,|\, q^\ast, q^{\dagger}, E^{-1}[\{(\psi^{-1}(q^\ast,q^\dagger),v)\}\cup \hat{C}_{2}]\right) = \overline{P}_1 \left(v^\ast \,|\, q^\ast, E^{-1}[\{(\psi^{-1}(q^\ast),v)\}\cup \hat{C}_{2}]\right)
\end{equation}
and so \cref{eq:abductionproof7s'} implies that
%\begin{linenomath}
\begin{multline}
    \sum_{v\in\mathcal{V}_{2}}P^{n}_{2}(v|\psi^{-1}(q^{\dagger}),(q^{*},v^{*}),\hat{C}_{2})\overline{P}_{1}(v^{\dagger}|q^{\dagger},(q^{*},v^{*}), E^{-1}[\{(\psi^{-1}(q^{\dagger}),v)]\}\cup\hat{C}_{2}]) \\ > \ \sum_{v\in\mathcal{V}_{2}}P^{n}_{2}(v|\psi^{-1}(q^{\dagger}),\hat{C}_{2})\overline{P}_{1}(v^{\dagger}|q^{\dagger}, E^{-1}[\{\psi^{-1}(q^{\dagger}),v)\}\cup\hat{C}_{2}])
\end{multline}
%\end{linenomath}
as claimed.
\end{proof}

\subsection{Proof of Prop.~\ref{prop:scienceabduction}}
\begin{proof}
By hypothesis, $\overline{P}_{1}$ satisfies the abductive premise for any $\hat{C}_{1}\in\mathcal{\hat{C}}_{1}$, meaning that for any $v\in\mathcal{V}_{2}$,
%\begin{linenomath}
\eq{
\overline{P}_1\left(v^\ast\,|\, q^\ast, (q^\dagger,v^{\dagger}), E^{-1}[\{(\psi^{-1}(q^{*},q^{\dagger}),v)\}\cup E[\hat{C}_{1}]]\right)  = \alpha \overline{P}_1 \left(v^\ast \,|\, q^\ast, E^{-1}[\{(\psi^{-1}(q^\ast),v)\}\cup E[\hat{C}_{1}]]\right),
\label{eq:abductionproof1s}
}
%\end{linenomath}
which allows us to derive
%\begin{linenomath}
\eq{
    \overline{P}_1\left(v^\dagger\,|\, q^\dagger, (q^\ast,v^{\ast}), E^{-1}[\{(\psi^{-1}(q^{*},q^{\dagger}),v)\}\cup E[\hat{C}_{1}]]\right)  = \alpha \overline{P}_1 \left(v^\dagger \,|\, q^\dagger, E^{-1}[\{(\psi^{-1}(q^\dagger),v)\}\cup E[\hat{C}_{1}]]\right).
\label{eq:abductionproof2s}
}
%\end{linenomath}
Since, also by hypothesis,
\begin{equation}
    \overline{P}_{1}\left(v^\ast \,|\, q^\ast, q^\dagger, E^{-1}[\{(\psi^{-1}(q^{*},q^{\dagger}),v)\}\cup E[\hat{C}_{1}]]\right) = \overline{P}_{1} \left(v^\ast \,|\, q^\ast, E^{-1}[\{(\psi^{-1}(q^\dagger),v)\}\cup E[\hat{C}_{1}]]\right),
\end{equation}
we are able to rewrite \cref{eq:abductionproof2s} as
\begin{multline}\label{eq:abductionproof2.5s}
    \overline{P}_1\left(v^\dagger, v^\ast\,|\, q^\dagger, q^\ast, E^{-1}[\{(\psi^{-1}(q^{*},q^{\dagger}),v)\}\cup E[\hat{C}_{1}]]\right)  \\ = \alpha
					\overline{P}_1 \left(v^\dagger \,|\, q^\dagger, E^{-1}[\{(\psi^{-1}(q^\dagger),v)\}\cup E[\hat{C}_{1}]\right) \overline{P}_1 \left(v^\ast \,|\, q^\ast, E^{-1}[\{(\psi^{-1}(q^\ast),v)\}\cup E[\hat{C}_{1}]]\right).
\end{multline}
%\begin{linenomath}
%\end{linenomath}
By hypothesis, $\varphi_{2}$ is embed-calibrated with $\varphi_{1}$ for $E[\hat{C}_{1}]$ for each of the three questions $\psi^{-1}(q^\ast)$, $\psi^{-1}(q^\dagger)$, and $\psi^{-1}(q^\ast,q^\dagger)$, and for $m=1,2$. Then choosing $m=2$ and applying \cref{eq:abductionproof2.5s}, we get
%\begin{linenomath}
\begin{multline}\label{eq:abductionproof4s}
		\sum_{v\in\mathcal{V}_2}P^{n}_{2}(v|\psi^{-1}(q^{*},q^{\dagger}),E[\hat{C}_{1}])D\Big[ \Psi(\psi(\psi^{-1}(q^{*},q^{\dagger})),v)(\mathcal{V}_{1},\mathcal{V}_{1}),  \\
		 \overline{P}_1\left( \mathcal{V}_{1}, \mathcal{V}_{1} \,|\, q^\ast, q^\dagger, E^{-1}[\{(\psi^{-1}(q^{*},q^{\dagger}),v)\}\cup E[\hat{C}_{1}]]\right) \Big]\le \epsilon.
\end{multline}
%\end{linenomath}
The convexity of $D[.,.]$ establishes
%\begin{linenomath}
\begin{multline}\label{eq:abductionproof5s}
    D\Big[\PPP(\mathcal{V}_{1},\mathcal{V}_{1}|\psi^{-1}(q^{*},q^{\dagger}),E[\hat{C}_{1}]),  \\
		 \sum_{v\in\mathcal{V}_2}P^{n}_{2}(v|\psi^{-1}(q^{*},q^{\dagger}),E[\hat{C}_{1}])\overline{P}_1\left(\mathcal{V}_{1}, \mathcal{V}_{1} \,|\, q^\ast, q^\dagger, E^{-1}[\{(\psi^{-1}(q^{*},q^{\dagger}),v)\}\cup E[\hat{C}_{1}]]\right) \Big]
  				              \le \epsilon.
\end{multline}
%\end{linenomath}
Go through the analogous reasoning for $m=1$ twice, once for each of the two distinct questions $q' \ne q, q'' \ne q$ that we assume exist, questions which (via $\psi(.)$) specify the single question $q^\ast$ and the single question $q^\dagger$, respectively. In these two cases calibration means that:
%\begin{linenomath}
\eq{
\label{eq:abductionproof6s}
		D\left[	\PPP(\mathcal{V}_{1} \,|\, q^\ast,E[\hat{C}_{1}]),  \;  \sum_{v\in\mathcal{V}_2}P^{n}_{2}(v|\psi^{-1}(q^{*}),E[\hat{C}_{1}])\overline{P}_1 \left(\mathcal{V}_{1} \,|\, q^\ast, E^{-1}[\{(\psi^{-1}(q^\ast),v)\}\cup E[\hat{C}_{1}]]\right) \right]
  				              &\le \epsilon,}
\eq{
\label{eq:abductionproof6.5s}
		D\left[	\PPP(\mathcal{V}_{1} \,|\, q^\dagger,E[\hat{C}_{1}]),  \;
				\sum_{v\in\mathcal{V}_2}P^{n}_{2}(v|\psi^{-1}(q^{\dagger}),E[\hat{C}_{1}])\overline{P}_1 \left(\mathcal{V}_{1} \,|\, q^\dagger, E^{-1}[\{(\psi^{-1}(q^\dagger),v)\}\cup E[\hat{C}_{1}]]\right) \right]
  				              &\le \epsilon.
}
%\end{linenomath}
(where we have extended the definition of $\PPP$ in the obvious way to the case where it has one claim as an argument rather than two).\par

By hypothesis, ${{D}}[., .]$ is a locally Lipschitz continuous function of its probability distribution arguments (where those distributions are considered as vectors in a Euclidean metric space) when evaluated for the distributions specified in ~\cref{eq:abductionproof5s} to \cref{eq:abductionproof6.5s}. Then since a divergence equals zero only if its arguments are identical, for small $\epsilon$ \cref{eq:abductionproof5s} implies:
%\begin{linenomath}
\begin{multline}\label{eq:abductionproof6.6s}
    \PPP(v^\ast,v^\dagger|\psi^{-1}(q^{*},q^{\dagger}),E[\hat{C}_{1}]) = \PPP(v^\dagger \;|\; q^\dagger, E[\hat{C}_{1}])  \PPP(v^\ast \;|\; q^\ast, E[\hat{C}_{1}]) \\ \simeq \ \sum_{v\in\mathcal{V}_2}P^{n}_{2}(v|\psi^{-1}(q^{*},q^{\dagger}),E[\hat{C}_{1}])\overline{P}_1\left(v^\ast,v^\dagger \,|\, q^\ast, q^\dagger, E^{-1}[\{(\psi^{-1}(q^{*},q^{\dagger}),v)\}\cup E[\hat{C}_{1}]]\right).
\end{multline}
%\end{linenomath}
\cref{eq:abductionproof6s} implies:
%\begin{linenomath}
\eq{
 \PPP(v^\ast \,|\, q^\ast,E[\hat{C}_{1}]) &\simeq \sum_{v\in\mathcal{V}_2}P^{n}_{2}(v|\psi^{-1}(q^{*}),E[\hat{C}_{1}])\overline{P}_1 \left(v^\ast \,|\, q^\ast, E^{-1}[\{(\psi^{-1}(q^\ast),v)\}\cup E[\hat{C}_{1}]]\right),
\label{eq:abductionproof6.7s}
}
%\end{linenomath}
and \cref{eq:abductionproof6.5s} implies:
%\begin{linenomath}
\eq{
 \PPP(v^\dagger \,|\, q^\dagger,E[\hat{C}_{1}]) &\simeq \sum_{v\in\mathcal{V}_2}P^{n}_{2}(v|\psi^{-1}(q^{\dagger}),E[\hat{C}_{1}])\overline{P}_1 \left(v^\dagger \,|\, q^\dagger, E^{-1}[\{(\psi^{-1}(q^\dagger),v)\}\cup E[\hat{C}_{1}]]\right).
\label{eq:abductionproof6.8s}
}
%\end{linenomath}
Together, \cref{eq:abductionproof6.6s} through \cref{eq:abductionproof6.8s} imply:
%\begin{linenomath}
\eq{
 \PPP(v^\ast,v^\dagger|\psi^{-1}(q^{*},q^{\dagger}),E[\hat{C}_{1}]) &\simeq \alpha
				\PPP(v^\ast \,|\, q^\ast,E[\hat{C}_{1}])
				\PPP(v^\dagger \,|\, q^\dagger,E[\hat{C}_{1}]).
\label{eq:abductionproof7s}
}
%\end{linenomath}
By hypothesis,
%\begin{linenomath}
\eq{
\PPP \left(v^\ast \,|\, q^\ast, q^\dagger, E[\hat{C}_{1}]\right) &= \PPP \left(v^\ast \,|\, q^\ast, E[\hat{C}_{1}]\right)
}
%\end{linenomath}
and so \cref{eq:abductionproof7s} implies that
%\begin{linenomath}
\begin{equation}
    \PPP(v^{\dagger}|q^{\dagger},(q^{*},v^{*}),E[\hat{C}_{1}])>\PPP(v^{\dagger}|q^{\dagger}, E[\hat{C}_{1}]),
\end{equation}
%\end{linenomath}
as claimed.\par

\end{proof}

\end{appendix}

\printbibliography

%\bibliography{sms}

\end{document}